\documentclass[12pt]{article}
\usepackage{fullpage,graphicx,psfrag,amsmath,amsfonts,verbatim,cleveref}
\usepackage[small,bf]{caption}
\usepackage{subcaption}
\usepackage{url}
\usepackage{booktabs}

\usepackage{algorithm}
\usepackage{algpseudocode}
\usepackage{algorithmicx}

\usepackage[suppress]{color-edits}
\addauthor{authors}{blue}
\addauthor{df}{blue}
\addauthor{ds}{blue}
\addauthor{sh}{blue}

\newcommand{\OneAndAHalfSpacedXI}{}
\newcommand{\DoubleSpacedXI}{}
\newcommand{\SingleSpacedXI}{}

\usepackage{amsmath,amsthm}
\newcommand{\Halmos}{{\qed}}
\newcommand{\halmos}{{\qed}}

\newtheorem{proposition}{Proposition}

\newtheorem{theorem}{Theorem}
\newtheorem{definition}{Definition}

\newtheorem{lemma}{Lemma}[section]

\newtheorem{example}{Example}

\usepackage{natbib}

\title{Minimizing Multimodular Functions and Allocating Capacity in Bike-Sharing Systems}
\author{Daniel Freund \and Shane G. Henderson \and David B. Shmoys}

\begin{document}
\maketitle

\begin{abstract}
The growing popularity of bike-sharing systems around the world has motivated recent attention to models and algorithms for their effective operation. Most of this literature focuses on their daily operation for managing asymmetric demand. In this work, we consider the more strategic question of how to (re-)allocate dock-capacity in such systems. We develop mathematical formulations for variations of this problem (either for service performance over the course of one day or for a long-run-average) and exhibit discrete convex properties in associated optimization problems. This allows us to design a practically fast polynomial-time allocation algorithm to compute an optimal solution for this problem, which can also handle practically motivated constraints, such as a limit on the number of docks moved in the system.

We apply our algorithm to data sets from Boston, New York City, and Chicago to investigate how different dock allocations can yield better service in these systems. Recommendations based on our analysis have led to changes in the system design in Chicago and New York City. Beyond optimizing for improved quality of service through better allocations, our results also provide a metric to compare the impact of strategically reallocating docks and the daily rebalancing of bikes.
\end{abstract}

\newpage
\tableofcontents
\newpage

\section{Introduction}

As bike-sharing systems become an integral part of the urban landscape, novel lines of research seek to model and optimize their operations. In many systems, such as New York City's Citi Bike, users can rent and return bikes at any station within the city. This flexibility makes the system attractive for commuters and tourists alike. From an operational point of view, however, this flexibility leads to imbalances when demand is asymmetric, as is commonly the case. The main contributions of this paper are to identify key questions in the {\em design} of operationally efficient bike-sharing systems, to develop a polynomial-time algorithm for the associated discrete optimization problems, to apply this algorithm on real usage data, and to investigate the effect this optimization has in practice.

The largest bike-sharing systems in the US are dock-based, meaning that they consist of stations, spread across a city, each of which has a number of docks in which bikes can be locked. If a bike is present in a dock, users can rent it and return it at any other station with an open dock. However, system imbalance often causes some stations to have only empty (or \emph{open}) docks and others to have only full docks (i.e., ones filled with bikes). In the former case, users need to find alternate modes of transportation, whereas in the latter they might not be able to end their trip at the intended destination. In many bike-sharing systems, this has been found to be a leading cause of customer dissatisfaction, e.g., \cite{cabisurvey2014}.

In order to meet demand in the face of asymmetric traffic, bike-sharing system operators seek to \emph{rebalance} the system by moving bikes from locations with too few open docks to locations with too few bikes. To facilitate these operations, a burst of recent research has investigated models and algorithms to increase their efficiency and increase customer satisfaction. While similar in spirit to some of the literature on rebalancing, in this work we use a different control to increase customer satisfaction. Specifically, we answer the question \emph{how should bike-sharing systems allocate dock capacity to stations within the system so as to minimize the number of dissatisfied customers}? 

\begin{figure}[ht]
    \centering
    \includegraphics[width=.5\textwidth]{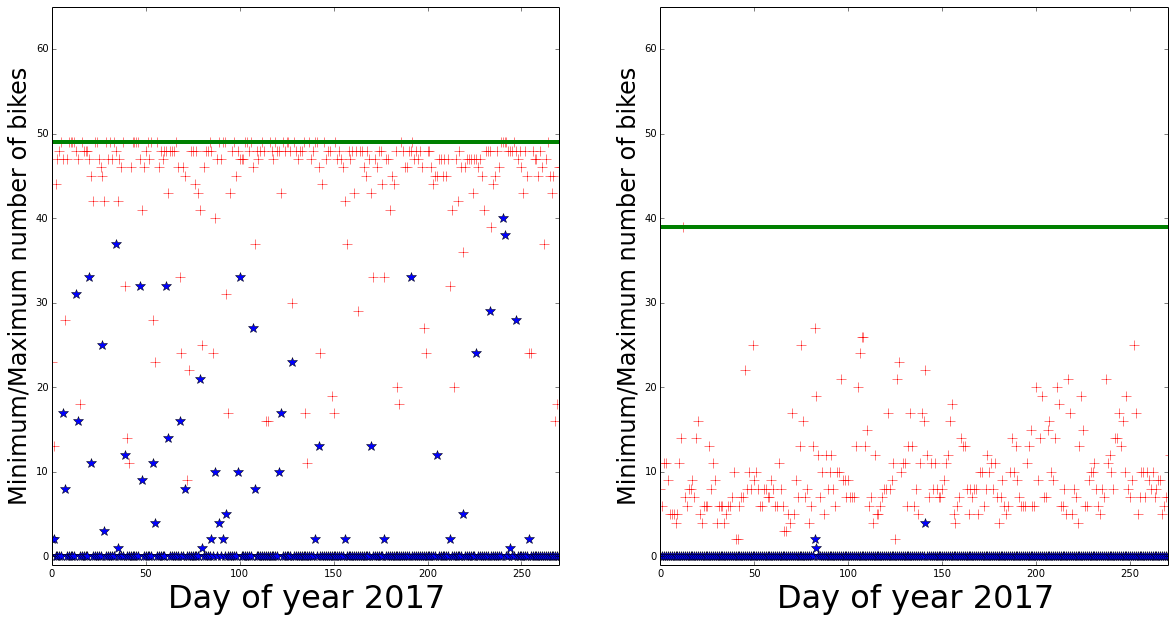}
    \caption{Minimum (blue) and maximum (red) number of bikes at two bike-sharing stations over the course of each day of January-September 2017. The green lines denote the capacities of the stations.}
    \label{fig:fill_levels}
\end{figure}
A superficial analysis of usage data reveals that there may be potential in reallocating capacity: some stations have spare capacity that users never or rarely use (see Figure \ref{fig:fill_levels}) whereas other stations have all of their capacity used on most days.
We give a more theoretically grounded answer to this question by developing two optimization models, both based on the underlying metric that system performance is captured by the expected number of customers that do not receive service. In the first model, we focus on planning one day, say 6am-midnight, where for each station we determine its allocation of bikes and docks; this framework assumes that there is sufficient rebalancing capacity overnight to restore the desired bike allocation by 6am the next morning. Since in practice this turns out to be quite difficult, the second model considers a set-up induced by a long-run average which assumes that no rebalancing happens at all; in a sense, this exhibits the opposite regime. The theory developed in this paper enabled extensive computational experiments on real data sets; through these we found that there are dock allocations that simultaneously perform well with respect to both models, yielding improvements to both (in comparison to the current allocation) of up to 20\%. {These results were leveraged by system operators in Chicago and New York City 
and led to 100 (200) docks being moved in New York City (Chicago). 
Convinced by the impact analysis in these cities, operators of other major US bike-sharing systems, including Blue Bikes in Boston and Capital Bikeshare in Washington, D.C., have run our analysis on their data to capture the potential of reallocated dock capacity as well.} 




\subsection{Our Contribution}\label{ssec:contribution}


\cite{raviv2013optimal} defined a \emph{user dissatisfaction function} (UDF) that measures the expected number of out-of-stock events at an individual bike-sharing station. To do so, they define a stochastic process on the possible number of bikes (between 0 and the capacity of the station). The stochastic process observes attempted rentals and returns of bikes over time; {this process is assumed to be exogenously given at each station and independent of our decisions/the availability of bikes and docks in other stations.} Each arrival triggers a change in the state, either decreasing (rental) or increasing (return) the number of available bikes by one. When the number of bikes is 0 and a rental is attempted, or when it equals the station capacity and a return is attempted, a customer experiences an out-of-stock event. Various follow-up papers, (\cite{schuijbroek2013inventory}, \cite{o2015smarter}, and \cite{parikh2014estimation}), have suggested different ways to compute the expected number of out-of-stock events $c_i(d_i,b_i)$ that occur over the course of one day at each station $i$ for a given allocation of $b_i$ bikes and $d_i$ empty docks (i.e., $d_i+b_i$ docks in total) at station $i$ at the start of the day.



We use the same UDFs to model the question of how to allocate dock capacity within the system. Given $c_i(\cdot,\cdot)\;\forall i$, our goal is to find an allocation of bikes and docks in the system that minimizes the total expected number of out-of-stock events within a system of $n$ stations, i.e., $\sum_{i=1}^n c_i(d_i,b_i)$. Since the number of bikes and docks is limited, we need to accommodate a \emph{budget constraint} $B$ on the number of bikes in the system and another on the number of docks $D+B$ in the system. Other constraints are often important, such as lower and upper bounds on the capacity for a particular station; furthermore, through our collaboration with Citi Bike in NYC it also became apparent that \emph{operational constraints} limit the number of docks moved from the current system configuration. Thus, we aim to minimize the objective among solutions that require at most some number of docks moved. Notice that $D$ and $B$ could either denote the inventory that is currently present in the system (in which case the question is how to reallocate it) or include new inventory (in which case the question is how to augment the current system design). 

After formally defining this model and discussing its underlying assumptions in Section \ref{sec: model}, we design in Section \ref{sec: algorithm} a discrete gradient-descent algorithm that provably solves the minimization problem with $O(n+B+D)$ oracle calls to evaluate cost functions and an (in practice, vastly dominated) overhead of $O((n+B+D)\log(n))$ elementary list operations. In Section \ref{sec: scaling} we show that scaling techniques, together with a subtle extension of the analysis of the gradient-descent algorithm, improve the running-time to $O(n\log(B+D))$ oracle calls and $O(\log(B+D)(n\log(n)))$ elementary list operations for the setting without operational constraints; in Appendix \ref{sec:appendix_proof_scaling} we include the proofs thereof as well as explanations of how operational constraints can be handled when aiming for running-time logarithmic in $B+D$. In Appendix \ref{sec: running_time}, we include a computational study to complement this theoretical analysis of the efficiency of our algorithms. 

The primary motivation of this analysis is to investigate whether the number of out-of-stock events in bike-sharing systems can be significantly reduced by a data-driven approach. In Section~\ref{sec: case_study}, we apply the algorithms to data sets from Boston, NYC, and Chicago to evaluate the impact on out-of-stock events. One shortcoming of that optimization problem is its assumption that we can perfectly restore the system to the desired initial bike allocation overnight. Through our collaboration with the operators of systems across the country, it has become evident that current rebalancing efforts overnight are vastly insufficient to realize such an optimal (or even near-optimal) allocation of bikes for the current allocation of docks. Thus, we consider in Section \ref{sec: long_run} the opposite regime, in which no rebalancing occurs at all. To model this, we define an extension of the cost function under a long-run average regime. In this regime, the assumed allocation of bikes at each station is a function of only the number of docks and the estimated demand at that station. Interestingly, our empirical results reveal that operators of bike-sharing systems can \emph{have their cake and eat it too}: optimizing dock allocations for one of the objectives (optimally rebalanced or long-run average) yields most of the obtainable improvement for the other. 

{Based on our recommendations the operators of Citi Bike in New York City agreed with the city's Department of Transportation to move 34 docks between 6 stations as part of a pilot program. We use these moves to evaluate the impact of reallocated capacity. Specifically, in Section \ref{sec: impact}, we prove that observing rentals and returns after capacity has been added provides a natural way to estimate the reduction in out-of-stock events (due to dock capacity added) that can be computed in a very simple manner. We apply this approach to the stations that were part of the pilot to derive estimates for the realized reduction in the number of stockouts at those stations.}  





\subsection{Related Work}
\label{sec: rel_work}

A recent line of work, including variations by \cite{raviv2013static}, \cite{forma20153}, \cite{kaspi2015bike}, \cite{ho2014solving}, and  \cite{freundSmarter}, considered static rebalancing problems, in which a capacitated truck (or a fleet of trucks) is routed over a limited time horizon. The truck may pick up and drop off bikes at each station, so as to minimize the expected number of out-of-stock events that occur after the completion of the route. These are evaluated by the same objective function of \cite{raviv2013optimal} that we consider as well.


In contrast to this line of work, \cite{o2015smarter} addressed the question of allocating both docks and bikes; he uses the UDFs (defined over a single interval with constant rental and return rates) to design a mixed integer program over the possible allocations of bikes and docks. Our work extends upon this by providing a fast algorithm for generalizations of that same problem and extensions thereof. 
The optimal allocation of bikes has also been studied by \cite{jianhenderson}, \cite{datner2015setting}, and by  \cite{jianetal}, with the latter also considering the allocation of docks (in fact, the idea behind the algorithm considered by \cite{jianetal} is based on an early draft of this paper). They each develop frameworks based on ideas from simulation optimization; while they also treat demand for bikes as being exogenous, their framework captures the downstream effects of changes in supply upstream. \cite{jianetal} found that these effects are mostly captured by decensoring piecewise-constant demand estimates (see Section \ref{sec: assumptions}).


Orthogonal approaches to the question of where to allocate docks have been taken by \cite{kabra2015bike} and \cite{wang2016applying}. The former considers demand as endogenous and aims to identify the station density that maximizes sales, whereas we consider demand and station locations as exogenously given and aim to allocate docks and bikes to maximize the amount of demand that is being met. The latter aims to use techniques from retail location theory to find locations for stations to be added to an existing system.

Further related literature includes a line of work on rebalancing triggered by  \cite{chemla2013bike}. Susbequent papers, e.g., by \cite{nair2013large},  \cite{dell2014bike}, \cite{erdougan2014static}, \cite{erdougan2015exact}, \cite{bruck2019static}, and \cite{li2020branch} solve variants of a routing problem with fixed numbers of bikes that need to be picked up/dropped off at each station -- \cite{de2016bike} extensively surveys these papers. Before rebalancing bike-sharing systems  became an object of academic study, the closely related traveling salesman problems with pickup and delivery had already been studied outside the bike-sharing domain since \cite{hernandez2004branch}.  Other approaches to rebalancing include for example the papers of \cite{liurebalancing}, \cite{ghosh2016robust}, \cite{rainer2013balancing}, \cite{shu2013models}, or more recently \cite{brinkmann2019dynamic}. We refer the readers to the surveys of
\cite{laporte2018shared},
\cite{freund2019bike}, and
\cite{shui2020review} for a wider overview of the rebalancing literature.   While all of these fall into the wide range of recent work on the operation of bike-sharing systems, they differ from our work in the controls and methodologies they employ.

Finally, a great deal of work has been conducted in the context of predicting demand. In this work, we assume that the predicted demand is given, e.g., using the methods of \cite{o2015data} or \cite{singhvi2015predicting}. Further methods to predict demand have been suggested by \cite{li2015traffic}, \cite{chen2016dynamic}, and \cite{zhang2016bicycle} among others. Our results can be combined with any approach that predicts demand at each station independently of all others.


\noindent\textbf{Relation to Discrete Convexity. } Our algorithms and analyses crucially exploit the property that the UDFs $c_i(\cdot,\cdot)$ at each station are multimodular (see Definition \ref{def: multimod}). 
This provides an interesting connection to the literature on discrete convex analysis. {Prior works connecting inventory management to discrete convexity include \cite{lu2005order,zipkin2008structure,li2014multimodularity} among others; we refer the reader to a recent survey by \cite{chen2021discrete} for an extensive overview.} In concurrent work by \cite{kaspi2015bike} it was shown that the number of out-of-stock events $F(b,U-d-b)$ at a bike-sharing station with fixed capacity $U$, $b$ bikes, and $U-d-b$ unusable bikes is $M^\natural$ (read \emph{M natural}) convex in $b$ and $U-d-b${; functions with such discrete convex properties, in particular $M$-convex and $M^\natural$ functions, were respectively introduced by \cite{murota1996convexity,murota1998discrete} and \cite{murota1999m} (see the book by \cite{murota2003discrete} for a complete overview of early results in discrete convexity).} 
Unusable bikes effectively reduce the capacity at the station, since they are assumed to remain in the station over the entire time horizon. A station with capacity $U$, $b$ bikes, and $U-b-d$ unusable bikes, must then have~$d$ empty docks; hence,  $c(d,b)=F(b,U-d-b)$ for $d+b\leq U$, which parallels our result that $c(\cdot,\cdot)$ is multimodular. Though this would suggest that algorithms to minimize $M^\natural$-convex functions could solve our problem optimally, one can show that $M^\natural$-convexity is not preserved, even in the version with only budget constraints: we provide in Appendix \ref{appendix_M} an example that shows both that  an $M^\natural$-convex function restricted to an $M^\natural$-convex set is not $M^\natural$-convex and that Murota's algorithm for $M^\natural$-convex function minimization can be suboptimal in our setting. In fact, when including the operational constraints even discrete midpoint convexity, a strict generalization of multimodularity studied for example by \cite{fujishige2000notes} and \cite{moriguchi2017discrete}, which is in turn much weaker than $M^\natural$ convexity, breaks down. We provide an example for this in Appendix \ref{appendix_dmcf}.  Surprisingly, we are nevertheless able to design fast algorithms; these exploit not only the multimodularity of each individual $c_i$, but also the separability of the objective function (w.r.t. the stations), that is, the fact that each $c_i$ is only a function of $d_i$ and $b_i$. This not only extends ideas from the realm of unconstrained discrete convex minimization to the constrained setting, but also yields algorithms that (for our special case) have significantly faster running times than those that would usually arise in the context of multimodular function minimization. Since the conference version of this paper appeared, \cite{shioura2018m} has taken our work as motivation to study M-convex function minimization under L1-distance constraints, a strict generalization of our objective. 
Finally, Shioura (private communication) pointed out an error in a preliminary version of this paper, and so, although all of the main elements of our proof of correctness of the discrete gradient-descent algorithm can be found in our preliminary version \citep{freund2017minimizing, freund2016minimizing}, the presentation here differs from that given earlier.



\section{Model}\label{sec: model}
The fundamental primitives of our model of a bike-sharing system are customers, bikes, docks, and stations. Below we formally define these primitives and the optimization problem that is based on them.
\subsubsection*{Model primitives.} A bike-sharing system consists of $n$ \emph{stations}. Each station $i$ is characterized by an exogenously given \emph{demand profile} $p_i$, where $p_i$ is a distribution over arrival sequences of \emph{customers} at~$i$ over the course of a time horizon (e.g., 6AM-12AM). Such arrival sequences are denoted $X=(X_1,X_2,\ldots, X_s)\in\{\pm1\}^s$ where $X_t=-1$ corresponds to a customer arriving to rent a bike, and~$X_t=1$ corresponds to a customer arriving to return a bike. The ability of a customer arriving at a station to rent, resp. return, a bike is dependent on the number of \emph{bikes}, resp. \emph{empty docks}, available at the station at the time of arrival: if no bikes, resp. empty docks, are available at the time of the customer's arrival, the customer is unable to rent, resp. return, a bike and disappears with an out-of-stock event. If instead a customer arrives at a station to rent a bike and a bike is available, then the number of bikes at the station decreases by 1, and the number of empty docks at the station increases by~1. Similarly, if a customer arrives to return a bike, and an empty dock is available at that time, then the number of bikes at the station increases by~1 and the number of empty docks decreases by 1. Notice that, throughout the time horizon, the total number of docks (empty and full) at station $i$ remains the same. This is because the number of empty docks increases by 1 if and only if the number of full docks decreases by 1 (and vice versa).  The respective number of bikes and empty docks at the station at the time of arrival of $X_t$ is based only on (i) the initial allocation of bikes and empty docks at the beginning of the time horizon, and (ii) the arrival sequence of customers up to $X_t$, which we denote $X(t-1)=(X_1,\ldots,X_{t-1})$.
\subsubsection*{Initial allocations of bikes and docks.} The decision variables in our optimization model are the initial number of empty docks and bikes allocated to each station $i$ at the beginning of the time horizon. We denote the initial number of empty docks at station $i$ by $d_i$, and the initial number of bikes (full docks) by $b_i$; combining these two we find that a station $i$ has an allocated capacity of $d_i+b_i$ docks in total.

\OneAndAHalfSpacedXI
\begin{figure}
\centering
\begin{subfigure}{.5\textwidth}
  \centering
  \includegraphics[width=.6\textwidth]{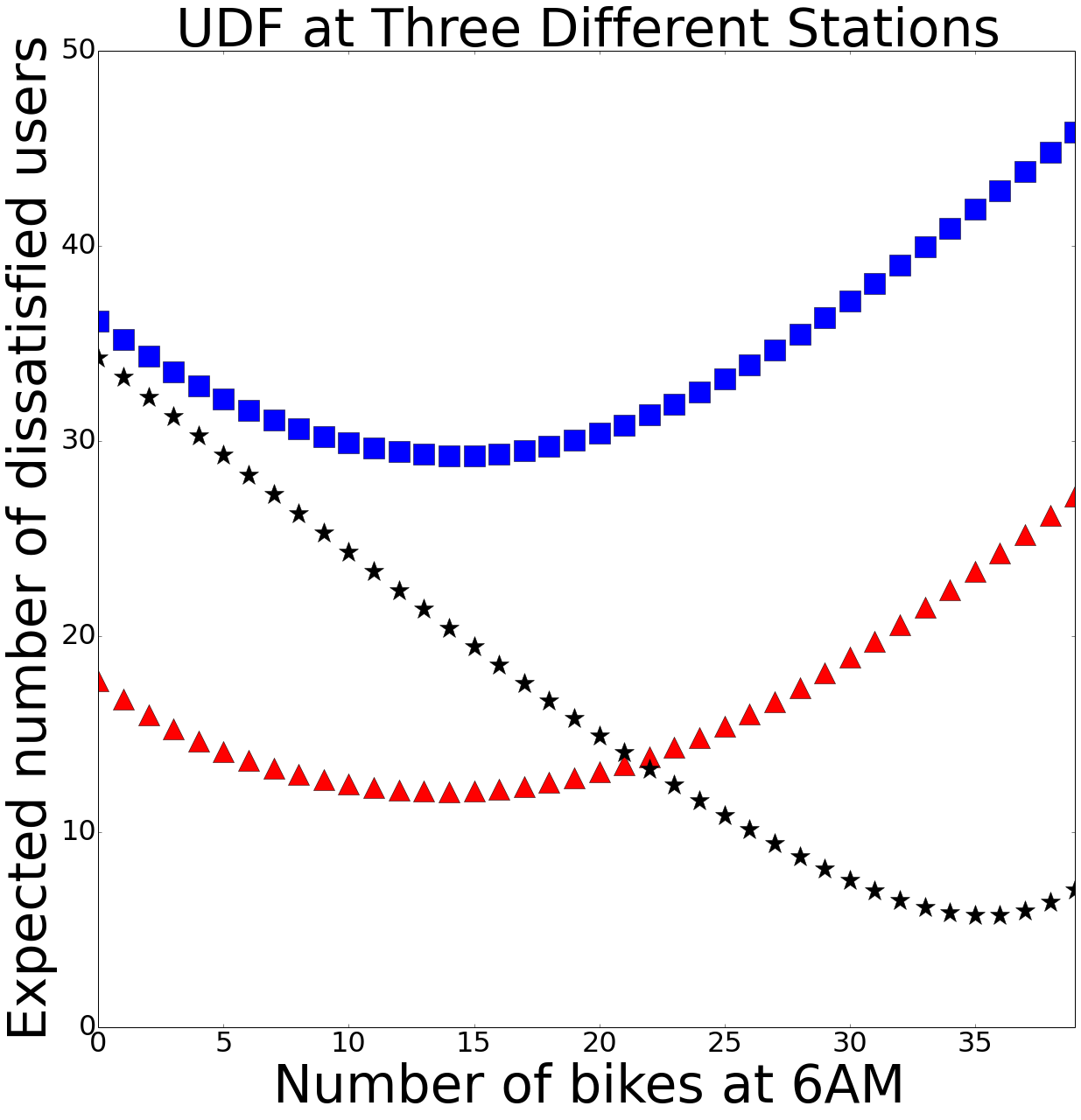}
  \caption{As a function of bikes for stations with capacity~39.}
  \label{fig:4UDFs}
\end{subfigure}%
\begin{subfigure}{.555\textwidth}
  \centering
  \includegraphics[width=.633\textwidth]{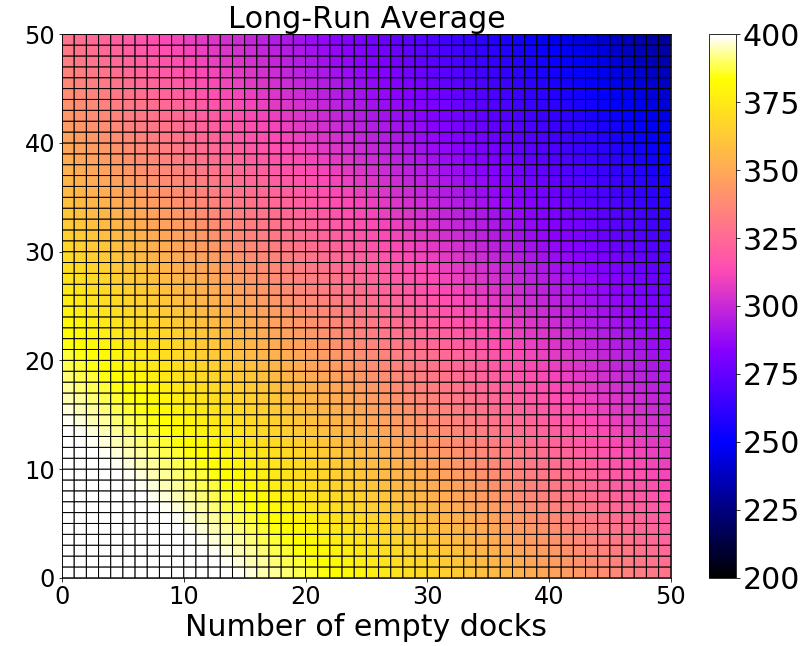}
  \caption{As a function of $(d,b)$ at a single station.}
  \label{fig:sub2}
\end{subfigure}
\caption{Visualizations of user dissatisfaction functions based on real data.}
\label{fig:UDF2D}
\end{figure}
\DoubleSpacedXI

\subsubsection*{User dissatisfaction function.} 
The UDF $c^X(d,b)$ maps the initial number of empty docks and bikes to the number of customers, among the sequence $X=(X_1,\ldots,X_s)$, that experience out-of-stock events (see Figure \ref{fig:UDF2D}). Then, the UDF at station $i$ is given by $c_i(d,b)=\mathbb{E}_{X\sim p_i}[c^X_i(d,b)]$. In an effort to keep notation concise in the main body of the text, we move a formal recursive definition of $c^X(\cdot,\cdot)$ to Appendix \ref{sec:appendix_lemma_multimod}. {UDFs are sometimes used with different weights for stock-outs depending on whether they occur at empty or full stations; while we focus throughout on the unweighted case, in which $c^X(d,b)$ is just a count of the stockouts, our results extend to the weighted case (see Appendix \ref{sec:appendix_lemma_multimod}).}
\begin{definition}
A function $f:\mathbb{Z}^2\to\mathbb{R}\cup\{\infty\}$ with
\OneAndAHalfSpacedXI
\begin{eqnarray}
f(d+1,b+1)-f(d+1,b)\geq f(d,b+1)-f(d,b)&;\\
f(d-1,b+1)-f(d-1,b)\geq f(d,b)-f(d,b-1) &;\\
f(d+1,b-1)-f(d,b-1)\geq f(d,b)-f(d-1,b) &;
\end{eqnarray}
%
\noindent for all $d,b$ is called \emph{multimodular} \citep{hajek1985extremal,altman2000multimodularity,murota2003discrete}\label{def: multimod}.
For future reference, we also define the following implied additional inequalities ((6) \& (1) are equivalent, (1) \& (2) imply (5), and (3) \& (6) imply (4)): 
%
\begin{eqnarray}
f(d+2,b)-f(d+1,b)\geq f(d+1,b)-f(d,b);\\
f(d,b+2)-f(d,b+1)\geq f(d,b+1)-f(d,b);\\
f(d+1,b+1)-f(d,b+1)\geq f(d+1,b)-f(d,b).
\end{eqnarray}
For $f$ evaluating to infinite values we assume the conventions that $\infty-\infty=\infty$, $\infty\geq x$, and $x\geq-\infty$ for every $x\in\mathbb{R}\cup\{-\infty, \infty\}$. We refer the reader to Figure \ref{fig:inequalities} in Appendix {\ref{sec:appendix_lemma_multimod}} for a visual illustration of the diminishing return properties described by these inequalities.
\end{definition}

\DoubleSpacedXI

\vspace{-.2in}

\subsubsection*{System-wide objective.}
Our goal is to minimize the combined number of out-of-stock events across the system, i.e., $\sum_ic_i(d_i,b_i)$. Writing $\vec{d}$ and $\vec{b}$ for the vectors that contain $d_i$ and $b_i$ in their~$i$th position we denote this sum by $c(\vec{d},\vec{b})$. We minimize $c(\vec{d},\vec{b})$ subject to four kinds of constraints that we introduce now.
\subsubsection*{Constraints.} Our optimization problem involves two kinds of budget constraints. The first is on the total number of bikes allocated, that is, $\sum_i b_i$, bounded by $B$. The second is on the total number of docks allocated in the system, that is, $\sum_id_i+b_i$, which is bounded by $D+B$. In addition, we have an \emph{operational constraint} that bounds the number of docks that can be reallocated within the system. To formally state this constraint it is useful to define the following.
\begin{definition}\label{def:ball}
{Consider two allocations $(\vec{d},\vec{b})$ and $(\vec{d}',\vec{b}')$ with $|\vec{d}+\vec{b}|_1=|\vec{d'}+\vec{b'}|_1$, i.e., the same number of docks allocated in total. The number of docks that need to be reallocated to get from $(\vec{d},\vec{b})$ to $(\vec{d}',\vec{b}')$ (ignoring the allocation of bikes) is $|\vec{d'}+\vec{b'}-\vec{d}-\vec{b}|_1/2$. 
For allocations $(\vec{d},\vec{b})$ and $(\vec{d}',\vec{b}')$ with the same total number of docks we define this as the \emph{dock-move distance} between them.} 
\end{definition}
Given an initial allocation $(\vec{\bar{d}},\vec{\bar{b}})$, for which we assume $\sum_i \bar{d}_i+\bar{b}_i\leq D+B$, the operational constraint is then of the form {$|\vec{\bar{d}}+\vec{\bar{b}}-\vec{d}-\vec{b}|_1/2\leq z$}
for some $z$ {(where this constraint is well-defined even when $|\vec{\bar{d}}+\vec{\bar{b}}|_1\neq|\vec{d}+\vec{b}|_1$)}. Finally, we have physical constraints that give lower and upper bounds on the number of docks allocated to each station $i$, where we assume that~${0\leq} l_i\leq \bar{d}_i+\bar{b}_i\leq u_i$ for every~$i$. The resulting optimization problem can be written as follows:

\OneAndAHalfSpacedXI
\begin{equation*}
\begin{aligned}
\mathtt{minimize}_{(\vec{d},\vec{b})\in{\mathbb{N}^n\times \mathbb{N}^n}} & c(\vec{d},\vec{b})\\
\mathtt{s.t.}  & \sum_i d_i+b_i &\leq D+B,\\
& \sum_i b_i & \leq B, \\
& 
{|\vec{\bar{d}}+\vec{\bar{b}}-\vec{d}-\vec{b}|_1/2\leq z},\\
& l_i\leq d_i+b_i \leq u_i & \forall i.
\end{aligned}
\tag{P1}\label{opt:original}
\end{equation*}
\DoubleSpacedXI

Following standard convention, we define $c_i(d_i,b_i)=\infty$ {for $d_i+b_i>u_i$ and $d_i+b_i<l_i$, which allows us to drop the last row of inequalities in \ref{opt:original}. We also set $c_i(d,b)=\infty$ for $d<0$ or~$b<0$. With these changes, $c_i(\cdot,\cdot)$ fulfills the inequalities in Definition \ref{def: multimod}.
%
\begin{lemma}\label{lemma: multimod}
The function $c_i(\cdot,\cdot)$ is multimodular.
\end{lemma}
The proof of the lemma is based on a coupling argument, and appears in Appendix \ref{sec:appendix_lemma_multimod}.}
In addition, we may add a {$(n+1)$st} dummy (``depot'') station $\mathcal{D}$ to guarantee that the first two constraints hold with equality in an optimal solution. Thus, we can transform our original optimization problem into one of the following form, which is our focus throughout the main body of the text:
\OneAndAHalfSpacedXI
\begin{equation*}
\begin{aligned}
\mathtt{minimize}_{(\vec{d},\vec{b})\in{\mathbb{Z}^{n+1}\times \mathbb{Z}^{n+1}}} & c(\vec{d},\vec{b}) &\\
\mathtt{s.t.} &  \sum_i d_i+b_i &= D+B,\\
& \sum_i b_i &= B, \\
& {|\vec{\bar{d}}+\vec{\bar{b}}-\vec{d}-\vec{b}|_1/2\leq z, \text{ where } |\vec{\bar{d}}+\vec{\bar{b}}|_1=D+B}.
\end{aligned}
\tag{P2}\label{opt:adapted}
\end{equation*}
\DoubleSpacedXI
Specifically, the reduction from Problem \ref{opt:original} to \ref{opt:adapted} is based on the following: let $\bar{D}=D+B-\sum_i \bar{d}_i+\bar{b}_i$, i.e., the number of docks that are not in the current allocation but can be added, $\bar{z}=z+\lfloor \frac{\bar{D}}{2}\rfloor$, and define station $\mathcal{D}$ with $l_{\mathcal{D}}=B, u_{\mathcal{D}}=2B+D$, $\bar{d}_{\mathcal{D}}+\bar{b}_{\mathcal{D}}=B+\bar{D}$, and $c_{\mathcal{D}}(d,b)=d+b-B$ when~$l_{\mathcal{D}}\leq d+b\leq u_{\mathcal{D}}$ --- observe that $c_{\mathcal{D}}$ fulfills the requirements of Definition \ref{def: multimod}{. Further, $c_{\mathcal{D}}$ has the property that its objective is increasing in the number of docks allocated to it --- whereas the objective at all other stations is non-increasing in the number of docks allocated. In the proof of the following proposition, this will be used to ensure that optimal solutions to \ref{opt:adapted} fulfill $d_\mathcal{D}+b_\mathcal{D}=l_\mathcal{D}$. It is worth noting that our algorithm/analysis for \ref{opt:adapted} does not rely on the $c_i$ being non-increasing in the number of docks allocated, i.e., $c_{\mathcal{D}}$ being decreasing will not affect our analysis.}. 
\begin{proposition}\label{prop:p1p2equiv}
If $(\vec{d},\vec{b})$ is optimal for \ref{opt:adapted} with stations $[n]\cup\{\mathcal{D}\}$,  bike budget $B$, dock budget~$D+2B$, and operational constraint $\bar{z}$, then restricting $(\vec{d},\vec{b})$ to $[n]$ is optimal for \ref{opt:original}.
\end{proposition}

The proof of the proposition is in Appendix \ref{appendix:tradeoff}. There, we also show how to optimally solve an optimization problem that involves an additional trade-off between the size of $D$ and the size of $z$, i.e., between the inventory cost of additional docks and the operational cost of reallocating docks. We are now ready to discuss the assumptions in our model before analyzing in Section \ref{sec: algorithm} an algorithm to optimally solve \ref{opt:adapted}.

\subsection{Discussion of Assumptions}
\label{sec: assumptions}
Before describing and analyzing the algorithm we use to solve the optimization problem in Section~\ref{sec: algorithm}, we discuss here the assumptions as well as the advantages that come along with them.

\noindent\textbf{Seasonality and Frequency of Reallocations. }
In contrast to bike rebalancing, the reallocation of docks is a strategic question that involves docks being moved at most annually. As such, a concern is that the recommendations for a particular month might not yield improvement for other times of the year. One way to deal with this is to explicitly distinguish, in the demand profiles, between different seasons, i.e., have $k$ different distributions for $k$ different types of days and then consider the expectation over these as the objective. Though the user dissatisfaction functions accommodate that approach, we find on real data (see Section \ref{ssec: objective}) that this is not actually necessary: the reallocations that yield greatest impact for the summer months of one year also perform very well for the winter months of another. This even held true in New York City, where the system significantly expanded year-over-year: despite the number of stations in the system more than doubling and total ridership increasing by around 70\% from 2015 to 2017, we find that the estimated improvement due to reallocated docks is surprisingly stable across these different months. In part this is due to the fact that the relative demand patterns at different stations strongly correlate between seasons, i.e., the demand of each station in each interval in one month is well-approximated by a constant multiple of demand in another. For example, the vectors of half-hourly demand estimates (either rentals or returns) for each New York City station in June and December 2018 have a Pearson correlation coefficient greater than $0.85$. Though this does not formally imply that the improvement in the UDFs would correlate, it gives some explanation for why it might.

\noindent\textbf{Cost of Reallocation. }
Rather than explicitly building in a cost for reallocations in our formulation, we instead bound the number of docks that are moved.  
This is mostly motivated by {our industry partner's practical considerations}: the cost of physically reallocating capacity from one location to another is negligible when compared to the administrative effort, a negotiation with city officials and other stakeholders, needed to reallocate capacity. {As part of these negotiations the operator will request that a limited number of docks be moved from the current system configuration. While we solve the problem assuming that this number is a known constant, in practice it is part of the negotiations.  
To hold these negotiations it is of utmost importance for the operator to know the value of reallocating a given (fixed) number of docks; thus, our results were used to help prepare the operator for these negotiations, in particular, to answer for different values of $z$ the crucial questions of \emph{how much benefit would the system derive from moving $z$ docks} and \emph{which docks would be among those $z$}. This then also implies that tactical questions of} how to carry out the reallocations is of minor importance in practice. Further, the cost of reallocating docks can be compared to the cost of rebalancing bikes: while the (one-off) reallocation of a single dock is about an order of magnitude more expensive than that of a single bike, the reallocated dock has daily impact on improved service levels (in contrast to the one-off impact of a rebalanced bike). Thus, the cost quickly amortizes; Citi Bike estimates in as little as 2 weeks. Finally, the cost to acquire new docks is orders of magnitudes higher than all of the aforementioned costs, leading us to focus only on reallocated capacity in our analysis; nevertheless, we show in {Appendix \ref{appendix:tradingoff_new_allocated}} that the algorithm also extends to capture the tradeoff between installing newly bought and reallocating existing docks. 

\noindent\textbf{Bike Rebalancing. }  The user dissatisfaction functions assume that no rebalancing takes place over the course of the planning horizon. System data indicates that this is close to reality at most stations; for example, in New York City, more than 60\% of all rebalancing is concentrated at just~28 out of 762 stations which justifies the assumption for the vast majority of stations.
{Now, consider the remaining few stations, at which almost all rebalancing is concentrated: perhaps unsurprisingly we find that none of these stations are identified by the optimization as having their capacity reduced. In general, rebalancing can always limit the number of dissatisfied users to 0: consider a station with 2 docks that is stocked with 1 bike; as long as rebalancing adds/removes a bike after each pickup/dropoff, users will not experience stockouts. Thus, reducing the number of dissatisfied customers at a station with no rebalancing is somewhat analogous to reducing rebalancing needs at a station with rebalancing. For illustrative purposes consider the following deterministic example: a station with 60 docks observes demand for 120 rentals in the morning and demand for 120 dropoffs in the afternoon. Suppose the station is full with bikes in the morning. Without rebalancing, the station observes $60-x$ stockouts in the morning, and $60-y$ stockouts in the afternoon, where $x,y\leq 60$ are respectively the number of bikes rebalancing drops off in the morning/picks up in the afternoon. With 15 docks added these quantities would turn into $45-x$ and $45-y$ for $x,y\leq 45$. Thus, the same amount of rebalancing, up to a smaller upper bound, would simply reduce the number of dissatisfied customers (by the amount captured by the UDFs); beyond that upper bound additional rebalancing is no longer needed. This example aligns with anecdotal experiences system operators have shared with us: stations that had dock capacity added to them subsequently required less rebalancing.}

Though we assume that no rebalancing occurs over the course of the planning horizon, the optimization model assumes that the initial number of bikes at each station is optimally allocated. We relax this assumption in Section \ref{sec: long_run} when we consider a regime in which no rebalancing occurs at all. Despite the fact that the two regimes can be viewed as polar opposites (optimally rebalanced overnight and no rebalancing overnight), our results indicate that they yield very similar recommendations for the operators. Our motivation to focus on these opposite extremes is simple: modeling a modest amount of rebalancing poses significant challenges. For example, unlike the effect of daily usage patterns, overnight rebalancing is affected by greater variability from external factors, ranging from the number of trucks to the supply of just-repaired bikes.



\noindent\textbf{Exogenous Rentals and Returns.}
The demand profiles assume that the sequences of arrivals are exogenous, i.e., there is a fixed distribution that defines the sequence of rentals and returns at each station. Before justifying this assumption, it is worth considering a setting in which it fails spectacularly: consider an allocation of bikes and docks that allocates no bikes at all. With no bikes, no attempted rental is ever successful and therefore no returns ever occur. As such, the sequence of arrivals of returns at one station are not independent of the allocations elsewhere. 

Another extreme arises where the stations never run out of bikes, and there is always capacity available to receive bike returns. In this case, bike rentals and returns proceed smoothly independent of allocations, and so can be viewed as exogenously given. This ideal case is the one to which we strive in our reallocation efforts. Of course, the assumption is never realized exactly in practice.

At this point it is helpful to discuss the stochastic model for rentals and returns that we use in our calculations. Suppose that at each station, potential bikers arrive according to a Poisson process that is independent of that at all other stations. Such a model is plausible because of the Palm-Khintchine theorem that states, roughly speaking, that the superposition of the bike rentals of a large number of independent users is well modeled by a Poisson process; see, e.g., p.~221 of \cite{kartay75}, \cite{cin72}, p.~107 of \cite{nel13}. 
Also, suppose that users select their destinations according to an origin-destination routing matrix, thereby splitting the Poisson incoming flows into independent biker flows. Assuming biking times between any fixed pair of stations are identically distributed, and are independent across all bikers and station pairs, it follows that the process of returning bikes at a destination station from a fixed origin station, being a delayed Poisson process, is again a Poisson process. But then, the overall bike-return process at the destination station, being a superposition of such flows from all origin stations, is again a Poisson process. Moreover, due to the splitting property of Poisson processes, the rental-return processes at each destination station are mutually independent. Thus, at each station it is reasonable to model the returns and rentals of bikes as Poisson processes, justifying the exogenous arrivals assumption. This modeling structure is approximate {for several reasons: (i) the Poisson flows entering a destination station are interrupted if an upstream station runs out of bikes, (ii) a destination station may observe additional returns due to nearby stations being full, and (iii) a station may observe additional demand due to nearby stations being empty}. As mentioned earlier, we attempt to minimize such shortages, so that we strive for conditions under which the approximation is close to reality, though it is still an approximation. Under this Poisson model, the rental and return processes at different stations are not independent. For example, a surge in rentals at one station may result in a surge of returns at a ``downstream'' station. Fortunately, our objective function is additively separable in stations, so independence at different stations is not required to compute the objective function; the ``marginal'' property that flows are Poisson at each station considered individually suffices.

Perhaps an even stronger justification for the exogenous-arrivals assumption comes from work by {\cite{jianetal} and \cite{datner2015setting} who both use simulation optimization approaches. 
\cite{datner2015setting} use their simulation optimization approach to identify only the optimal allocation of bikes. They endogenize (i)-(iii) above and compare their results to optimizing with the UDF. Though they focus on a slightly different objective (\emph{total user travel time}, where stockouts may lead to pickups/dropoffs at other stations or to users walking), they also report the fraction of rides affected by stockouts, which is what we/the UDFs aim to minimize; for this objective, their solutions improve upon the UDFs by only 1.2\% on average (across 6 scenarios).
Similarly to us, \cite{jianetal} aim to find the configuration of bikes and docks across the system that minimizes the number of out-of-stock events over the course of the day. In contrast to the user dissatisfaction functions, decensoring the demand data for their simulation required additional modeling decisions that allow them to endogenize (i) and (ii) above. }
 While this simulation approach still assumed that demand for rentals was exogenous, it endogenized returns, excluding (at least) the example suggested above. However, it causes the resulting simulation optimization problem to be non-convex in an unbounded fashion. Indeed, for any bound $L$, one can construct highly contrived examples in which there exists an initial allocation $(\vec{d},\vec{b})$ and stations $i$ and $j$ such that when starting at allocation $(\vec{d},\vec{b})$ it is the case that (a) moving two bikes from $i$ to $j$ improves the objective by at least $L$ and (b) moving one bike from $i$ to $j$ gives a solution that is worse than $(\vec{d},\vec{b})$. Such examples not only show that the objective function in that model is non-convex, they also show that solutions from such a framework are harder to interpret.
  \cite{jianetal} proposed a range of different gradient-descent algorithms as heuristics to find good solutions, including adaptations of the algorithms we present and analyze here. Despite the simulation adding key complexities to the system, the heuristics gave only limited improvements -- approximately 3\% -- when given the solution found by our algorithms as a starting point. {Thus, there exists substantial data-driven evidence to justify the use of UDFs.} 
  
Finally, the assumptions that rentals and returns are exogenous, and the objective is separable across stations, are quite common in the rebalancing literature. This includes for example \cite{raviv2013optimal}, \cite{raviv2013static}, \cite{di2013hybrid}, \cite{rainer2013balancing}, \cite{raidl2013balancing}, \cite{ho2014solving}, \cite{kloimullner2014balancing}, \cite{kaspi2015bike}, \cite{forma20153},  \cite{alvarez2016optimizing}, and \cite{schuijbroek2013inventory}, most of whom make the assumption implicitly.


\noindent\textbf{Out-of-stock Events and Demand Profiles. }
In practice, we cannot observe attempted rentals at empty stations nor can we observe attempted returns at full stations. Worse still, given that most bike-sharing systems have mobile apps that allow customers to see real-time information about the current number of bikes and empty docks at each station, there may be customers who want to rent a bike at a station, see on the app that the station has few bikes available presently, and decide against going to the station out of concern that by the time they {arrive, the remaining bikes will} already have been taken by someone else. Should such a case be considered an out-of-stock event (respectively, an attempted rental)? The user dissatisfaction functions assume that such events do not occur as the definition relies on out-of-stock events occurring only when stations are either entirely empty or entirely full.

Further, in order to compute the user dissatisfaction functions, we need to be able to estimate the demand profiles: using only observed rentals and returns is insufficient as it ignores latent demand at empty/full stations. To get around this, we mostly apply a combination of approaches by \cite{o2015data}, \cite{henomashm16}, and \cite{parikh2014estimation}: we estimate Poisson arrival rates (independently for rentals and returns) for each 30 minute interval and use a formula developed by \cite{henomashm16} to compute, for any initial condition (in number of bikes and empty docks) the expected number of out-of-stock events over the course of the interval. We plug these into a stochastic recursion suggested by \cite{parikh2014estimation} to obtain the expected number of out-of-stock events over the course of a day as a function of the number of bikes and empty docks at 6AM. This is far from being the only approach to compute user dissatisfaction functions; for example, in Section \ref{sec: impact} we explicitly combine empirically observed arrivals with estimated rates for times when rentals/returns are censored at empty/full stations. 



\noindent\textbf{Advantages of User Dissatisfaction Functions. }
The user dissatisfaction functions yield several advantages over a more complicated model such as the simulation. First, they provide a computable metric that can be used for several different operations: in Section \ref{sec: algorithm} we show how to optimize over them for reallocated capacity and in Section \ref{sec: impact} we use them to evaluate the improvement from already reallocated capacity. \cite{chung2018} used them to study an incentive program operated by Citi Bike in New York City, and they have been used extensively for motorized rebalancing (see Section \ref{sec: rel_work}). As such, the user dissatisfaction functions provide a single metric on which to evaluate different operational efforts to improve service quality, which adds value in itself. Second, for the particular example of reallocating dock capacity that we study here, they yield a tractable optimization problem, which we prove in Section \ref{sec: algorithm}. Third, for the reallocation of dock capacity, the discrete convexity properties we prove imply that a partial implementation of the changes suggested by the optimization (see Section \ref{sec: case_study}) is still guaranteed to yield improvement. Finally, given a solution to the optimization problem, it is easy to track the partial contribution to the objective from changed capacity at each station, making solutions interpretable.

\section{A Discrete Gradient-Descent Algorithm}\label{sec: algorithm}
We begin this section by examining the mathematical structure of Problem \ref{opt:adapted} that allows us to develop efficient algorithms.
In Section~\ref{sec: algorithm_opt}, we define a natural neighborhood structure on the set of feasible allocations and define a discrete gradient-descent algorithm on this neighborhood structure. We prove in Section~\ref{sec: algopt} that for the problem without operational constraints (\ref{opt:adapted} with $z=\infty$), solutions that are locally optimal with respect to the neighborhood structure are also globally optimal; since our algorithm continues to make local improvements until it finds a local optimum, this proves that the solution returned by our algorithm must be globally optimal. Finally, in Section~\ref{sec: kstepopt}, we prove that the algorithm takes at most $z$ iterations to find the best allocation obtainable by moving at most $z$ docks within the system (see Definition \ref{def:neighborhood} for a formal definition of a move of a dock); this not only proves that the gradient-descent algorithm optimally solves the minimization problem when including operational constraints, but also guarantees that doing so requires at most $D+B$ iterations.






%
%

\subsection{Algorithm}\label{sec: algorithm_opt}


We now present our algorithm before analyzing it for settings without the operational constraints. Intuitively, in each iteration our discrete gradient-descent algorithm picks  one dock and at most one bike within the system and moves them from one station to another. 
It chooses the dock, and the bike, so as to maximize the reduction in objective value, i.e., in a discrete sense it executes a gradient-descent step. To formalize this notion, we define the \emph{movement of a dock} via the following transformations. Denote by $e_i=(0,\ldots,1,\ldots,0)$ the canonical unit vector that has a 1 in its $i$th position and 0s elsewhere. 

\begin{definition}\label{def:neighborhood}
A \emph{dock-move} from $i$ to $j$ corresponds to one of the following transformations of feasible solutions:
\begin{enumerate}
\item $o_{ij}(\vec{d},\vec{b})=\big(\vec{d}-e_i+e_j,\vec{b}\big)$  -- Moving one (empty) dock from $i$ to $j$;
\item $e_{ij}\big(\vec{d},\vec{b})=\big(\vec{d},\vec{b}-e_i+e_j\big)$ -- Moving one dock and one bike from $i$ to $j$, i.e., one full dock;
\item $E_{ijh}(\vec{d},\vec{b})=\big(\vec{d}-e_i+e_h, \vec{b}+e_j-e_h\big)$ -- Moving {one} dock from $i$ to $j$ and {one} bike from $h$ to~$j$;
\item $O_{ijh}(\vec{d},\vec{b})= \big(\vec{d}+e_j-e_h, \vec{b}-e_i+e_h\big)$ -- Moving {one} bike from $i$ to $h$ and {one} dock from $i$ to~$j$. 
\end{enumerate}
We often refer to the first kind as \emph{moving an empty dock} from $i$ to $j$ and to the second kind as \emph{moving a full dock} from $i$ to $j$ to indicate that the dock is moved by itself (empty) or with a bike (full). Without qualification, the movement of a dock can refer to any of the above.
Further, we define the \emph{neighborhood} $N(\vec{d},\vec{b})$ as the set of allocations that are one dock-move away from $(\vec{d},\vec{b})$: 
$$
N(\vec{d},\vec{b}):=\{o_{ij}(\vec{d},\vec{b}), e_{ij}(\vec{d},\vec{b}), E_{ijh}(\vec{d},\vec{b}), O_{ijh}(\vec{d},\vec{b}): i,j,h\in[n]\}.
$$
Notice that $(\vec{d},\vec{b})\in N(\vec{d}',\vec{b}')$ implies that
\dfedit{$|\vec{d}+\vec{b}-\vec{d}'-\vec{b}'|_1=2$}
, i.e., a dock-move distance of~1; the converse however does not hold true as the dock-move distance between two allocations does not take into account their allocation of bikes.
\end{definition}
Throughout the paper we also sometimes refer to the move of a bike from $i$ to $j$, by which we mean a transformation from $(\vec{d},\vec{b})$ to $(\vec{d}+e_i-e_j,\vec{b}-e_i+e_j)$. This changes the allocations of bikes to stations while keeping the number of docks at each station constant.

The above-defined neighborhood structure  gives rise to a very simple algorithm (see Algorithm~\ref{alg:gradient_descent} in Appendix \ref{appendix:pseudocode}): we first find an optimal allocation of bikes for the current allocation of docks, i.e., when each station $i$ is restricted to have $\bar{d}_i+\bar{b}_i$ docks allocated to it (see Algorithm~\ref{alg:bike_optimal}).  The convexity of each $c_i$ in the number of bikes (see Fact 1 in Appendix \ref{appendix:bike_opt}), with fixed number of docks, implies that this can be done greedily by taking out all the bikes and then adding them one by one (see e.g., first page of \cite{hochbaum1994lower}). Denote this allocation by $(\vec{d}^0,\vec{b}^0)$. Then, iterate over $r\in\{1,\ldots,z\}$, i.e., through $z$ periods, by either setting $(\vec{d}^{r},\vec{b}^{r})$ to be the best allocation in the neighborhood of $(\vec{d}^{r-1},\vec{b}^{r-1})$, or if that allocation is no better than $(\vec{d}^{r-1},\vec{b}^{r-1})$, returning $(\vec{d}^{r-1},\vec{b}^{r-1})$ (see Algorithm~\ref{alg:gradient_descent_step}). After iteration $z$, the algorithm returns $(\vec{d}^{z},\vec{b}^{z})$.


\subsection{Optimality without Operational Constraints}\label{sec: algopt}
We first prove that Algorithm \ref{alg:gradient_descent} returns an optimal solution {to the problem without operational constraints}.
Specifically, we analyze the following:
\OneAndAHalfSpacedXI
\begin{equation*}
\begin{aligned}
\mathtt{minimize}_{(\vec{d},\vec{b})} & c(\vec{d},\vec{b})\\
\mathtt{s.t. } &  {|\vec{d}+\vec{b}|_1} &= D+B,\\
& {|\vec{b}|_1} & = B
\end{aligned}
\tag{P3}\label{opt:no_ops}
\end{equation*}
\DoubleSpacedXI

We show that an allocation  $(\vec{d},\vec{b})$ that is locally optimal with respect to $N(\cdot,\cdot)$ must also be globally optimal with respect to \ref{opt:no_ops}. Thus, if Algorithm \ref{alg:gradient_descent}, initialized with $z=\infty$, finds a solution for which there is no better solution in the neighborhood, then it returns an optimal solution to \ref{opt:no_ops}. However, since that solution may {have dock-move distance greater $z$ to $(\vec{\bar{d}},\vec{\bar{b}})$}, 
global optimality of the algorithm only follows for \ref{opt:no_ops}, not for \ref{opt:adapted}. Before we prove Lemma~\ref{lemma: neighborhoods} to establish this, we first define an allocation of bikes and docks as \emph{bike-optimal} if it minimizes the objective among allocations with the same number of docks at each station.
\begin{definition}
$(\vec{d},\vec{b})$ is \emph{bike-optimal} if~$(\vec{d},\vec{b}) \in \arg \min_{(\vec{d'},\vec{b'}):\forall i, d_i+b_i=d'_i+b'_i, \; {|\vec{b'}|_1} =B}\{c(\vec{d'},\vec{b'})\}$.
\end{definition}
The following lemma ensures that our analysis can, for the most part, focus only on bike-optimal solutions. 
\begin{lemma}\label{lemma: bike_opt}
\noindent
Suppose $(\vec{d},\vec{b})$ is bike-optimal. Then given any $i$ and $j$, the allocation resulting from the best dock-move from $i$ to~$j$ is bike-optimal.
\end{lemma}
The proof of the lemma is in Appendix \ref{appendix:bike_opt}. We highlight three implications of the lemma:
\begin{enumerate}
    \item Since Algorithm \ref{alg:gradient_descent} finds a bike-optimal solution initially and picks the best dock-move in each iteration, bike-optimality is an invariant of the algorithm, i.e.,  $(\vec{d}^0,\vec{b}^0),(\vec{d}^1,\vec{b}^1), (\vec{d}^2,\vec{b}^2),\ldots$ are all bike-optimal. 
    \item Consider a bike-optimal allocation $(\vec{d},\vec{b})$ and an allocation $(\vec{d}',\vec{b}')$ at dock-move distance~1 that is not in $N(\vec{d},\vec{b})$; then there exists a dock-move from $(\vec{d},\vec{b})$ that creates an allocation with (i) the same allocation of docks at each station as $(\vec{d}',\vec{b}')$, i.e., $d_i'+b_i'$ at each station $i$, and (ii) an objective no worse than $(\vec{d}',\vec{b}')$.
    \item To prove optimality of Algorithm \ref{alg:gradient_descent} for \ref{opt:no_ops} it  suffices to prove that bike-optimal solutions that are locally optimal w.r.t.\ our neighborhood structure are also globally optimal.
\end{enumerate}

We formalize the last of these in Lemma~\ref{lemma: neighborhoods}, the proof of which we defer to Appendix \ref{sec:appendix_lemma_3}. That appendix also contains the proof of Lemma \ref{lemma: exchange}, of which Lemma~\ref{lemma: neighborhoods} is a corollary.



\begin{lemma}\label{lemma: neighborhoods}
Suppose $(\vec{d},\vec{b})$ is bike-optimal, but not optimal for {either \ref{opt:adapted} or \ref{opt:no_ops}}. Let $(\vec{d}^\star,\vec{b}^\star)$ denote a better solution.  Then there exists {$(\vec{d'},\vec{b'})\in N(\vec{d},\vec{b})$ such that $(\vec{d'},\vec{b'})$ has both a lower objective value and a smaller dock-move distance to $(\vec{d}^\star,\vec{b}^\star)$ than $(\vec{d},\vec{b})$ does.}
\end{lemma}

{In the statement of Lemma \ref{lemma: neighborhoods} the allocation $(\vec{d'},\vec{b'})$ is not restricted to fulfill the operational constraints, {i.e., it may not be feasible for Problem \ref{opt:adapted}} (else, this would already imply that Algorithm~\ref{alg:gradient_descent} always, eventually, terminates with an optimal solution). Thus, optimality only follows for \ref{opt:no_ops}. The allocations identified in the next lemma, $(\vec{d}',\vec{b}')$ and $(\vec{d}^{\star\star},\vec{b}^{\star\star})$, also need not satisfy the operational constraints.}


\begin{lemma}\label{lemma: exchange}
Consider any bike-optimal solution $(\vec{d},\vec{b})$ and a better     allocation $(\vec{d}^\star,\vec{b}^\star)$; let $j$ and $k$ denote stations with $d_j+b_j<d_j^\star+b_j^\star$ and $d_k+b_k>d_k^\star+b_k^\star$. Then either there exist $(\vec{d'},\vec{b'})\in N(\vec{d},\vec{b})$ with $d_j'+b_j' = d_j+b_j+1$ and $(\vec{d}^{\star\star},\vec{b}^{\star\star})\in N(\vec{d}^\star,\vec{b}^\star)$ with $d_j^{\star\star}+b^{\star\star}_j=d^\star_j+b^\star_j-1$ or there exist $(\vec{d'},\vec{b'})\in N(\vec{d},\vec{b})$ with $d_k'+b_k' = d_k+b_k-1$ and $(\vec{d}^{\star\star},\vec{b}^{\star\star})\in N(\vec{d}^\star,\vec{b}^\star)$ with $d_k^{\star\star}+b^{\star\star}_k=d^\star_k+b^\star_k+1$  such that 
\begin{enumerate}
\item $c(\vec{d},\vec{b})-c(\vec{d'},\vec{b'})\geq c(\vec{d}^{\star\star},\vec{b}^{\star\star})-c(\vec{d}^{\star},\vec{b}^{\star})$
\item the dock-move distance from $(\vec{d'},\vec{b'})$ to $(\vec{d}^\star,\vec{b}^\star)$ is less than from $(\vec{d},\vec{b})$ and the dock-move distance from $(\vec{d}^{\star\star},\vec{b}^{\star\star})$ to $(\vec{d},\vec{b})$ is less than from $(\vec{d}^\star,\vec{b}^\star)$
{\item the dock-move from $(\vec{d},\vec{b})$ to $(\vec{d'},\vec{b'})$ yields $(\vec{d}^\star,\vec{b}^\star)$ when applied to $(\vec{d}^{\star\star},\vec{b}^{\star\star})$ (so, e.g., if $e_{kj}(\vec{d},\vec{b}) = (\vec{d'},\vec{b'})$, then $e_{kj}(\vec{d}^{\star\star},\vec{b}^{\star\star})=(\vec{d^\star},\vec{b}^\star)$,  or equivalently $e_{jk}(\vec{d^\star},\vec{b}^\star)=(\vec{d}^{\star\star},\vec{b}^{\star\star})$).}
\end{enumerate}
\end{lemma}

\emph{Remark}: In the discrete convexity literature, a rewriting of the objective allows this to be interpreted as the exchange property of $M^\natural$-convex functions; 
this connection has been explored in the follow-up work of \cite{shioura2018m}.



\subsection{Operational Constraints \& Running Time}\label{sec: kstepopt}

In this section, we show that Algorithm \ref{alg:gradient_descent} is optimal for \ref{opt:adapted} by proving that{, for any $r$, in $r$ iterations it finds the best allocation obtainable by moving at most $r$} docks. We thereby also provide an upper bound on the running-time of the algorithm, since an optimal solution can be at most $\min\{D+B, z\}$ dock-moves apart from $(\vec{d}^0,\vec{b}^0)$.



Our proof works inductively. We begin by showing (Lemma \ref{lemma:const_neighborhood}) that, assuming that $(\vec{d}^{r},\vec{b}^{r})$ minimizes the objective {among solutions at dock-move distance at most $r$ to $(\vec{d}^0,\vec{b}^0)$}, 
$(\vec{d}^{r+1},\vec{b}^{r+1})$ {must be a local optimum 
among solutions at dock-move distance at most $r+1$ to $(\vec{d}^0,\vec{b}^0)$}.
{This local optimality in the $({r}+1)$st iteration} 
guarantees a particular structural property (see Lemma \ref{lemma: iexists} in Appendix \ref{appendix:proof_theorem} for details). 
{In} Theorem~\ref{thm: koptimal} we {use} 
this structural property, together with the optimality of the solution in the ${r}$th iteration and the gradient-descent step, {to show} that $(\vec{d}^{r+1},\vec{b}^{r+1})$ is globally optimal among solutions {with dock-move distance at most $r+1$ to $(\vec{d}^{0},\vec{b}^0)$.} 
The proofs of Lemma~\ref{lemma:const_neighborhood}, Lemma~\ref{lemma: iexists}, 
and Theorem~\ref{thm: koptimal} can be found in Appendix~\ref{appendix_3.4}.

\begin{lemma}\label{lemma:const_neighborhood}
Suppose $(\vec{d}^{{r}},\vec{b}^{{r}})$ minimizes the objective among solutions {with $|\vec{d}+\vec{b}-\vec{d}^{{r}}-\vec{b}^{{r}}|_1\leq 2r$} 
and let $(\vec{d}^{{{r}}+1},\vec{b}^{{{r}}+1})$ denote  the next choice of the gradient-descent algorithm, i.e., an allocation in the neighborhood of $(\vec{d}^{{r}},\vec{b}^{{r}})$ that minimizes the objective. Then $(\vec{d}^{{{r}}+1},\vec{b}^{{{r}}+1})$ is a local optimum {among solutions with $|\vec{d}+\vec{b}-\vec{d}^{{r}}-\vec{b}^{{r}}|_1\leq 2(r+1)$, 
that is, there is no solution in $N(\vec{d}^{{{r}}+1},\vec{b}^{{{r}}+1})$ that is at dock-move distance at most $r+1$ to $(\vec{d},\vec{b})$ and has a lower objective than $(\vec{d}^{{{r}}+1},\vec{b}^{{{r}}+1})$}.
\end{lemma}



\begin{theorem}\label{thm: koptimal}
{Initialized at $(\vec{d},\vec{b})$} Algorithm \ref{alg:gradient_descent} finds in {iteration} $r$ an allocation that minimizes the objective among {those at dock-move distance at most equal to $r$}.
\end{theorem}

\emph{Remark}: In an earlier manuscript, as well as the proceedings version of the paper, we had falsely stated that local optima {are globally optimal among allocations with dock-move distance at most~$r$ to~$(\vec{d},\vec{b})$, }
i.e., that Lemma~\ref{lemma: neighborhoods} continues to hold in the setting where {allocations at dock-move distance~$r+1$ } 
are infeasible, and immediately derived optimality from that and Lemma \ref{lemma:const_neighborhood}. An example by \cite{shioura2018m} shows that this is false: there exist solutions that are locally optimal with respect to our neighborhood structure and {at dock-move distance at most $r$ to $(\vec{d},\vec{b})$, despite not being globally optimal among the solutions at dock-move distance at most $r$. } 
\cite{shioura2018m} also provides an alternative proof of correctness for our algorithm.


 






\section{Scaling Algorithm}\label{sec: scaling}
We now extend our analysis in Section \ref{sec: algorithm} to adapt our algorithm to a scaling algorithm that provably finds an optimal allocation of bikes and docks {for Problem \ref{opt:no_ops}, i.e.,} the setting without \emph{operational} constraints, in~$O\big(n\log(B+D)\big)$ iterations. 

The idea underlying the scaling algorithm (see Algorithm \ref{alg:scaling} in Appendix \ref{appendix:pseudocode}) is to proceed in $\lfloor\log_2(B+D) \rfloor+1$ phases, where in the $k$th phase each move involves $\alpha_k=2^{\lfloor\log_2(B+D) \rfloor+1-k}$ bikes/docks rather than just one. The $k$th phase is prefaced by finding the bike-optimal allocation of bikes (given the constraints of only moving $\alpha_k$ bikes at a time), and terminates when no move of $\alpha_k$ docks yields improvement. {We first observe} that the multimodularity of $c(d,b)$ implies multimodularity of $c(\alpha_kd,\alpha_kb)$ for all $k$ (e.g., Table 1.3 in \cite{murotacoss2018}). Thus, our analysis in the last section implies that in the $k$th phase, the scaling algorithm finds an optimal allocation among all that differ in a multiple of $\alpha_k$ in each coordinate from $(\vec{\bar{d}},\vec{\bar{b}})$. Further, since $\alpha_{\lfloor\log_2(B+D) \rfloor+1}=1$, it finds the globally optimal allocation in phase $\lfloor\log_2(B+D) \rfloor+1$. 


 \begin{theorem}\label{thm:scaling}
Algorithm \ref{alg:scaling} finds an optimal allocation for \ref{opt:no_ops} in $O\big(n\log(B+D)\big)$ iterations.
 \end{theorem}

Notice that Theorem \ref{thm:scaling} can provide a non-trivial speedup, relative to $O(n+D+B)$ for the gradient-descent algorithm, when $B$ and $D$ are large relative to $n$, i.e., when stations have many docks on average. Otherwise it may create unnecessary overhead; in Appendix \ref{sec: running_time} we observe this on real data when comparing the different algorithms on data sets from different cities. Motivated by this insight, we then also define a hybrid algorithm (see Algorithm \ref{alg:hybrid}) that proceeds like the scaling algorithm but skips some of the phases. {The proof of Theorem \ref{thm:scaling}, included  
in Appendix~\ref{sec:appendix_proof_scaling}, relies on a proximity result that bounds the distance (in dock-move distance) between optimal solutions in consecutive phases; while we prove such a result for Problem \ref{opt:adapted} (Lemma \ref{lem:scaling_bound}), the running time bound in Theorem \ref{thm:scaling} relies on a bound from \cite{shioura2018m} that is smaller but less general (as it applies only to \ref{opt:no_ops} but not to \ref{opt:adapted}). We then extend our scaling algorithm (see Algorithm \ref{alg:scaling_ops}) to also work for {Problem \ref{opt:adapted}}, and use our own proximity bound to prove the following theorem.}
{
\begin{theorem}\label{thm:scaling_ops}
Algorithm \ref{alg:scaling_ops} runs in time polynomial in $n$ and $\log(B+D)$.
\end{theorem}
We remark that a significantly faster scaling algorithm for \ref{opt:adapted} has been developed in the concurrent work of \cite{shioura2018m}}.


\section{Case Studies}\label{sec: case_study}



%

In this section we present the results of case studies based on data from three different bike-sharing systems: Citi Bike in NYC, Blue Bikes in Boston, and Divvy in Chicago. Some of our results are based on an extension of the user dissatisfaction function which we first define in Section \ref{sec: long_run}. Thereafter, in Section \ref{ssec: data} we describe the data sets underlying our computation. Finally, in Section \ref{ssec: objective} we describe the insights obtained from our analysis. {While some of the results presented in this section are based on proprietary data, we discuss in our electronic supplement \ref{sec:ec_online} which can be reproduced using a data set and the source code we are making public.}

\vspace{-.05in}

\subsection{Long-Run-Average Cost}\label{sec: long_run}
A topic that has come up repeatedly in discussions with operators of {bike-sharing systems} is the fact that their means to rebalance overnight do not usually suffice to begin the day with the bike-optimal allocation. In some cities, like Boston, no rebalancing at all happens overnight. As such, it is desirable to optimize for reallocations that are robust with respect to the amount of overnight rebalancing. To capture such an objective, we define the long-run average of the user dissatisfaction function. Rather than mapping an initial condition in bikes and empty docks to the expected number of out-of-stock events over the course of one day, the long-run average maps to the average number of out-of-stock events over the course of infinitely many days. Notice that (under a weak ergodicity assumption discussed below) in this model the initial allocation of bikes is irrelevant, and so is the total number of bikes allocated, as the long-run distribution of bikes present is determined solely by the distribution of arrival sequences. Formally, denoting by $X\oplus Y$ the concatenation of arrival sequences~$X$ and $Y$, i.e., $(X_1,\ldots, X_t, Y_1,\ldots, Y_s)$, we define the long-run average of a station $i$ with demand profile $p_i$ as follows.
\begin{definition}
The long-run-average of the user dissatisfaction function at station $i$ with demand profile $p_i$ is 
$$
c^\pi_i(d,b) = \lim_{T\to\infty}\frac{\mathbb{E}_{Y_j\sim p_i (i.i.d.)} [c^{Y_1\oplus Y_2\oplus\cdots\oplus Y_T}(d,b)]}{T}.
$$
\end{definition}
We can compute $c{_i}^\pi(d,b)$ by computing for a given demand profile $p_i$ the transition probabilities $\rho_{xy}:=\sum_X p_i(X) \mathbf{1}_{\delta_X(d_i+b_i-x,x)=y}$, that is the probability of station $i$ having $y$ bikes at the end of a day, given that it had $x$ at the beginning, and given that each sequence of arrivals $X$ occurs with probability $p_i(X)$. Given the resulting transition probabilities, we define a discrete Markov chain on $\{0,\ldots,d_i+b_i\}$ and denote by $\pi_{p_i}^{d_i+b_i}$ its stationary distribution. This permits us to compute $c{_i}^\pi(d,b)=\sum_{k=0}^{d+b}\pi_{p_i}^{d+b}(k)c_i(d+b-k,k)$. Furthermore, from the definition of $c{_i}^\pi(\cdot,\cdot)$ it is immediately clear that $c{_i}^\pi(\cdot,\cdot)$ is also multimodular; as such all results proven in the previous sections about $c{_i}(\cdot,\cdot)$ also extend to $c{_i}^\pi(\cdot,\cdot)$. In addition, we observe that, as long as the discrete Markov chain with transition probabilities $\rho_{xy}$ is ergodic (e.g., with demand based on non-zero Poisson rates for both bike rentals and returns), $c{_i}^\pi(\cdot,\cdot)$ depends only on the sum of its two arguments but not on the value of each (as the initial number of bikes does not influence the steady-state number of bikes). Before comparing the results of optimizing over $c{_i}^\pi(\cdot,\cdot)$ and over $c{_i}(\cdot,\cdot)$, we now give some intuition for why the long-run average provides a contrasting regime.

\subsubsection*{Intuition for the Long-run Average. } It is instructive to consider an example to illustrate where optimizing over the long-run average deviates from optimizing over a single day. 
{Consider two stations $i$ and $j$:
 at station $i$ demand consists, determistically, of $k$ rentals every day;
at station~$j$, with probability $P<<\frac{1}{2}$, there are $k$ rentals followed by $k$ returns, and with probability~$1-P$ there are no rentals at all.
At station $i$ the user dissatisfaction function decreases by 1 for each of the first $k$ full docks added; however, its long-run average objective remains constant at $k$: No matter how many docks and bikes are added, in the long-run the station is empty at the beginning of the day and therefore all $k$ customers experience out-of-stock events. At station $j$, the first $k$ full docks added only decrease the user dissatisfaction function by $2P<<1$ each, but the long-run average is also decreased by $2P$ for each dock added. Thus, optimally placing $k$ docks and bikes at the two stations yields fundamentally different solutions depending on whether we optimize for one objective or the other. Furthermore, optimizing for the long-run average only gives a fraction $2P$ of the optimal improvement for a single day, while optimizing for a single day gives no improvement at all for the long-run average objective. 
}
{Two lessons can be derived from this example. First,  optimizing over one regime can, in theory, return solutions that are very bad in the other.
Second, stations at which demand is antipodal (rentals in the morning, returns in the afternoon or vice-versa) make better use of additional capacity in the long-run average regime.}



\subsection{Data Sets}\label{ssec: data}

We use data sets from the bike-sharing systems of three major American cities to investigate the effect different allocations of docks might have in each city. The three cities, New York City, Boston, and Chicago, vary widely in the sizes of their systems. When the data was collected from each system's open data feed (summer 2016), Boston had 1300 bikes and 2805 docks across 158 stations, Chicago had~4700 bikes and~9987 docks across~581 stations, and NYC had 6750 bikes and 15274 docks across 455 stations (given that the feeds only provide the number of bikes in each station, they do not necessarily capture the entire fleet size, e.g., in New York City a significant number of bikes is kept in depots over night).

For each station (in each system), we compute piecewise constant Poisson arrival rates to inform our demand profiles. To be precise, we take all weekday rentals/returns in the month of June 2016, bucket them in the 30-minute interval of the day in which they occur, and divide the number of rentals/returns at each station within each half-hour interval by the number of minutes during which the station was nonempty/nonfull. We compute the user dissatisfaction functions assuming that the demand profiles stem from these Poisson arrivals (see \citealt{henomashm16} and \citealt{parikh2014estimation}). Some of our results in this section rely on the same procedure with data collected from other months. 

Given that (in practice) we do not usually know the lower and upper bounds on the size of each station, we set the lower bound to be the current minimum capacity within the system and the upper bound to be the maximum one.  Furthermore, we assume that $D+B$ is equal to the current allocated capacity in the system, i.e., we only reallocate existing docks.




\subsection{Impact on Objective}\label{ssec: objective}


We summarize our results in Table \ref{tab:summary}. The columns Present, OPT, and 150-moved compare the objective with (i) the allocation before any docks are moved, (ii) the optimal allocation of bikes and docks, and (iii) the best allocation of bikes and docks that can be achieved by moving at most 150 docks from the current allocation. The columns headed $c$ contain the bike-optimal objective for a given allocation of docks, the columns headed $c^\pi$ the long-run-average objective (for the same dock allocation). Two interesting observations can be made. First, though the optimizations are done over bike-optimal allocations without regard to the long-run average, the latter improves greatly in all cases. Second, in each of the cities, moving 150 docks yields a large portion of the total possible improvement. This stands in contrast to the large number of moves needed to find the actual optimum (displayed in the column Moves to OPT) and is due to diminishing returns of the moves. 

\OneAndAHalfSpacedXI
\begin{table}
    \centering
    \begin{tabular}{|c||r|r|r|r|r|r|r|}
    \hline
    & \multicolumn{2}{c}{Present} \vline
    & \multicolumn{2}{c} {OPT} \vline
    & \multicolumn{2}{c}{150-moved}
    & \multicolumn{1}{|c|}{Moves to OPT}\\
    \hline
    City & $c$ & $c^\pi$
        & $c$ & $c^\pi$
        & $c$ & $c^\pi$ & \\
    \hline
    Boston  & 831 & 1,118 & 607 & 945  & 672 & 985 & 412   \\
    Chicago & 1,462 & 2,340 & 763 & 1,847  & 1,224 & 2,123 & 1,556 \\
    NYC & 8,251 & 10,937 & 6,499 & 9,232  & 7,954 & 10,643 & 2,821  \\
    \hline
    \end{tabular}
    \caption{Summary of main computational results with $c$ denoting bike-optimal, $c^\pi$ the long-run-average cost {-- obtained from Algorithm \ref{alg:gradient_descent} applied to data sets from June 2016}.}
    \label{tab:summary}
\end{table}
\DoubleSpacedXI

\begin{figure}
    \centering
\includegraphics[height=.25\textwidth]{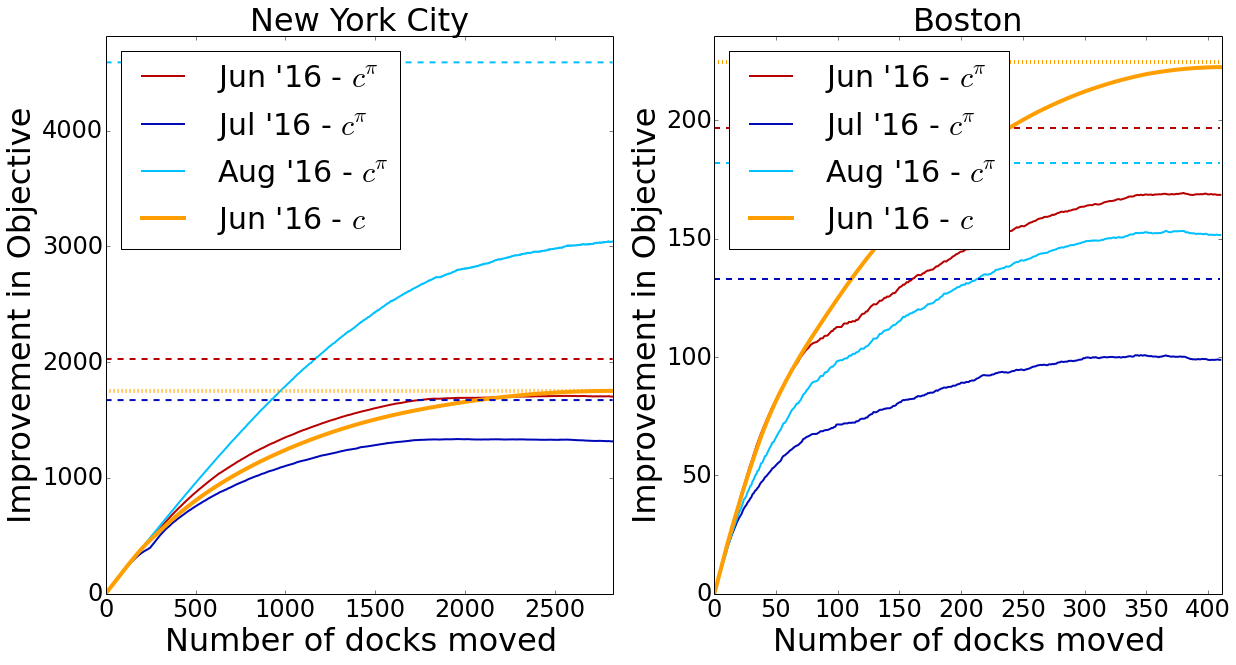}
    \caption{{Improvement in objectives (evaluated with different data sets and $c^\pi$) for moves to the optimal allocation for June '16 data and objective $c$. Dashed lines denote the respective optimal objective.}}
    \label{fig:kmoves}
\end{figure}

A more complete picture of these insights is given in Figure \ref{fig:kmoves}. The $x$-axis shows the number of docks moved starting from the present allocation, the $y$-axis shows the cumulative \emph{improvement in objective}, i.e., the difference between the initial objective and the objective after moving $x$ docks. Each of the solid lines corresponds to different demand estimates being used to evaluate the same allocation of docks. The dotted lines (in the same colors) represent the maximum improvement, for each of the demand estimates, that can be achieved by reallocating docks; while these are not achieved through the dock moves suggested by the estimates based on June 2016 data, significant improvement is made towards them in every case. In particular, the initial moves yield approximately the same improvement for the different objectives/demand estimates.
Thereafter, the various improvements diverge, especially for the NYC data from August 2016. This may be partially due to the system expansion in NYC that occurred in the summer of 2016, but does not contradict that all allocations corresponding to values on the $x$-axis are optimal in the sense of Theorem \ref{thm: koptimal}.


\OneAndAHalfSpacedXI
\DoubleSpacedXI

\subsubsection*{Seasonal Effects. } As we mentioned in Section \ref{sec: model} we also consider the impact of seasonal effects. In Table \ref{table: seasons} we show the improvement in objective when optimizing the movement of 200 docks in New York City based on demand estimates in June 2016 and evaluate the objective with the long-run average based on demand estimates based on March and November 2017. The estimated improvements suggest that optimizing with respect to June yields notable improvement with respect to any other.


\OneAndAHalfSpacedXI
\addtolength{\tabcolsep}{3pt}    
\begin{table}[ht]
\centering
\begin{tabular}{c|r|r|r|}
\cline{2-4} 
                                    & June 2016 & March 2017 & November 2017 \\ \hline
\multicolumn{1}{|c|}{New York City} &  358.7         &        260.3    &       294.6        \\ \hline
\end{tabular}
\caption{{Improvement of 200 docks moved based on the long-run average objective $c^{\pi}$ evaluated with demand estimates from June 2016, evaluated with the long-run average objective $c^{\pi}$ and demand estimates from 2017.}}
\label{table: seasons}
\end{table}
\DoubleSpacedXI

\subsubsection*{Operational Considerations. }
It is worth comparing the estimated improvement realized through reallocating docks to the estimated improvement realized through current rebalancing efforts. {According to its monthly report \citep{citibikereport}, Citi Bike rebalanced an average 3,452 bikes per day in June 2016:
this number counts the average number of \emph{rebalancing actions}, meaning that each pickup/dropoff counts as one bike rebalanced. A simple coupling argument implies that a single pickup/dropoff} yields at most a change of 1 in the user dissatisfaction function (see Figure~\ref{fig:UDF2D}); thus, rebalancing reduced out-of-stock events by \emph{at most} 3,452 per day (assuming that each rebalanced bike actually has that much impact is extremely optimistic). Contrasting that to the estimated impact of strategically moving, for example, 500 docks diminishes the estimated number of out-of-stock events \emph{by more than a fifth} of Citi Bike's (daily) rebalancing efforts.

Second, discussions with operators uncovered an additional operational constraint that can arise due to the physical design of the docks. Since these usually come in triples or quadruples, the exact moves suggested may not be feasible; e.g., it may be necessary to move docks in multiples of 4. By running the scaling algorithm {(see Algorithm \ref{alg:scaling} in Appendix \ref{appendix:pseudocode}) only with $k\geq2$}, we can find an allocation in which docks are only moved in multiples of 4. With that allocation, the objective of the bike-optimal allocation is 640, 848, and 6573 in Boston, Chicago, and NYC respectively, suggesting that despite this additional constraint, when compared to the column headed by~OPT and~$c$, almost all of the improvements can be realized.

\DoubleSpacedXI



\section{A Posteriori Evaluation of Impact}\label{sec: impact}


{In this section, we show how one can use the user dissatisfaction functions to estimate (after the fact) the impact of reallocated capacity, and apply this approach to the 6 stations that were part of the pilot program mentioned in Section \ref{ssec:contribution}. 
One way to evaluate the impact} would be to estimate new demand rates after docks have been reallocated, compute new user dissatisfaction functions for stations with added (decreased) capacity, and evaluate for those stations and the new demand rates the decrease (increase) between the old and the new number of docks. A drawback of such an approach is the heavy reliance on the assumed underlying stochastic process. Instead, we present here a data-driven approach with only little reliance on estimated underlying demand profiles.

Throughout this section, we denote by $d$ and $b$ the number of empty docks and bikes at a station after docks were reallocated, whereas $d'$ and $b'$ denote the respective numbers before docks were reallocated. Notice that while $d+b$ and $d'+b'$ are known (capacity before and after docks were moved) and $b$ can be found on any given morning (number of bikes in the station at 6AM), we rely on some assumed value for $b'$ --- for that, in our implementation, we picked both $\min\{d'+b', b\}$ and $b\times\frac{(d'+b')}{d+b}$, that is, either the same number of bikes (unless that would be larger than the old capacity before docks were added) or the same proportion of docks filled with bikes. 

\vspace{-.1in}

\subsection{Arrivals at Stations with Increased Capacity}\label{ssec:impact_regimes} In earlier sections, we assumed a known distribution for the sequence of arrivals based on which we compute the user dissatisfaction functions. In contrast, in this section we rely exclusively on observed arrivals (without any assumed knowledge of the underlying stochastic process) to analyze stations with increased capacity. This is motivated by a coupling argument to justify that censoring need not be taken care of explicitly in this case. To formalize our argument, we need to introduce some additional notation for the arrival sequences. Recall from Section \ref{sec: model} that we denoted by $X=(X_1,\ldots, X_S)$ a sequence of customers arriving at a bike-sharing station to either rent or return a bike and that $X$ included failed rentals and returns, which in practice would not be observed because they are censored. Which $X_i$ are censored depends on the (initial) number of bikes and docks at the station. Let us denote by $X(d,b)$ the subsequence of $X$ that only includes those customers whose rentals/returns are successful (not censored) at a station that is initialized with $d$ empty docks and $b$ bikes, i.e., the ones who do not experience out-of-stock events. Given the notation $c^X(\cdot,\cdot)$ used in Section \ref{sec: model} for a particular sequence of arrivals, we can then compute $c^{X(d,b)}(\cdot,\cdot)$. In particular, denoting by $d',b'$ the number of empty docks and bikes without the added capacity, we may compute $c^{X(d,b)}(d',b')$. The following proposition then motivates the notion that censoring may be ignored at stations with added capacity.

\begin{proposition}\label{prop:estimate}
For any $X$, $d'\leq d$ and $b'\leq b$, we have 


$$
c^{X}(d',b')-c^{X}(d,b)
=
c^{X(d,b)}(d',b')-c^{X(d,b)}(d,b)=c^{X(d,b)}(d',b').
$$
\end{proposition}


\emph{Proof. }
 The proof of the second equality follows immediately from $X(d,b)$ including exactly those customers among $X$ that are not censored, when a station is initialized with $d$ empty docks and $b$ bikes, so $c^{X(d,b)}(d,b)=0$. Now, on the LHS, we can inductively go through all customers among $X$ that are out-of-stock events when the station is initialized with $d$ empty docks and $b$ bikes. Since $d\geq d'$ and $b\geq b'$, each one of those increases both terms in the difference by 1. Thus, taking them out of $X$ does not affect the value of the difference. But then, we are left with only $X(d,b)$. 


\subsubsection*{Extension to Stations with Decreased Capacity. } Proposition \ref{prop:estimate} does not apply to stations with decreased capacity: suppose $d<d'$ and $b=b'$; once the station (initialized with $d$ empty docks and $b$ bikes) becomes full, $X(d,b)$ observes no further returns even though these would be part of $X(d',b')$. To account for out-of-stock events occurring in that way, we fill in the censored periods with demand estimates. This does not usually require knowledge of the full demand-profile; for example, for a station that is non-empty and non-full over the course of the day, no estimates are needed at all. Further, for periods of time in which the station is full, we only need to estimate the number of intended returns -- rentals over that period of time would not be censored.

\subsubsection*{Extension to Rebalancing.} Based on our reasoning in Section \ref{sec: assumptions}, our analysis of the user dissatisfaction functions and the resulting optimization problem (see Sections \ref{sec: model} and \ref{sec: algorithm}) so far did not consider the rebalancing of bikes. In contrast, in the a posteriori analysis, we are able to take rebalancing into account.

To simplify the exposition, we restrict ourselves here to rebalancing that adds bikes to a station, though the reasoning extends to rebalancing that removes bikes. The simplest approach to treat bikes added through rebalancing is to just treat them as returns and thus include them (as virtual customers) in the sequence of arrivals $X$. However, this may cause an unreasonable increase to the value of $c^{X(d,b)}(d',b')$ (when the number of bikes added is greater than the number of empty docks would have been at that point in time if the station had initially had $d'$ empty docks). In that case, the virtual customers (corresponding to rebalanced bikes) would incur out-of-stock events and thereby increase the value of the user dissatisfaction function. A more optimistic method that also treats rebalanced bikes as virtual customers is to redefine the user dissatisfaction function in such a way that out-of-stock events are only incurred by returns that correspond to non-rebalanced bikes. This, in essence, decouples the user dissatisfaction functions into subsequences, each of which are evaluated independently. Our analysis below applied the latter, more optimistic method.


\subsection{Impact of Initial Dock Moves in NYC}
\begin{figure}[ht]
    \centering
\includegraphics[height=.25\textwidth]{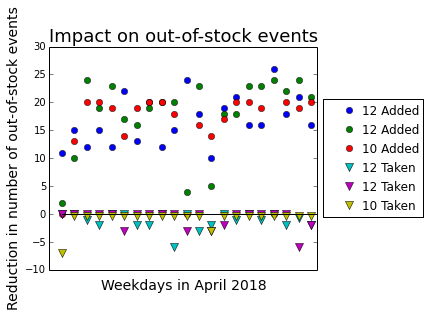}
    \caption{Evaluation of impact at stations with increased and decreased capacity.}
    \label{fig:evaluation}
\end{figure}

We consider 3 stations at which capacity was increased and 3 stations at which it was decreased based on our recommendations. For two of the stations at which capacity was increased, 12 docks were added, for one of them the capacity was increased by 10; the decreases were by the same amounts, so in total this involved reallocating 34 docks. In Figure \ref{fig:evaluation} we present the estimated impact for each weekday in April 2018 (without the extension to rebalancing). For stations with added capacity we set $d$ and $b$ according to the number of bikes at 6AM. We evaluated $c^{X(d,b)}(d',b')$ for stations with docks added (see Proposition \ref{prop:estimate}) using the observed arrivals $X(d,b)$ for each day. For the stations with docks taken away we estimated $X$ by assuming a Poisson number of rentals (returns) whenever the station was empty (full), where the rate is based on decensored estimated demand from the same month. We use that to compute $c^{X(d,b)}(d',b')-c^{X(d,b)}(d,b)$ for these stations. The resulting values for different implementations are summarized in Table \ref{table: impact}; aggregated over the entire month, the estimated net reduction in out-of-stock events varies between 831 and 1121, i.e., about 1.2--1.6 fewer dissatisfied users per day and dock moved. Translating this into reduced rebalancing costs, and comparing it to the cost of reallocating docks, strategically reallocating docks amortizes (depending on some system idiosyncrasies) in 2--5 weeks.





\OneAndAHalfSpacedXI
\begin{table}[ht]
\centering
\begin{tabular}{c|r|r|r|r|}
\cline{2-5}
                                                        & \multicolumn{2}{c|}{No Rebalancing}                   & \multicolumn{2}{c|}{Rebalancing}                \\ \cline{2-5} 
                                                        & $\min\{b, d'+b'\}$    & $b\times\big(\frac{d'+b'}{d+b}\big)$ & $\min\{b, d'+b'\}$    & $b\times\big(\frac{d'+b'}{d+b}\big)$ \\ \hline
\multicolumn{1}{|l|}{Decrease where capacity was added} &                831.0       &         1121.0                   &      882.0                 &           1027.0                 \\ \hline
\multicolumn{1}{|l|}{Increase where capacity was taken} &                  0     &              58.7              &              0         &     59.7                       \\ \hline
\multicolumn{1}{|l|}{Net Reduction}                     & \multicolumn{1}{r|}{831.0}  & \multicolumn{1}{r|}{1062.3}      & \multicolumn{1}{r|}{882.0} & \multicolumn{1}{r|}{967.3}      \\ \hline
\end{tabular}
\caption{Estimated cumulative changes at stations affected by dock reallocations based on the different evaluations described in Section \ref{ssec:impact_regimes}.}
\label{table: impact}
\end{table}
\DoubleSpacedXI

\section{Conclusion}\label{sec: conclusion}

We have considered several models that capture central questions in the design of dock-based bike-sharing systems, as are currently prevalent in North America. These models gave rise to new algorithmic discrete optimization questions, and we have demonstrated that they have sufficient mathematical structure to permit their efficient solution, thereby also extending existing theory in discrete convexity. We have focused on the (re-)allocation of docks throughout the footprint of a bike-sharing system, capturing aspects of both better positioning of existing docks, and the optimal augmentation of an existing system with additional docks. These algorithms and models have been employed by systems within the United States with the desired effect of improving their day-to-day performance.
 
An alternative to optimizing dock allocations is to abandon the need to do so at all, by means of adopting a so-called dockless system. This approach has become prevalent in China, and is gradually being implemented in North America on a much smaller scale (both in comparison to the systems in China, and to the dock-based systems in North America); the management of these systems has its own challenges, and it remains to be seen whether these challenges can be overcome.  Hybrid systems in which differential pricing enables centralized docking/parking areas that work in concert with dockless bikes may provide another path forward, as is done, for example, in Portland's Biketown system. Extensions of the methods we developed here will likely see continued use in this new setting as well.

\newpage
\bibliography{bib}
\newpage

\appendix
\section{Pseudocode of Algorithms}\label{appendix:pseudocode}
In this appendix we provide pseudocode for each of the algorithms described in the main body of the text and add some more details about their efficient implementation. Notice that Algorithm \ref{alg:bike_optimal} can easily be adapted with scaling techniques to run in $O(\log(B))$ phases with each one requiring at most $n$ iterations.


\SingleSpacedXI
\begin{algorithm}
\caption{Greedy Routine to Find Bike-Optimal Allocation}\label{alg:bike_optimal}
\begin{algorithmic}[1]
\Require UDFs $c_i(\cdot,\cdot)$, current allocation $(\vec{d},\vec{b})$, integer $\alpha$
\State $Flag \gets True$
\While{$Flag$}
\State $i^\star = \arg\max_i \{c_i(d_i,b_i)-c_i(d_i-\alpha,b_i+\alpha)\}$ 
\State $j^\star = \arg\max_j \{c_j(d_j,b_j)-c_j(d_j+\alpha,b_j- \alpha)\}$ 
\If{$c_{i^\star}(d_{i^\star},b_{i^\star})+c_{j^\star}(d_{j^\star},b_{j^\star})  > c_{i^\star}(d_{i^\star}-\alpha,b_{i^\star}+\alpha)+ c_{j^\star}(d_{j^\star}+\alpha,b_{j^\star}-\alpha)$} 
\State $d_{i^\star} \gets d_{i^\star}-\alpha$, $b_{i^\star} \gets b_{i^\star}+\alpha$
\State $d_{j^\star} \gets d_{j^\star}+\alpha$, $b_{j^\star} \gets b_{j^\star}-\alpha$
\Else 
\State $Flag\gets False$ 
\EndIf
\EndWhile
\State \textbf{return} $(\vec{d},\vec{b})$\Comment{If $\alpha=1$, then $(\vec{d},\vec{b})$ is bike-optimal for the current allocation of docks $(\vec{\bar{d}},\vec{\bar{b}})$}
\end{algorithmic}
\end{algorithm}

\DoubleSpacedXI
For each of our algorithms to find the optimal reallocation of docks we require the following notion of a gradient-descent step.

\SingleSpacedXI

\begin{algorithm}
\caption{Gradient-descent Step}\label{alg:gradient_descent_step}
\begin{algorithmic}[1]
\Require UDFs $c_i(\cdot,\cdot) \forall i$, allocation $(\vec{d},\vec{b})$, and neighborhood-definition $N(\cdot,\cdot)$ 
\State $(\vec{d}',\vec{b}') = {\arg\min_{(\vec{d}',\vec{b}')\in N(\vec{d},\vec{b})} c(\vec{d}',\vec{b})}$
\If{$c(\vec{d}',\vec{b}')<c(\vec{d},\vec{b})$} 
\State \textbf{ return $(\vec{d}',\vec{b}')$} \Comment{$(\vec{d}',\vec{b}')$ is the best solution in the neighborhood of $(\vec{{d}},\vec{{b}})$}
\Else
\State \textbf{ return } $(\vec{{d}},\vec{{b}})$ \Comment{$(\vec{{d}},{\vec{b}})$ is at least as good as any other solution in its neighborhood}
\EndIf
\end{algorithmic}
\end{algorithm}

\DoubleSpacedXI

Algorithm \ref{alg:gradient_descent_step} requires taking a minimum over a neighborhood of size ${\Theta(n^3)}$. This is because two of the four moves include three stations and we may pick any three stations to carry out these moves. However, the four different moves can only yield six different changes at each station:  an added bike, a removed bike, an added empty dock, a removed empty dock, an added full dock, or a removed full dock. By precomputing in a sorted heap, in $\Theta(n\log(n))$ steps, the six possible ways in which the objective at each station can be affected, we can then find the minimum in Algorithm \ref{alg:gradient_descent_step} in time $O(1)$. Rather than computing the heaps in each iteration, we update only the values in the heaps that correspond to the stations involved in the dock-move in Algorithm \ref{alg:gradient_descent_step}. Precomputing then requires $O(n\log(n))$ whereas each iteration thereafter runs in $O(\log(n))$.

\SingleSpacedXI

\begin{algorithm}
\caption{Gradient-descent algorithm}\label{alg:gradient_descent}
\begin{algorithmic}[1]
\Require UDFs $c_i(\cdot,\cdot)$ and current allocation $(\vec{\bar{d}},\vec{\bar{b}})$, constraint~$z$
\State $(\vec{d}^0,\vec{b}^0) \gets $ Algorithm \ref{alg:bike_optimal} applied to the current allocation $(\vec{\bar{d}},\vec{\bar{b}})$ with $\alpha=1$
\For{$r$ in $\{1,\ldots, z\}$}
\State $(\vec{d}^{r},\vec{b}^{r}) \gets$ Algorithm \ref{alg:gradient_descent_step} applied to $(\vec{d}^{r-1},\vec{b}^{r-1})$ with $N(\cdot,\cdot)$ as in Definition \ref{def:neighborhood}
\State \textbf{if} $(\vec{d}^r,\vec{b}^r)=(\vec{d}^{r-1},\vec{b}^{r-1})$ \textbf{then} \textbf{return} $(\vec{d}^r,\vec{b}^r)$ \Comment{$(\vec{d}^{r},\vec{b}^{r})$ is optimal with regard to $N(\cdot,\cdot)$}
\EndFor
\State \textbf{return} $(\vec{d}^z,\vec{b}^z)$
\end{algorithmic}
\end{algorithm}

\DoubleSpacedXI

Below we define the scaling algorithms described in Section \ref{sec: scaling}. The analysis of these algorithms is found in Appendix \ref{sec:appendix_proof_scaling}. We define $N^\alpha(\vec{d},\vec{b})$ as the set of allocations that can be obtained by applying the same dock-move $\alpha$ times to $(\vec{d},\vec{b})$, i.e.,
$$
N^\alpha(\vec{d},\vec{b}):=\{o_{ij}^\alpha(\vec{d},\vec{b}), e_{ij}^\alpha(\vec{d},\vec{b}), E_{ijh}^\alpha(\vec{d},\vec{b}), O_{ijh}^\alpha(\vec{d},\vec{b}): i,j,h\in[n]\},
$$
where $f^\alpha(\vec{d},\vec{b}) = f(f^{\alpha-1}(\vec{d},\vec{b}))$ and $f^0(\vec{d},\vec{b})=(\vec{d},\vec{b})$.
\SingleSpacedXI

\begin{algorithm}
\caption{Scaling Algorithm without Operational Constraints}\label{alg:scaling}
\begin{algorithmic}[1]
\Require UDFs $c_i(\cdot,\cdot)$, $(\vec{\bar{d}},\vec{\bar{b}})$
\State $\vec{d}\gets \vec{\bar{d}},\; \vec{b} \gets  \vec{\bar{b}}$
\For{$k$ in $\{0, 1,2, \ldots, \lfloor\log_2(B+D) \rfloor, \lfloor\log_2(B+D) \rfloor+1\}$}
\State $\alpha_k \gets 2^{\lfloor\log_2(B+D) \rfloor+1-k}$
\State $(\vec{d},\vec{b}) \gets $ Algorithm 1 applied to the current allocation $(\vec{d},\vec{b})$ with $\alpha=\alpha_k$
\While{The gradient-descent step yields improvement}
\State $(\vec{d},\vec{b}) \gets$ Gradient-descent step applied to $(\vec{d},\vec{b})$ with $N^{\alpha_k}(\cdot,\cdot)$
\EndWhile
\EndFor
\State \textbf{return} $(\vec{d},\vec{b})$
\end{algorithmic}
\end{algorithm}

\DoubleSpacedXI

In practice, the scaling algorithm need not use $k>6$ because we tend to have $u_i-l_i<2^6=64\forall i$; further, since initializing each phase requires some overhead in recomputing the binary heaps, we can then define the following version of the scaling algorithm (see Algorithm \ref{alg:hybrid}) that skips some powers of 2.

\SingleSpacedXI

\begin{algorithm}
\caption{Hybrid Algorithm without Operational Constraints}\label{alg:hybrid}
\begin{algorithmic}[1]
\Require UDFs $c_i(\cdot,\cdot)$, $(\vec{\bar{d}},\vec{\bar{b}})$
\State $\vec{d}\gets \vec{\bar{d}},\; \vec{b} \gets  \vec{\bar{b}}$
\For{$k$ in $\{0, 2,3\}$}
\State $\alpha_k \gets 2^k$
\State $(\vec{d},\vec{b}) \gets $ Algorithm 1 applied to the current allocation $(\vec{d},\vec{b})$ with $\alpha=\alpha_k$
\While{The gradient-descent step yields improvement}
\State $(\vec{d},\vec{b}) \gets$ Gradient-descent step applied to $(\vec{d},\vec{b})$ with $N^{\alpha_k}(\cdot,\cdot)$
\EndWhile
\EndFor
\State \textbf{return} $(\vec{d},\vec{b})$
\end{algorithmic}
\end{algorithm}

\begin{algorithm}
\caption{Scaling Algorithm with Operational Constraints}\label{alg:scaling_ops}
\begin{algorithmic}[1]
\Require UDFs $c_i(\cdot,\cdot)$, $(\vec{\bar{d}},\vec{\bar{b}})$, constraint~$z$
\State $\vec{d}\gets \vec{\bar{d}},\; \vec{b} \gets  \vec{\bar{b}}$
\For{$k$ in $\{0, 1,2, \ldots \lceil\log(D+B)\rceil\}$}
\State $\alpha_k \gets 2^{\lfloor\log_2(B+D) \rfloor+1-k}$
\State $M\gets 16n^6(n+4)\alpha_k$
\State $(\vec{d},\vec{b}) \gets $ Algorithm 1 applied to the current allocation $(\vec{d},\vec{b})$ with $\alpha=\alpha_k$
\State $(\vec{d},\vec{b})\gets$ output of Algorithm \ref{alg:find_d_b} with input $(\vec{\bar{d}},\vec{\bar{b}})$, $(\vec{d},\vec{b})$, $M$, $\alpha_k$
\While{The gradient-descent step yields a solution that (i) has lower objective and (ii) fulfills the operational constraints}
\State $(\vec{d},\vec{b}) \gets$ Gradient-descent step applied to $(\vec{d},\vec{b})$ with $N^{\alpha_k}(\cdot,\cdot)$
\EndWhile
\EndFor
\State \textbf{return} $(\vec{d},\vec{b})$
\end{algorithmic}
\end{algorithm}

\begin{algorithm}
\caption{Subroutine in-between phases}\label{alg:find_d_b}
\begin{algorithmic}[1]
\Require Current capacity $\bar{d}_i+\bar{b}_i$ for each station $i$, allocation $(\vec{d}^k,\vec{b}^k)$, integers $M$ and $\alpha_k$
\State $\mathcal{S}^+\gets\{i:d_i^k+b_i^k\geq\bar{d}_i+\bar{b}_i\}, \mathcal{S}^-\gets\{j:d_j^k+b_j^k<\bar{d}_j+\bar{b}_j\}$
\State $\vec{d}\gets\vec{d}^k, \vec{b}\gets\vec{b}^k$
\While{there exists $i\in \mathcal{S}^+$ with $d_i+b_i>\max\{d_i^k+b_i^k-M,\bar{d}_i+\bar{b}_i\}$ 
}
\If{there exists $j\in \mathcal{S}^-$ with $d_j+b_j<\min\{d_j^k+b_j^k+M,\bar{d}_j+\bar{b}_j\}$.}
\State Update $(\vec{d},\vec{b})$ by carrying out a best-possible dock-move of $\alpha_k$ docks and bikes from such $i$ to such $j$
\EndIf
\State \textbf{else} pick $j\in \mathcal{S}^-$ with $d_j+b_j<\min\{d_j^k+b_j^k+nM,\bar{d}_j+\bar{b}_j\}$, and update $(\vec{d},\vec{b})$ by carrying out a best-possible dock-move of $\alpha_k$ docks and bikes from such $i$ to $j$
\EndWhile
\While{there exists $j\in \mathcal{S}^-$ with $d_j+b_j<\min\{d_j^k+b_j^k+M,\bar{d}_j+\bar{b}_j\}$.
}
\State pick $i\in \mathcal{S}^+$ with $d_i+b_i>\max\{d_i^k+b_i^k-nM,\bar{d}_i+\bar{b}_i\}$, and update $(\vec{d},\vec{b})$ by carrying out a best-possible dock-move of $\alpha_k$ docks and bikes from $i$ to such $j$
\EndWhile
\State \textbf{return} $(\vec{d},\vec{b})$
\end{algorithmic}
\end{algorithm}


\DoubleSpacedXI


\section{Omitted Proofs for Section \ref{sec: model}}

\subsection{Proof of Lemma \ref{lemma: multimod}}
\label{sec:appendix_lemma_multimod}

\subsubsection*{Recursive UDF definition.} Throughout this proof we denote by $X(t)$ the subsequence of $X=(X_1,\ldots,X_s)$ that only includes the first~$t$ customers, i.e., $X(t)=(X_1,\ldots,X_t)$, and by $\delta_X(d,b)$ the number of open docks at a station with an initial allocation of $d$ empty docks and $b$ bikes (full docks) after the arrival of sequence $X$. Then we can derive $\delta_{X(t)}(d,b)$ inductively via the equations $\delta_{X(0)}(d,b)=d$ and
$$\delta_{X(t)}(d,b):=\max\{0,\min\{d+b,\delta_{X(t-1)}(d,b)-X_t\}\}.$$

\vspace{-.1in}

Recall that, since arrivals of customers do not affect the total number of docks (empty and full) at the station, the number of bikes after the arrival of the $t$th customer $X_t$ is the number of docks that are not open, i.e., $d+b-\delta_{X(t)}(d,b)$.
Our objective is based on the number of out-of-stock events. In accordance with the arrival model described in Section \ref{sec: model}, customer $t$ experiences an out-of-stock event if and only if she fails to either rent or return a bike, i.e., if and only if $\delta_{X(t)}(d,b)=\delta_{X(t-1)}(d,b)$. This, in turn, occurs if and only if $X_t=1$ and $\delta_{X(t-1)}(d,b)=0$ or $X_t=-1$ and $\delta_{X(t-1)}(d,b)=d+b$. 
As a function of the initial number of open docks and available bikes, we can thus write our cost-function (for $d,b$ yielding finite values) as
$$c^{X}(d,b)=|\{\tau: X_\tau=1, \delta_{X(\tau-1)}(d,b)=0 \}|+|\{\tau: X_\tau=-1, \delta_{X(\tau-1)}(d,b)=d+b \}|.$$
{Notice that in the weighted case, in which a stock-out at an empty station is penalized with $\eta\geq0$, the first term is multiplied by $\eta$.}
Equivalently {(in the unweighted case)}, as a recursion, we can write $c^{X(0)}(d,b)=0$ and $$c^{X(t)}(d,b)=c^{X(t-1)}(d,b)+\mathbf{1}_{\{X_t=1,\delta_{X(t-1)}(d,b)=0\}\cup\{X_t=-1,\delta_{X(t-1)}(d,b)=d+b\}}.$$

\begin{figure}[H]
    \centering
    \includegraphics[width=.5\textwidth]{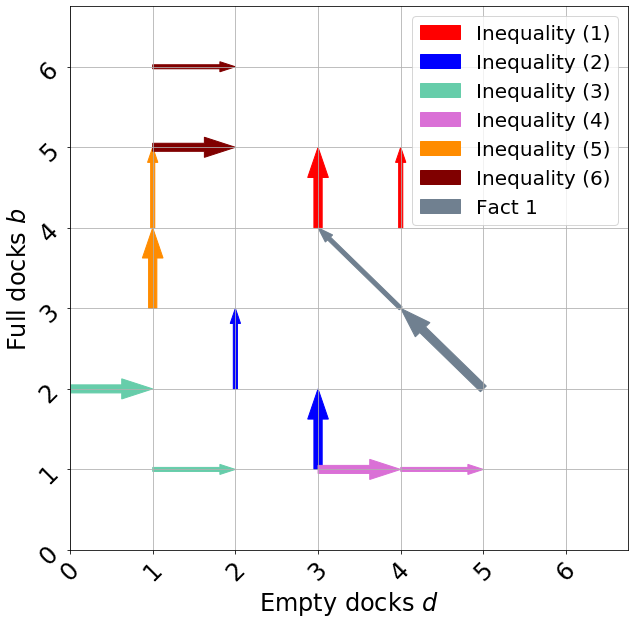}
    \caption{Visual illustration of the inequalities in Definition \ref{def: multimod}, wherein a thicker arrow implies greater improvement, e.g., Inequality (2) implies that the UDF reduction due to an additional full dock at $(3,1)$ is greater-equal to what it would be at $(2,2)$.}
    \label{fig:inequalities}
\end{figure}

Before we prove the inequalities, we refer the reader to a visual illustration of them in Figure~\ref{fig:inequalities}. This should serve as a reminder of the physical interpretation of each inequality, where thicker arrows (of the same color) correspond to greater improvement, e.g., Inequality (1) says that the cost of removing (improvement of adding) a full dock is smaller if an additional empty dock has been added --- this may be easier to see in $f(d+1,b)-f(d+1,b+1)\leq f(d,b)-f(d,b+1)$, i.e., after multiplying (1) by $-1$. Notice also the symmetry in the inequalities in $d$ and $b$, e.g., Inequality~(4) says about adding empty docks what Inequality (5) says about adding full docks.

\subsubsection*{Proof of lemma.}

{We prove the lemma for the weighted case in which the penalty of a stock-out is~$\eta\geq0$ if the station is empty (customer cannot get a bike), and normalized to $1$ when the station is full (customer cannot get a dock).}


We start by considering the cases that, for any of the three inequalities, one of the terms involved evaluates to $\infty$. This can happen (i) when the smallest number of docks allocated, $d+b$ in (1) and $d+b-1$ in (2) and (3), is below $l_i$, (ii) when the largest number allocated, $d+b+2$ in (1), $d+b$ in (2) and (3), is above $u_i$, or (iii) when one of the arguments, $d$ or $b$ in (1), $d-1$ or $b-1$ in (2) and (3), is negative. In all of the cases the inequalities hold vacuously true with the stated conventions that we have $\infty-\infty=\infty$ as well as $\infty\geq x\geq-\infty$ for every $x\in\mathbb{R}\cup\{-\infty, \infty\}$.

For the case where none of the terms evaluate to $\infty$ we prove the inequalities by induction, showing that $c^{X(t)}(\cdot,\cdot)$ is multimodular for all $t$. Since $c_i$ is an expectation over different arrival sequences, the statement of the lemma follows by linearity of expectation. With $t=0$, by definition, $c^{X(t)}(\cdot,\cdot)=0$ and thus there is nothing to show. Suppose that $c^{X(0)}(\cdot,\cdot)$ through $c^{X(t-1)}(\cdot,\cdot)$ are all multimodular. We prove that $c^{X(t)}(\cdot,\cdot)$ is then multimodular as well.

\textbf{Inequality (1).}  We start with $c^{X(t)}(d+1,b+1)-c^{X(t)}(d+1,b)\geq c^{X(t)}(d,b+1)-c^{X(t)}(d,b)$. If 
\vspace{-.1in}
\[\max\{c^{X(1)}(d+1,b+1),c^{X(1)}(d+1,b),c^{X(1)}(d,b+1),c^{X(1)}(d,b)\}=0,\] 

\vspace{-.1in}
\noindent we can use that Inequality (1), by inductive assumption, holds after $t-1$ customers by simply considering $t=1$ as the start of time. Else, we argue below that the inductive assumption on inequalities (4) and~(5) for up to $t-1$ arrivals imply Inequality~(1). 

If $X_1=1$, it must be the case that $d=0$ since otherwise the first customer would be able to return the bike for all of the initial conditions. But then both sides of Inequality (1), evaluated after the first customer, are 0 and $\delta_{X(1)}(d+1,b+1)=0$, $\delta_{X(1)}(d+1,b)=0$, $\delta_{X(1)}(d,b+1)=0$, and $\delta_{X(1)}(d,b)=0$. 
In that case, we may use the inductive assumption on Inequality (5) applied to the remaining $t-1$ customers. 

If instead $X_1=-1$, then it must be the case that $b=0$ since otherwise the first customer would be able to rent a bike for all of the initial conditions. But then both sides of Inequality (1), evaluated after the first customer, are {$-\eta$} and we have $\delta_{X(1)}(d+1,b+1)=d+b+2$, $\delta_{X(1)}(d+1,b)=d+b+1$, $\delta_{X(1)}(d,b+1)=d+b+1$, and $\delta_{X(1)}(d,b)=d+b$, so we may apply Inequality (4) inductively to the remaining $t-1$ customers.

\textbf{Inequality (2).} It remains to prove inequalities (2) and (3) where we may assume $b,d\geq 1$ based on the first paragraph in the proof of the lemma. We restrict ourselves to Inequality (2) as the proof for Inequality~(3) is symmetric with each $X_i$ replaced by $-X_i$ and the coordinates of each term exchanged. Thus, we aim to show $c^{X(t)}(d-1,b+1)-c^{X(t)}(d-1,b)\geq c^{X(t)}(d,b)-c^{X(t)}(d,b-1)$. As before, if 
\[\max\{c^{X(1)}(d-1,b+1),c^{X(1)}(d-1,b), c^{X(1)}(d,b),c^{X(1)}(d,b-1)\}=0,\] 

\noindent the inductive assumption on Inequality (2) applies. If instead $X_1=1$ and the maximum is positive, then the LHS and the RHS of the inequality both evaluate to 0 after the first customer and we have $\delta_{X(1)}(d-1,b+1)=0$, $\delta_{X(1)}(d-1,b)=0$, $\delta_{X(1)}(d,b)=0$, $\delta_{X(1)}(d,b-1)=0$. In that case, both sides of the inequality are subsequently coupled and the inequality holds with equality. 

In contrast, if $X_1=-1$ and the maximum is positive, then $b=1$, the RHS is {$-\eta$}, and the LHS is 0. In this case we have
$\delta_{X(1)}(d-1,b+1)=d$, $\delta_{X(1)}(d-1,b)=d$, $\delta_{X(1)}(d,b)=d+1$, $\delta_{X(1)}(d,b-1)=d$.
Let~$\hat{t}$ denote the index of the next customer that is dissatisfied when the system is initialized with any one of the four initial conditions~$(d-1,b+1), (d-1,b), (d,b)$, or $(d,b-1)$. 

If $X_{\hat{t}}=1$, then the terms that correspond to having $d$ open docks after the first customer all increase by 1. Thus, both terms on the LHS of Inequality (2) increase by 1, whereas only the negative term on the RHS increases. Thereafter, the inequality holds with $0\geq{-1-\eta}$. Moreover, after customer $\hat{t}$, we have $\delta_{X(\hat{t})}(d-1,b+1)=\delta_{X(\hat{t})}(d,b)=0$, and $\delta_{X(\hat{t})}(d-1,b)=\delta_{X(\hat{t})}(d,b-1)=0$. Thus, both sides of the inequality are again coupled.

Finally, if $X_{\hat{t}}=-1$, then both terms on the RHS, but only the negative term on the LHS, increase by~{$\eta$} with customer $\hat{t}$. Thus, thereafter both sides are again equal. In this case as well, both sides remain coupled since $\delta_{X(\hat{t})}(d-1,b+1)=\delta_{X(\hat{t})}(d,b)=d+b$, and $\delta_{X(\hat{t})}(d-1,b)=\delta_{X(\hat{t})}(d,b-1)=d+b-1$.  
This concludes the proof of the Lemma. \hfill \qed

\subsection{Proof of Proposition \ref{prop:p1p2equiv}}\label{appendix:tradeoff}

Our proof crucially relies on the fact that $c_i(\cdot,\cdot)$ is non-increasing in each argument for $i\in[n]$ (this is immediate from the recursive definition in Appendix \ref{sec:appendix_lemma_multimod}). On a high-level we proceed as follows: we first argue that we may assume without loss of generality that $D+B\leq\min\{\sum_i u_i-1,2z+{|\vec{\bar{d}}+\vec{\bar{b}}|_1}\}$. Then we argue that there exists an optimal solution to \ref{opt:original} that uses all of the docks available. Next, we argue that such an optimal solution translates, with the same objective, into a feasible solution to our construction of~\ref{opt:adapted}. This implies that an optimal solution to~\ref{opt:adapted} has objective no worse than the optimal solution to~\ref{opt:original}. Finally, we argue that an optimal solution to~\ref{opt:adapted}, without loss of generality, translates (with same objective) into a feasible solution to \ref{opt:original}.

First, we claim that we may assume without loss of generality that $D+B\leq 2z+|\vec{\bar{d}}+\vec{\bar{b}}|_1$ and $D+B<\sum_i u_i$. The first holds because we know that for any feasible solution we have ${|\vec{d}+\vec{b}|_1\leq 2z+|\vec{\bar{d}}+\vec{\bar{b}}|_1}$, so additional available docks beyond that cannot be allocated no matter what. If the second does not hold (after adjusting $D$ for the first), then $c_i(\cdot,\cdot)$ being non-increasing in both arguments means that an optimal allocation is to allocate $u_i$ docks to each station $i$ and use Algorithm \ref{alg:bike_optimal} to find the bike-optimal allocation for that allocation of docks.

{With those assumptions, consider an optimal solution $(\vec{d}^\star,\vec{b}^\star)$ to \ref{opt:original} that has finite objective and maximizes the number of docks used, i.e., it maximizes {$|\vec{d}^\star+\vec{b}^\star|_1$}.}

\textbf{Claim 1.} We have ${|\vec{d}^\star+\vec{b}^\star|_1}=D+B$.\\
\emph{Proof of claim.}
Suppose we have ${|\vec{d}^\star+\vec{b}^\star|_1}< D+B$. Since $c_i(\cdot,\cdot)$ is non-increasing in both arguments, by the definition of $(\vec{d}^\star,\vec{b}^\star)$ it must be infeasible to add an empty dock to any ~$i$ with~$d_i^\star+b_i^\star<u_i$ (at least one such $i$ exists since ${|\vec{d}^\star+\vec{b}^\star|_1}<D+B<\sum_j u_j$). Adding such an empty dock would not violate ${|\vec{d}^\star+\vec{b}^\star|_1}\leq D+B$ nor would it violate either of ${|\vec{b}^\star|_1} \leq B$ or $d_i^\star+b_i^\star\leq u_i$, so it must violate the operational constraint {$|\vec{d}^\star+\vec{b}^\star-\vec{\bar{d}}-\vec{\bar{b}}|_1\leq 2z$}. Thus, we must have ${ |\vec{d}^\star+\vec{b}^\star - \vec{\bar{d}}-\vec{\bar{b}}|_1=2z}$. If there is a station $k$ with $d_k^\star+b_k^\star < \bar{d}_k+\bar{b}_k\leq u_k$, then an empty dock can feasibly be added to $k$, so such~$k$ cannot exist and we have $d_i^\star+b_i^\star \geq \bar{d}_i+\bar{b}_i$ for every $i$. But then we get~${|\vec{d}^
\star+\vec{b}^\star|_1 = 2z+| \vec{\bar{d}}+\vec{\bar{b}}|_1}$, and thus a contradiction in
$$
D+B\leq 2z+| \vec{\bar{d}}+\vec{\bar{b}}|_1 = |\vec{d}^
\star+\vec{b}^\star|_1< D+B.\hfill\Halmos
$$

Recall from Section \ref{sec: model} that we defined $\bar{D}=D+B-| \vec{\bar{d}}+\vec{\bar{b}}|_1$, i.e., the number of docks that are not in the current allocation but can be added, $\bar{z}=z+\lfloor \frac{\bar{D}}{2}\rfloor$, and a new depot station $\mathcal{D}$ with $l_{\mathcal{D}}=B, u_{\mathcal{D}}=2B+D$, $d_{\mathcal{D}}+b_{\mathcal{D}}=B+\bar{D}$, and $c_{\mathcal{D}}(d,b)=d+b-B$ when $l_{\mathcal{D}}\leq d+b\leq u_{\mathcal{D}}$. Observe that with this additional station, a bike budget $B$ and a dock budget $D+2B$, we have $| \vec{\bar{d}}+\vec{\bar{b}}|_1$ equal the dock budget as required in the definition of \ref{opt:adapted}.  Further, the following is true for $(\vec{d}^\star,\vec{b}^\star)$.

{\textbf{Claim 2.} An allocation of $d_i^\star+b_i^\star$ docks at each station $i\in [n]$ and $B$ docks at $\mathcal{D}$ fulfills the constraints
$|\vec{d}+\vec{b}|_1=2B+D$ and $|\vec{d}+\vec{b}-\vec{\bar{d}}-\vec{\bar{b}}|_1\leq 2\bar{z}$, as well as for each $i\in[n]\cup\{\mathcal{D}\}: l_i\leq d_i+b_i\leq u_i$. 
\\
\emph{Proof of claim.}
By Claim 1 we have a total of $B+\sum_{i\in[n]}d_i^\star+b_i^\star=B+(D+B)=2B+D$ docks allocated. Since $(\vec{d}^\star,\vec{b}^\star)$ has finite objective for \ref{opt:original} and $l_{\mathcal{D}}=B$ this fulfills the required physical constraints. We need to argue that the operational constraint is also fulfilled. We have $$\sum_{i\in[n]}|d_i^\star+b_i^\star-\bar{d}_i-\bar{b}_i|\leq 2z$$ from the feasibility of $(\vec{d}^\star,\vec{b}^\star)$. We also have $\bar{d}_{\mathcal{D}}+\bar{b}_{\mathcal{D}} - d_{\mathcal{D}}^\star-b_{\mathcal{D}}^\star=\bar{D}$. If $\bar{D}$ is even, then this implies
$$
\sum_{i\in[n]}|d_i^\star+b_i^\star-\bar{d}_i-\bar{b}_i|+|\bar{d}_{\mathcal{D}}+\bar{b}_{\mathcal{D}} - d_{\mathcal{D}}^\star-b_{\mathcal{D}}^\star|\leq 2z+2\left\lfloor\frac{\bar{D}}{2}\right\rfloor=2\bar{z}
$$
and the operational constraint is fulfilled. If $\bar{D}$ is odd, notice that it must be possible to get from $(\vec{\bar{d}},\vec{\bar{b}})$ to $(\vec{d}^\star,\vec{b}^\star)$ by adding ~$\bar{D}$ docks and reallocating some number $z'\in\mathbb{N}_0$ of existing docks. Thus, we have $\sum_{i\in[n]}|d_i^\star+b_i^\star-\bar{d}_i-\bar{b}_i|=\bar{D}+2z'\leq 2z$. With~$\bar{D}$ being odd the inequality must hold strictly, so the operational constraint is fulfilled with
$$
\bar{D}+\sum_{i\in[n]}|d_i^\star+b_i^\star-\bar{d}_i-\bar{b}_i|
\leq \bar{D}+2z-1
= 2\left(z +  \frac{\bar{D}}{2}-\frac{1}{2}\right)=2\left(z + \left\lfloor \frac{\bar{D}}{2}\right\rfloor\right)=2\bar{z}.\hfill\Halmos
$$
}

{Observe that the objective of $(\vec{d}^\star,\vec{b}^\star)$ in \ref{opt:original} is the same in \ref{opt:adapted} when allocating: for each $i\in[n]:$ $d_i^\star$ empty docks and $b_i^\star$ full docks; for $\mathcal{D}$, $B$ docks, of which $B-\sum_{i\in[n]}b_i^\star$ are full. Thus,  Claim 2 implies in particular that an optimal solution $(\vec{d}',\vec{b}')$ to the instance of Problem \ref{opt:adapted} we described above, that fulfills $d'_{\mathcal{D}}+b'_{\mathcal{D}}=l_{\mathcal{D}}=B$, has objective no worse than $(\vec{d}^\star,\vec{b}^\star)$ when restricted to stations in $[n]$. We next show that such a solution to \ref{opt:adapted}, restricted to stations in~$[n]$, is feasible for  \ref{opt:original}.}

{
\textbf{Claim 3}. For a feasible solution $(\vec{d}',\vec{b}')$ to the Problem \ref{opt:adapted} as described before, with finite objective, and $d'_{\mathcal{D}}+b'_{\mathcal{D}}=l_{\mathcal{D}}$, an allocation of $d_i'+b_i'$ docks to each station $i\in[n]$ is feasible for  \ref{opt:original}. \\
\emph{Proof of claim.}
The physical and the budget constraints in \ref{opt:original} hold true by feasibility and finite objective of $(\vec{d}',\vec{b}')$ for \ref{opt:adapted}. For the operational constraints, start with the feasibility for \ref{opt:adapted} implying 
$$2z+2\left\lfloor\frac{\bar{D}}{2}\right\rfloor=2\bar{z}\geq \sum_{i} |d_i'+b_i'-(\bar{d}_i+\bar{b}_i)|=\bar{D}+\sum_{i\in[n]} |d_i'+b_i'-(\bar{d}_i+\bar{b}_i)|.$$
Subtracting $\bar{D}$ on both sides we find that $2z\geq \sum_{i\in[n]} |d_i'+b_i'-(\bar{d}_i+\bar{b}_i)|$ which concludes the proof of the claim. \hfill\Halmos
}

What is left to show is that optimal solutions $(\vec{d},\vec{b})$ to \ref{opt:adapted} necessarily fulfill $d_{\mathcal{D}}+b_{\mathcal{D}}=l_{\mathcal{D}}$. This follows from a simple exchange argument and the fact that $c_i(\cdot,\cdot)$ is non-increasing for all $i\in[n]$ whereas it is increasing for $i=\mathcal{D}$. Therefore, for any solution $(\vec{d},\vec{b})$ not fulfilling this property, it {would be} better (and feasible) to remove an additional dock from $\mathcal{D}$ and add an additional dock to any station $i$ with $d_i+b_i<\bar{d}_i+\bar{b}_i$. 
{Suppose no such $i$ exists; then we must have
$$
|(\vec{d}+\vec{b}) - (\vec{\bar{d}}+\vec{\bar{b}})|_1 = \sum_{i\in[n]} \left[(d_i+b_i)-(\bar{d}_i+\bar{b}_i)\right] +  (\bar{d}_{\mathcal{D}}+\bar{b}_{\mathcal{D}})-(d_{\mathcal{D}}+b_{\mathcal{D}})$$
$$= 2 \left[(\bar{d}_{\mathcal{D}}+\bar{b}_{\mathcal{D}})-(d_{\mathcal{D}}+b_{\mathcal{D}})\right] <2\bar{D}\leq 2\bar{z},
$$
so moving a dock from $\mathcal{D}$ to some other station does not violate the operational constraint, and since $\sum_{i\in [n]}u_i > D+B\geq \sum_{i\in[n]} d_i+b_i$, there must exist $i\in[n]$ with $d_i+b_i<u_i$. Thus, in this case as well, there exist a station $i$ such that moving a dock from $\mathcal{D}$ to $i$ would improve the objective and be feasible. It follows that an optimal solution to our instance of P2 fulfills $d_{\mathcal{D}}+b_{\mathcal{D}}=l_{\mathcal{D}}$ and consequently has the same objective when restricted to stations in $[n]$.} This concludes the proof of Proposition \ref{prop:p1p2equiv}.\hfill\Halmos

\subsection{Trading off reallocated and new docks}\label{appendix:tradingoff_new_allocated}
We now discuss a variation of the problem, wherein rather than having fixed budgets that capture the {total} number of docks and the number of docks that may be reallocated, we consider a setting in which there is a joint budget on both. To do so we introduce two new parameters. The parameter~$k$ captures how much more expensive it is to acquire new docks rather than reallocate existing ones. The parameter $C$ bounds the joint cost of reallocating existing and acquiring new docks. In the systems we have worked with, $k$ is so large that the optimal solution would rarely ever acquire new docks. However, it is conceivable that in other systems the cost of reallocating would be larger and thus $k$ would be smaller; for such settings we provide here a variation of the optimization problem. Below~$\bar{D}$ is a new decision variable that captures the number of newly acquired docks and~$z$, previously a parameter, becomes a decision variable. 

\OneAndAHalfSpacedXI
\begin{eqnarray*}
\mathtt{minimize}_{(\vec{d},\vec{b}),z,\bar{D}} & c(\vec{d},\vec{b})\\
s.t. &  \sum_i d_i+b_i &\leq D+B+\bar{D},\\
& \sum_i b_i & \leq B, \\
&{|\vec{\bar{d}}+\vec{\bar{b}}-\vec{d}-\vec{b}|_1}
& \leq 2z+\bar{D},\\
\forall i\in[n]: & l_i\leq \; d_i+b_i &\leq u_i\\
& z+k\bar{D} & \leq C.
\end{eqnarray*}
\DoubleSpacedXI

For each fixed pair of values of $z$ and $\bar{D}$, the discrete-gradient descent algorithm finds an optimal solution by the analysis in Section \ref{sec: algorithm}. Furthermore, it is easily observed that for each value of $\bar{D}$, it is optimal to set $z=C-k\bar{D}$. Hence, one can find an optimal solution by trying all $\lfloor\frac{C}{k}\rfloor$ feasible values of $\bar{D}$ and corresponding values of $z$.



\section{Proofs from Section \ref{sec: algorithm}}\label{appendix_3.4}
\subsection{Proof of Lemma \ref{lemma: bike_opt}}
\label{appendix:bike_opt}

In the proof of Lemma \ref{lemma: bike_opt} we apply the following two auxiliary results.

\noindent\textbf{Fact 1.} For any demand profile $p_i$ and any fixed $d+b$ with $l_i\leq d+b\leq u_i$ the function $c_i(d+b-k,k)$ is convex in $k\in\{0,1,\ldots,d+b\}$.

\noindent\textbf{Fact 2.} Consider $f(\vec{x})=\sum_i f_i(x_i)$ over non-negative integers $x_i$ such that each $f_i$ is convex. If~$\vec{x}$ does not minimize $f(\vec{x})$ subject to $\sum_i x_i=B$, then there exist $i,j$ with $f(\vec{x}+e_i-e_j)<f(\vec{x})$.

We remark that Fact 1 was proven by \cite{raviv2013optimal} and follows by adding inequalities~(2) and~(3) from Definition \ref{def: multimod}, whereas Fact 2 follows from the marginal value analysis of  \cite{fox1966discrete}: compare a suboptimal solution $\vec{x}$ to an optimal solution $\vec{x}^\star$ that minimizes $|\vec{x}^\star-\vec{x}|_1$. Pick~$i,j$ such that $x_i<x_i^\star$, $x_j>x_j^\star$. Then convexity (increasing first derivative) of $f_i$ and $f_j$ implies that $f(\vec{x}+e_i-e_j)-f(\vec{x})=f_i(x_i+1)-f_i(x_i)+f_j(x_j-1)-f_j(x_j)$ $$\leq f_i(x_i^\star)-f_i(x_i^\star-1)+f_j(x_j^\star)-f_j(x_j^\star+1)
=f(\vec{x}^\star)-f(\vec{x}^\star-e_i+e_j)<0,$$
where the last inequality holds since $\vec{x}^\star$ minimizes $|\vec{x}^\star-\vec{x}|$ among optimal solutions.

\noindent\emph{Proof of Lemma.} The bike-optimality condition compares an allocation $(\vec{d}^\star,\vec{b}^\star)$ with other allocations~$(\vec{d}',\vec{b}')$ that fulfill $d_\ell'+b_\ell'=d_\ell^\star+b_\ell^\star$ for all $\ell$ as well as $|\vec{b}'|_1=|\vec{b}^\star|_1$, i.e., for $(\vec{d}^\star,\vec{b}^\star)$ to not be bike-optimal, there must be a better allocation of \emph{bikes} to the same allocation of docks to each station. By Fact~1, we know that (for that same allocation of docks) the cost at each station is a convex function of the number of bikes; combined with Fact~2 we know that if an allocation $(\vec{d}^\star,\vec{b}^\star)$ is not bike-optimal, then there exist two stations such that moving a bike from one to the other improves the objective. 

Now, consider the best allocation out of  $\{o_{ij}(\vec{d},\vec{b}),e_{ij}(\vec{d},\vec{b}),E_{ijh}(\vec{d},\vec{b}),O_{ijh}(\vec{d},\vec{b})\}_{h\in[n]}$, i.e., the one that minimizes~$c(\cdot,\cdot)$. For simplicity, we prove the result for the case that $o_{ij}(\vec{d},\vec{b})$ is that allocation, the other cases follow the same reasoning. If $o_{ij}(\vec{d},\vec{b})$ is bike-optimal there is nothing to show. Else, the reasoning in the last paragraph implies that there exists a bike-move in $o_{ij}(\vec{d},\vec{b})=(\vec{d}-e_i+e_j,\vec{b})$ that improves the objective. That move must be either from or to one of $i$ or $j$; otherwise, it yields the same change in objective to $o_{ij}(\vec{d},\vec{b})$ as it does to $(\vec{d},\vec{b})$. Since the latter is assumed to be a bike-optimal allocation a move of a bike cannot improve the objective.  Let $h$ denote an arbitrary third station. Then we may focus on a bike-move from~$h$ to~$j$, from~$j$ to~$h$, from~$i$ to~$j$, from~$j$ to~$i$, from~$i$ to~$h$, or from~$h$ to~$i$. 
Consider a bike moved from $h$ to $j$. The resulting allocation is $$(\vec{d}-e_i+e_j+e_h-e_j,\vec{b}+e_j-e_h)=(\vec{d}-e_i+e_h,\vec{b}+e_j-e_h)=E_{ijh}(\vec{d},\vec{b}).$$
Since~$o_{ij}(\vec{d},\vec{b})$  minimizes $c(\cdot,\cdot)$ within~$\{e_{ij}(\vec{d},\vec{b}),o_{ij}(\vec{d},\vec{b}),E_{ijh}(\vec{d},\vec{b}),O_{ijh}(\vec{d},\vec{b})\}_{h\in[n]}$, we infer that the move of a bike from $h$ to $j$ does not yield improvement in $o_{ij}(\vec{d},\vec{b})$. Similarly, moving a bike from $i$ to $j$ or from $i$ to $h$ yields the allocations $e_{ij}(\vec{d},\vec{b})$ and $O_{ijh}(\vec{d},\vec{b})$, and thus does not yield a lower objective.



It remains to show that moving a bike from $h$ to $i$, $j$ to $h$, or $j$ to $i$ yields no improvement. These all follow from bike-optimality of $(\vec{d},\vec{b})$ and multimodular Inequality (3). Specifically, since 

\vspace{-.25in}
\begin{eqnarray*}
c_i(d_i-1,b_i)-c_i(d_i-2,b_i+1) \leq c_i(d_i,b_i)-c_i(d_i-1,b_i+1)\\
c_j(d_j+2,b_j-1)-c_j(d_j+1,b_j) \geq c_j(d_j+1,b_j-1)-c_j(d_j,b_j)
\end{eqnarray*}
\noindent an added bike at $i$ improves less and a removed bike at $j$ costs more in $o_{ij}(\vec{d},\vec{b})$ than in $(\vec{d},\vec{b})$.\hfill\halmos

\subsection{Proofs of Lemmas \ref{lemma: neighborhoods} and \ref{lemma: exchange}}\label{sec:appendix_lemma_3}

For the proofs of lemmas \ref{lemma: neighborhoods}--\ref{lemma:const_neighborhood} it is particularly useful to adopt the \emph{diminishing return} interpretation of the multimodular inequalities in Figure \ref{fig:inequalities}. 
For example, Inequality (6) implies that when starting from $(1,5)$ the improvement in the UDF at a station is greater when an empty dock is added than when starting from $(1,6)$. Similarly, Inequality (3) means that among two stations with equal capacity, one with more bikes benefits more from an empty dock. The proofs of the two lemmas are building on these properties by comparing the improvements of the same change for different allocations.
We first prove Lemma \ref{lemma: exchange}, and then derive Lemma \ref{lemma: neighborhoods} as a corollary.

\subsubsection{Proof of Lemma \ref{lemma: exchange}.}
We go through the following cases that exhaust the scenarios in which we have $d_j+b_j<d_j^\star+b_j^\star$ and $d_k+b_k>d_k^\star+b_k^\star$, i.e., when $j$ has \emph{too few} docks and $k$ has \emph{too many} in $(\vec{d},\vec{b})$ relative to $(\vec{d}^\star,\vec{b}^\star)$:
\begin{enumerate}
\item $j$ has too few empty docks in $(\vec{d},\vec{b})$ and $k$ has too many ($d_j<d_j^\star$ and $d_k>d_k^\star$);
\item $j$ has too few full docks in $(\vec{d},\vec{b})$ and $k$ has too many ($b_j<b_j^\star$ and $b_k>b_k^\star$);
\item $j$ has too few empty, enough full docks, $k$ has too many full docks ($d_j<d_j^\star$, $b_j\geq b_j^\star$, and $b_k>b_k^\star$)
\begin{enumerate}
\item and there exists $\ell$ with $d_\ell+b_\ell< d_\ell^\star+b_\ell^\star$, $b_\ell<b_\ell^\star$;
\item and there exists $\ell$ with $d_\ell+b_\ell\geq d_\ell^\star+b_\ell^\star$, $b_\ell<b_\ell^\star$;
\item for all $\ell\not\in\{j,k\}$, we have $b_{\ell}\geq b_{\ell}^\star$, so $|\vec{b}|_1>|\vec{b}^\star|_1$;
\end{enumerate}
\item $j$ has enough empty, too few full docks, $k$ has too few full docks $d_j\geq d_j^\star$, $b_j<b_j^\star$, and $b_k\leq b_k^\star$
\begin{enumerate}
\item and there exists $\ell$ with $d_\ell+b_\ell>d_\ell^\star+b_\ell^\star$ and $b_\ell>b_\ell^\star$;
\item and there exists $\ell$ with $d_\ell+b_\ell\leq d_\ell^\star+b_\ell^\star$ and $b_\ell>b_\ell^\star$;
\item for all $\ell\not\in\{j,k\}$, we have $b_\ell\leq b_\ell^\star$, so $|\vec{b}|_1<|\vec{b}^\star|_1$.
\end{enumerate}
\end{enumerate}

\noindent 
We construct for each case vectors $(\vec{d}',\vec{b}')$ and $(\vec{d}^{\star\star},\vec{b}^{\star\star})$ that fulfill the requirements of the lemma. Our dock-moves from $(\vec{d},\vec{b})$ to construct $(\vec{d}',\vec{b}')$ always move a dock from a station, either $k$ or $\ell$, that has too many docks in $(\vec{d},\vec{b})$ to a station, either $j$ or $\ell$, that has too few docks in $(\vec{d},\vec{b})$, where \emph{too few} is again defined as relative to $(\vec{d}^{\star},\vec{b}^{\star})$; the allocation $(\vec{d}^{\star\star},\vec{b}^{\star\star})$ is always defined by inverting the dock-move from $(\vec{d},\vec{b})$ to $(\vec{d}',\vec{b}')$ and applying that inversion to $(\vec{d}^{\star},\vec{b}^{\star})$. Such a construction immediately implies the second and third property; the first property will be harder to show.

We begin with Case 1.), proving that $(\vec{d}',\vec{b}')=o_{kj}(\vec{d},\vec{b})$ and $(\vec{d}^{\star\star},\vec{b}^{\star\star})=o_{jk}(\vec{d}^\star,\vec{b}^\star)$ fulfill the requirements of the lemma. The second and third property hold by the reasoning above. Before proving the necessary inequality, we summarize the other cases. Case 2.) is symmetric to Case 1.) with~$e_{kj}$ and~$e_{jk}$ replacing~$o_{kj}$ and~$o_{jk}$, and due to that symmetry we omit the proof. Further, the results in~3.a) and~4.a) both reduce to Case~2.). Specifically,~$\ell$ in~3.a) fulfills the conditions of~$j$ in~2.), and~$k$ in~3.a) fulfills~$b_k>b_k^\star$. Thus, assuming Case 2.), $(\vec{d}',\vec{b}')=e_{k\ell}(\vec{d},\vec{b})$ and  $(\vec{d}^{\star\star},\vec{b}^{\star\star})=e_{\ell k}(\vec{d}^\star,\vec{b}^\star)$ fulfill the requirements in~3.a). In 4.a) $\ell$ fulfills the conditions of $k$ in~2.), whereas $j$ in 4.a) fulfills $b_j<b_j^\star$. Assuming again Case 2.), in 4.a) we know that $(\vec{d}',\vec{b}')=e_{\ell j}(\vec{d},\vec{b})$ and  $(\vec{d}^{\star\star},\vec{b}^{\star\star})=e_{j\ell}(\vec{d}^\star,\vec{b}^\star)$ fulfill the requirements. 
Finally, the proofs for cases 3.b) and 4.b) are also symmetric with $(\vec{d}',\vec{b}')=O_{kj\ell}(\vec{d},\vec{b})$ and $(\vec{d}^{\star\star},\vec{b}^{\star\star})=E_{jk\ell}(\vec{d}^{\star},\vec{b}^{\star})$  in the former and $(\vec{d}',\vec{b}')=E_{kj\ell}(\vec{d},\vec{b})$ and $(\vec{d}^{\star\star},\vec{b}^{\star\star})=O_{jk\ell}(\vec{d}^{\star},\vec{b}^{\star})$ fulfilling the requirements in the latter; thus, we only present the proof for 3.b). Cases 3.c) and 4.c) violate the constraint $|\vec{b}|_1=B=|\vec{b}^\star|_1$ and can thus be excluded. Thus, all that is left to show are the inequalities for Case 1.) and Case 3.b). 

The required inequality in Case 1.) is 
\begin{eqnarray*}
c(\vec{d},\vec{b})-c(\vec{d'},\vec{b'}) =\Big(c_j(d_j,b_j)-c_j(d_j+1,b_j)\Big) + \Big(c_k(d_k,b_k) - c_k (d_k-1,b_k)\Big) \\
\geq \Big(c_j(d_j^{\star}-1,b_j^{\star})-c_j(d_j^{\star},b_j^{\star})\Big) + \Big(c_k(d_k^{\star}+1,b_k^{\star})-c_k(d_k^{\star},b_k^{\star})\Big)
= c(\vec{d}^{\star\star},\vec{b}^{\star\star})-c(\vec{d}^{\star},\vec{b}^{\star}).
\end{eqnarray*}

We prove the inequality by comparing the respective $j$- and $k$-terms, starting with

\begin{figure}[ht]
    \centering
    \includegraphics[width=.5\textwidth]{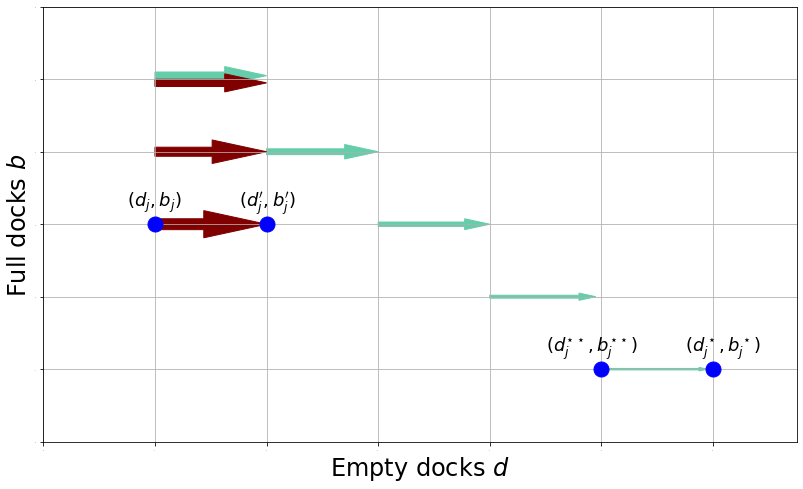}
    \caption{Visual illustration of the inequalities argued for station~$j$ in the proof of Lemma \ref{lemma: exchange}, Case~1.):\\ since $d_j<d_j^\star$ and $d_j+b_j<d_j^\star+b_j^\star$ Inequality (3), in acquamarine, and (6), in brown, bound the difference $c_j(d_j,b_j)-c_j(d_j+1,b_j)$ from below by~$c_j(d_j^{\star}-1,b_j^{\star})-c_j(d_j^{\star},b_j^{\star})$.}
    \label{fig:lemma_3}
\end{figure}

\vspace{-.25in}
\[
c_j(d_j,b_j)-c_j(d_j+1,b_j)\geq c_j(d_j^{\star}-1,b_j^{\star})-c_j(d_j^{\star},b_j^{\star}).
\]
We refer the reader to Figure \ref{fig:lemma_3} for a graphical presentation of the sequence of inequalities from Definition \ref{def: multimod} we apply (thicker arrows correspond to greater improvement of an added empty dock).


\noindent First, apply Inequality~(3) $t$ times to the RHS to get an upper bound of $$c_j(d_j^\star-1-t,b_j^\star+t)-c_j(d_j^\star-t,b_j^\star+t).$$ Setting $t=d_j^\star-d_j-1\geq 0$, we find that the RHS is bounded above by 
\vspace{-.15in}
\[c_j(d_j,b_j^\star+d_j^\star-d_j-1)-c_j(d_j+1,b_j^\star+d_j^\star-d_j-1).\]
\vspace{-.45in}

\noindent Then, apply Inequality~(6) repeatedly to the LHS to show that $\forall s\geq 0$ the LHS is bounded below by~$c_j(d_j,b_j+s)-c_j(d_j+1,b_j+s)$. Hence, by setting $s=b_j^\star+d_j^\star-d_j-b_j-1$, which is non-negative since $b_j+d_j<b_j^\star+d_j^\star$, we bound the LHS from below by
\[
c_j(d_j,b_j+b_j^\star+d_j^\star-d_j-b_j-1)-c_j(d_j+1,b_j+b_j^\star+d_j^\star-d_j-b_j-1).
\] 

\noindent This equals the upper bound on the RHS and thus proves the desired inequality.

Similarly, to show 
\begin{eqnarray}\label{eqn:above}
c_k (d_k-1,b_k)-c_k(d_k,b_k) \leq c_k(d_k^\star,b_k^\star)-c_k(d_k^\star+1,b_k^\star),
\end{eqnarray}

\noindent we apply Inequality~(3) $d_k-d_k^\star-1$ times to bound the LHS in (\ref{eqn:above}) from above by
$$c_k(d_k^\star,b_k+d_k-d_k^\star-1)-c_k(d_k^\star+1,b_k+d_k-d_k^\star-1).$$
Thereafter, we apply Inequality~(6) $(b_k+d_k)-(d_k^\star+b_k^\star) -1\geq 0$ times to obtain the desired bound.

\begin{figure}[ht]
    \centering
    \includegraphics[width=.35\textwidth]{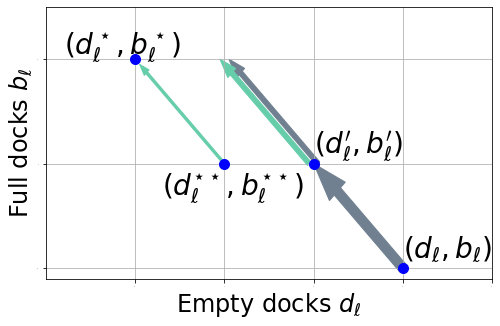}
    \caption{Visual illustration of the inequalities argued for station~$\ell$ in the proof of Lemma \ref{lemma: exchange}, Case~3.b).}
    \label{fig:lemma_3_ell}
\end{figure}

It remains to show that in Case 3.b), with $(\vec{d}',\vec{b}')=O_{kj\ell}(\vec{d},\vec{b})$ and $(\vec{d}^{\star\star},\vec{b}^{\star\star})=E_{jk\ell}(\vec{d}^{\star},\vec{b}^{\star})$, we have the required inequality~$c(\vec{d},\vec{b})-c(\vec{d'},\vec{b'})\geq c(\vec{d}^{\star\star},\vec{b}^{\star\star})-c(\vec{d}^{\star},\vec{b}^{\star})$. Notice first that all terms not involving $j,k,$ and~$\ell$  cancel out on each side. Further, the terms involving $j$ can be bounded the same way as in Case 1.); the bounds for $k$, receiving a full dock now instead of an empty dock in Case 1.), are symmetric to the ones in Case 1.) with (1) and (2) replacing (3) and (6). Thus, we only need to derive
\begin{eqnarray*}
c_\ell (d_\ell,b_\ell)-c_\ell(d_\ell-1,b_\ell+1) \geq c_\ell(d_\ell^\star+1,b_\ell^\star-1)-c_\ell(d_\ell^\star,b_\ell^\star)
\end{eqnarray*}
\noindent
We obtain this (see Figure \ref{fig:lemma_3_ell}) by bounding with $s=b_\ell^\star-b_\ell-1\geq 0$
$$c_\ell (d_\ell,b_\ell)-c_\ell(d_\ell-1,b_\ell+1) \geq 
c_\ell (d_\ell-s,b_\ell+s)-c_\ell(d_\ell-1-s,b_\ell+1+s),$$
which follows from Fact 1, and then bounding the resulting term as
$$c_\ell (d_\ell-b_\ell^\star+b_\ell+1,b_\ell^\star-1)
-c_\ell(d_\ell-b_\ell^\star+b_\ell,b_\ell^\star)\geq c_\ell(d_\ell^\star+1,b_\ell^\star-1)-c_\ell(d_\ell^\star,b_\ell^\star),$$
which follows from rewriting Inequality (3) as $f(d+1,b-1)-f(d,b)\geq f(d,b-1)-f(d-1,b) $ and applying it $(b_\ell+d_\ell)-(b_\ell^\star+d_\ell^\star)\geq0$ times. \hfill \halmos

\subsubsection{Proof of Lemma \ref{lemma: neighborhoods}.}
We first argue that we may assume without loss of generality that there is no solution $(\vec{d}^{\star\star},\vec{b}^{\star\star})$ in the neighborhood of $(\vec{d}^\star,\vec{b}^\star)$ that has better objective than $(\vec{d},\vec{b})$ does and is closer to $(\vec{d},\vec{b})$ in dock-move distance than $(\vec{d}^\star,\vec{b}^\star)$ is. Indeed, if such $(\vec{d}^{\star\star},\vec{b}^{\star\star})$ exists, then it must be component-wise in between $(\vec{d},\vec{b})$ and $(\vec{d}^\star,\vec{b}^\star)$, i.e., we have either $d_i+b_i \leq d_i^{\star\star}+b_i^{\star\star}\leq d_i^{\star}+b_i^{\star}$ or~$d_i+b_i \geq d_i^{\star\star}+b_i^{\star\star}\geq d_i^{\star}+b_i^{\star}$ for every~$i$. But then, a dock-move from $(\vec{d},\vec{b})$ towards $(\vec{d}^{\star\star},\vec{b}^{\star\star})$ also reduces the dock-move distance to $(\vec{d}^{\star},\vec{b}^{\star})$. Thus, we may relabel such $(\vec{d}^{\star\star},\vec{b}^{\star\star})$ as $(\vec{d}^{\star},\vec{b}^{\star})$, and argue the result for $(\vec{d},\vec{b})$ and this new $(\vec{d}^{\star},\vec{b}^{\star})$ instead (implying the result for the original $(\vec{d}^{\star},\vec{b}^{\star})$).

Now, without loss of generality, suppose there is no solution $(\vec{d}^{\star\star},\vec{b}^{\star\star})$ in the neighborhood of $(\vec{d}^\star,\vec{b}^\star)$ that has better objective than $(\vec{d},\vec{b})$.
Then we can apply Lemma~\ref{lemma: exchange} to~$(\vec{d},\vec{b})$ and $(\vec{d}^\star,\vec{b}^\star)$ to find~$(\vec{d}^{\star\star},\vec{b}^{\star\star})$ and~$(\vec{d}',\vec{b}')$ as in the lemma, with the additional property that $$c(\vec{d}^\star,\vec{b}^\star)< c(\vec{d},\vec{b})\leq c(\vec{d}^{\star\star},\vec{b}^{\star\star}), \text{ implying in particular }
c(\vec{d},\vec{b})-c(\vec{d'},\vec{b'}) \geq c(\vec{d}^{\star\star},\vec{b}^{\star\star})-c(\vec{d}^{\star},\vec{b}^{\star})>0. \halmos
$$

\subsection{Proof of Lemma \ref{lemma:const_neighborhood}}\label{appendix_proof_kopt}


We know by Lemma \ref{lemma: bike_opt} that $(\vec{d}^{{{r}}+1},\vec{b}^{{{r}}+1})$ must be bike-optimal. Let~$i$ and~$j$ respectively be the stations from which a dock is taken and to which it is moved in iteration~$r+1$. If that move involves a third station we denote it by $h$ (recall that a dock-move from $i$ to $j$ can take an additional bike from $i$ to a third station $h$ or take one from $h$ to $j$.) We first argue that none of the following yield a better solution {at dock-move distance at most $r+1$ to $(\vec{d},\vec{b})$} 
when applied to $(\vec{d}^{r+1},\vec{d}^{r+1})$:
\begin{enumerate}
\item Any dock-move not involving any of $i$, $j$, and $h$;
\item Any dock-move from station $i$ to some station $\ell$;
\item Any dock-move from some station $\ell$ to station $j$;
\item Any dock-move in which $i$ receives a dock from some station $\ell$;
\item Any dock-move in which some station $\ell$ receives a dock from $j$. 
\end{enumerate}
{In order to reason about these dock-moves, it helps to define the  notation
$$
S_r(\vec{d},\vec{b})=\left\{(\vec{d}',\vec{b}'):|\vec{d}'+\vec{b}'-\vec{d}-\vec{b}|_1\leq 2r\right\}
$$
to describe the set of allocations at dock-move distance at most $r$ from an allocation $(\vec{d},\vec{b})$, where we at times drop the argument $(\vec{d},\vec{b})$ when clear from context.
}

{Now, for the first kind of move we} observe that for each station $\ell$ involved we have $d_\ell^r=d_\ell^{r+1}, b_\ell^r=b_\ell^{r+1}$; thus, it yields the same improvement when applied to $(\vec{d}^r,\vec{b}^r)$ as when applied to $(\vec{d}^{r+1},\vec{b}^{r+1})$, and has the same effect on the dock-move distance to $(\vec{d},\vec{b})$, i.e., if such a dock-move yields a solution within $S_{{{r}}+1}(\vec{d},\vec{b})$ when applied to $(\vec{d}^{r+1},\vec{b}^{r+1})$, then applying it to $(\vec{d}^r,\vec{b}^r)$ also yields a solution within~$S_{{{r}}}(\vec{d},\vec{b})$. Hence, by the assumption that $(\vec{d}^r,\vec{b}^r)$ is optimal within $S_r(\vec{d},\vec{b})$, it cannot yield any improvement.

We can argue similarly about the second and third: applying such moves to  $(\vec{d}^{r+1},\vec{b}^{r+1})$ yields solutions in $S_{{{r}}+1}(\vec{d},\vec{b})$ only if they yield solutions in $S_{{{r}}}(\vec{d},\vec{b})$ when applied to $(\vec{d}^r,\vec{b}^r)$. However (see Figure~\ref{fig:inequalities}), inequalities (1), (4), (5), and (6) in the definition of multimodularity imply that such moves increase the objective at $i$ more (decrease the objective at $j$ less) when applied to $(\vec{d}^{r+1},\vec{b}^{r+1})$ than when applied to $(\vec{d}^r,\vec{b}^r)$. For example, suppose in period $r+1$ an empty dock is taken away from $i$: Then Inequality (4) lower bounds the cost of taking an empty dock from $i$, and Inequality~(1) lower bounds the cost of taking a full dock from $i$, in $(\vec{d}^{r+1},\vec{b}^{r+1})$ relative to those same actions in~$(\vec{d}^r,\vec{b}^r)$. Since $(\vec{d}^r,\vec{b}^r)$ is assumed optimal within $S_r(\vec{d},\vec{b})$, such moves yield no improvement.

The fourth and fifth kind of move can be excluded because the allocation of docks (perhaps with different allocations of bikes) resulting from them is among the choices considered by the gradient-descent algorithm in iteration~$r+1$. Since Lemma \ref{lemma: bike_opt} implies that  among such allocations considered by the algorithm, there is a bike-optimal one, the algorithm's choice to pick $(\vec{d}^{r+1},\vec{b}^{r+1})$ over this alternative implies that $(\vec{d}^{r+1},\vec{b}^{r+1})$ is no worse.

We are left with dock-moves either from or to $h$ as well as dock-moves that involve one of the three stations $i,j$, and $h$ only via a bike being moved. We go through these in a sequence of case distinctions. For the rest of the proof we assume the move in period $r+1$ is ~$E_{ijh}$, i.e., $(\vec{d}^{r+1},\vec{b}^{r+1})=(\vec{d}^r-e_i+e_h,\vec{b}^r+e_j-e_h)$. The arguments for $e_{ij},o_{ij}, $ and $O_{ijh}$ are symmetric. We begin with dock-moves that only involve $i,j,$ or $h$ via a bike being moved.


\subsubsection*{A bike to $i,j$.}  
\begin{figure}[ht]
    \centering
    \includegraphics[width=.35\textwidth]{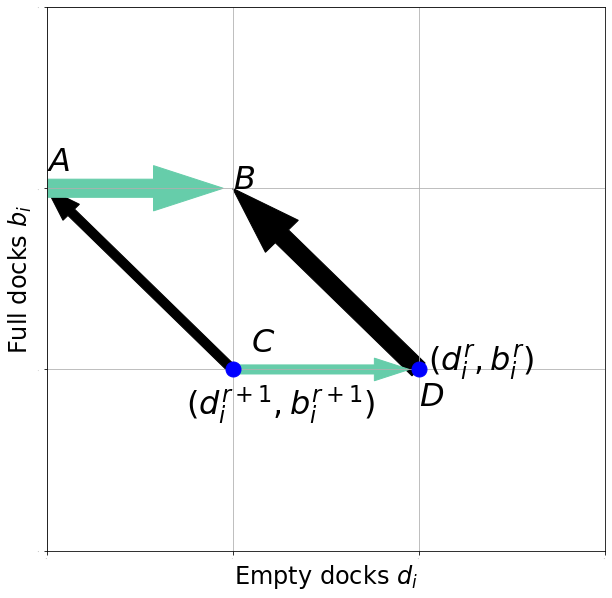}
    \caption{Visual illustration of why a bike being added to $i$ yields no more improvement to $(d_i^{r+1},b_i^{r+1})$ than to $(d_i^{r},b_i^{r})$ after $E_{ijh}$ in period $r+1$: Inequality~(3) gives, with $A,B,C,D$ representing the respective value of $c_i(\cdot,\cdot)$: $A-B\geq C-D \Longleftrightarrow C-A\geq D-B$.}
    \label{fig:bike_moved_to_i}
\end{figure}

Consider a move that involves a bike moved to $i$. Inequality (3) implies (see Figure~\ref{fig:bike_moved_to_i}) that the improvement at~$i$ is no more in $(\vec{d}^{r+1},\vec{b}^{r+1})$  than it would be in $(\vec{d}^{r},\vec{b}^{r})$. Inequality (2) implies the same for a move of a bike to station~$j$.

\subsubsection*{A bike from $i$ or $j$.} We next consider moves that take a bike from $i$ or $j$, i.e., $E_{\ell mi}$ or $E_{\ell mj}$. The algorithm's choice to carry out $E_{ijh}$ rather than $o_{ij}$ or $e_{ij}$ in iteration $r+1$ implies that taking a bike from either $i$ or $j$ in allocation $(\vec{d}^{r+1},\vec{b}^{r+1})$ has cost no less than the cost of taking a bike from $h$ in  allocation~$(\vec{d}^{r},\vec{b}^{r})$. Consider~$\ell$ and $m$ such that~$E_{\ell mh}(\vec{d}^{r+1},\vec{b}^{r+1})\in S_{r+1}$; for such $\ell$ and $m$ we must also have $E_{\ell mh}(\vec{d}^{r},\vec{b}^{r})\in S_{r}$. By the inductive assumption, that move would not yield an improvement to $(\vec{d}^{r},\vec{b}^{r})$; since taking a bike from $i$ or $j$ in allocation $(\vec{d}^{r+1},\vec{b}^{r+1})$ has a cost that is no less, by the argument above, we infer that the moves  $E_{\ell mi}(\vec{d}^{r+1},\vec{b}^{r+1})$ and  $E_{\ell mj}(\vec{d}^{r+1},\vec{b}^{r+1})$ also yield no improvement.

\subsubsection*{A bike from or to $h$.} 
For a move of a bike from $h$, Fact 1 (see Figure \ref{fig:inequalities}) implies the cost at $(\vec{d}^{r+1},\vec{b}^{r+1})$ is greater-equal to the cost at $(\vec{d}^r,\vec{b}^r)$. 

Consider a bike being added to $h$, i.e., $O_{\ell m h}$ for some $\ell,m$ such that~$O_{\ell m h}(\vec{d}^{r+1},\vec{b}^{r+1})\in S_{r+1}$. For $O_{\ell m h}(\vec{d}^{r+1},\vec{b}^{r+1})$ to be in $S_{r+1}$  at least one of $d_m^{r}+b_m^{r}<d_m+b_m$ or $d_\ell^{r}+b_\ell^{r}>d_\ell+b_\ell$ must hold true. 

Consider first $d_m^{{r}}+b_m^{{r}}<d_m+b_m$. Then any dock-move from $i$ to $m$ in  $(\vec{d}^{{r}},\vec{b}^{{r}})$ yields an allocation in $S_{{r}}(\vec{d},\vec{b})$, implying in particular that $c(o_{im}(\vec{d}^{{r}},\vec{b}^{{r}}))\geq c(\vec{d}^{{r}},\vec{b}^{{r}})$. Further, observe that the gradient-descent algorithm chooses $E_{ijh}$ rather than $e_{\ell j}$ in iteration $r+1$; as the impact of both at $j$ is the same (one added full dock), this implies that the increase in objective in taking a dock and a bike from $\ell$ is at least the increase at $i$ (from taking an empty dock) and $h$ (from taking a bike) in iteration~$r+1$. Thus, 
$$c(O_{\ell mh}(\vec{d}^{{r+1}},\vec{b}^{{r+1}}))-c(\vec{d}^{{r+1}},\vec{b}^{{r+1}})\geq c(o_{im}(\vec{d}^{{r}},\vec{b}^{{r}}))- c(\vec{d}^{{r}},\vec{b}^{{r}})\geq 0.$$

Next, consider $d_\ell^{{r}}+b_\ell^{{r}}>d_\ell+b_\ell$. We observe that $e_{\ell j}(\vec{d}^{{r}},\vec{b}^{{r}})\in S_{{r}}(\vec{d},\vec{b})$, 
and that the increase in objective at $\ell$ is the same for~$O_{\ell mh}(\vec{d}^{{r+1}},\vec{b}^{{r+1}})$ and~$e_{\ell j}(\vec{d}^{{r}},\vec{b}^{{r}})$. That increase is bounded below by the decrease in objective at $j$ due to  $e_{\ell j}(\vec{d}^{{r}},\vec{b}^{{r}})$ because $c(\vec{d}^r,\vec{b}^r)\leq c(e_{\ell j}(\vec{d}^r,\vec{b}^r))$ holds by the inductive assumption with $e_{\ell j}(\vec{d}^{{r}},\vec{b}^{{r}})\in S_{{r}}(\vec{d},\vec{b})$. Since the gradient-descent algorithm chooses~$E_{ijh}$ rather than~$o_{im}$ in iteration~$r+1$, that decrease at $j$ is greater-equal to the combined decrease in objective at~$m$ and $h$ due to $O_{\ell mh}(\vec{d}^{{r+1}},\vec{b}^{{r+1}})$. Combining all of the above we find that the increase in objective at $\ell$ due to $O_{\ell mh}(\vec{d}^{{r+1}},\vec{b}^{{r+1}})$ is greater than the decrease at $m$ and $h$.

This leaves us with dock-moves from and to $h$.

\subsubsection*{Full dock from or empty dock to $h$.}
\begin{figure}[ht]
    \centering
    \includegraphics[width=.35\textwidth]{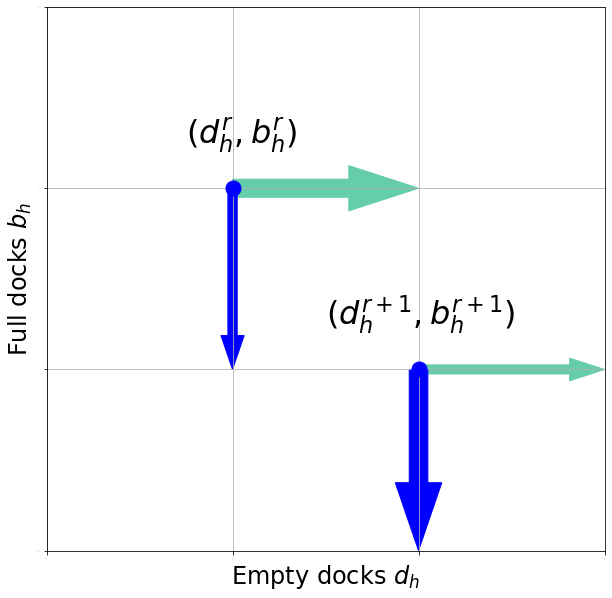}
    \caption{Bounding the improvement of an added empty dock (cost of a removed full dock) at $h$ in $(\vec{d}^{r+1},\vec{b}^{r+1})$.}
    \label{fig:lemma_5}
\end{figure}

 In this case (see blue arrows in Figure \ref{fig:lemma_5}), by Inequality~(2), a move of a full dock from~$h$, i.e., $e_{h\ell}$ or $O_{h\ell m}$ for some~$\ell,m\not\in\{i,j\}$, increases the objective at $h$ by at least as much (and decreases at $\ell$ and $m$ by the same amount) in $(\vec{d}^{r+1},\vec{b}^{r+1})$ as in $(\vec{d}^{r},\vec{b}^{r})$. Similarly, by Inequality~(3) (see aquamarine arrows in Figure \ref{fig:lemma_5}) the move of an empty dock to~$h$, i.e.,~$o_{\ell h}$ or~$O_{\ell hm}$ for any~$\ell$ and $m$, has no more improvement to $(\vec{d}^{r+1},\vec{b}^{r+1})$ than to $(\vec{d}^{r},\vec{b}^{r})$. Since in both of these cases we have $d_s^r+b_s^r=d_s^{r+1}+b_s^{r+1}$ for $s\in\{h,\ell,m\}$ it must be the case that if such dock-moves yield allocations in $S_{r+1}(\vec{d},\vec{b})$ when applied to $(\vec{d}^{r+1},\vec{b}^{r+1})$, then they must also yield allocations in~$S_{r}(\vec{d},\vec{b})$ when applied to~$(\vec{d}^{r},\vec{b}^{r})$. It follows that none of $e_{h\ell}$, $O_{h\ell m}$, $o_{\ell h}$, or~$O_{\ell hm}$ can yield improvement within~$S_{r+1}(\vec{d},\vec{b})$  when applied to $(\vec{d}^{r+1},\vec{b}^{r+1})$.

\subsubsection*{Empty dock from or full dock to $h$.}
Moves of an empty dock from $h$ (or a full dock to $h$) have a lower cost (greater improvement) at $h$ in $(\vec{d}^{r+1},\vec{b}^{r+1})$ than in $(\vec{d}^{r},\vec{b}^{r})$ and require a more careful argument. Suppose $o_{h\ell}$ yielded an improvement within $S_{r+1}(\vec{d},\vec{b})$ when applied to $(\vec{d}^{r+1},\vec{b}^{r+1})$; the cases for $E_{h\ell m}$,~$e_{\ell h}$, and $E_{\ell h m} $ are similar. If~$d^{{r}}_h+b^{{r}}_h>d_h+b_h$ and~$d^{{r}}_\ell +b^{{r}}_\ell < d_\ell +b_\ell $ held true, then moving a dock from~$h$ to~$\ell$ reduces the dock-move distance to~$(\vec{d},\vec{b})$, so~$o_{h\ell }(E_{ijh}(\vec{d}^{{r}},\vec{b}^{r}))\in S_{r}(\vec{d},\vec{b})$, and thus~$c(o_{h\ell}(\vec{d}^{{r+1}},\vec{b}^{{r+1}}))=c(o_{h\ell }(E_{ijh}(\vec{d}^{{r}},\vec{b}^{{r}})))<c(\vec{d}^{{r}},\vec{b}^{{r}})$ which contradicts the inductive assumption. Thus, for $o_{h\ell }(E_{ijh}(\vec{d}^{{r}},\vec{b}^{{r}}))$ to yield an improvement, it must be the case that~$o_{h\ell }(E_{ijh}(\vec{d}^{{r}},\vec{b}^{{r}}))\in S_{{{r}}+1}(\vec{d},\vec{b})\setminus S_{{{r}}}(\vec{d},\vec{b})$; it follows that either
\begin{enumerate}
\item $d^{{r}}_h+b^{{r}}_h>d_h+b_h$ and $d^{{r}}_\ell +b^{{r}}_\ell \geq d_\ell +b_\ell $ or
\item $d^{{r}}_h+b^{{r}}_h\leq d_h+b_h$ and $d^{{r}}_\ell +b^{{r}}_\ell <d_\ell +b_\ell $.
\end{enumerate}
Indeed, if these did not hold, then a dock-move from $h$ to $\ell$ would either yield a solution in $S_{{r}}$ or one not in $S_{{{r}}+1}$. 
We now show that in both cases the following expression is non-negative:

\vspace{-.45in}
\begin{eqnarray*}
& c(o_{h\ell }(\vec{d}^{{{r}}+1},\vec{b}^{{{r}}+1})) - c(\vec{d}^{{{r}}+1},\vec{b}^{{{r}}+1}) =
& c_h(d^{{r}}_h,b^{{r}}_h-1)-c_h(d^{{r}}_h+1,b^{{r}}_h-1)+c_\ell(d^{{r}}_\ell +1,b^{{r}}_\ell )-c_\ell(d^{{r}}_\ell ,b^{{r}}_\ell ).
\end{eqnarray*}

\vspace{-.15in}
When $d^{{r}}_h+b^{{r}}_h>d_h+b_h$ and $d^{{r}}_\ell +b^{{r}}_\ell \geq d_\ell +b_\ell $, observe that $e_{hj}(\vec{d}^{{{r}}},\vec{b}^{{{r}}})\in S_{{r}}(\vec{d},\vec{b})$. Thus, the inductive assumption implies that $c_h(d^{{r}}_h,b^{{r}}_h-1)+c_j(d^{{r}}_j,b^{{r}}_j+1)\geq c_h(d^{{r}}_h,b^{{r}}_h)+c_j(d^{{r}}_j,b^{{r}}_j)$, or equivalently: $$c_h(d^{{r}}_h,b^{{r}}_h){\stackrel{\star}{\leq}} c_h(d^{{r}}_h,b^{{r}}_h-1)+c_j(d^{{r}}_j,b^{{r}}_j+1) - c_j(d^{{r}}_j,b^{{r}}_j).$$ Further, by the choice of the {gradient-descent} algorithm in iteration $r+1$, we know that $c(E_{ijh}(\vec{d}^{{{r}}},\vec{b}^{{{r}}}))\leq c(o_{i\ell}(\vec{d}^{{{r}}},\vec{b}^{{{r}}}))$, i.e., an additional empty dock at~$\ell$ has no more improvement than an additional dock and an additional bike at~$j$ minus the cost of taking a bike from~$h$. Formally,

\vspace{-.3in}
\OneAndAHalfSpacedXI
\begin{eqnarray*}
 c_\ell (d^{{r}}_\ell +1,b^{{r}}_\ell )-c_\ell (d^{{r}}_\ell ,b^{{r}}_\ell )
\leq \left[c_j(d^{{r}}_j,b^{{r}}_j)-c_j(d^{{r}}_j,b^{{r}}_j+1)\right] + \left[c_h(d^{{r}}_h,b^{{r}}_h) - c_h(d^{{r}}_h+1,b^{{r}}_h-1)\right].
\end{eqnarray*}
\DoubleSpacedXI

Plugging in the upper bound (starred inequality) on $c_h(d_h^r,b_h^r)$ the $j$-terms cancel out and we are left with $c_h(d^{{r}}_h,b^{{r}}_h-1)-c_h(d^{{r}}_h+1,b^{{r}}_h-1)$, giving us $c(o_{h\ell }(\vec{d}^{{{r}}+1},\vec{b}^{{{r}}+1})) - c(\vec{d}^{{{r}}+1},\vec{b}^{{{r}}+1})\geq 0$.


When $d^{{r}}_h+b^{{r}}_h\leq d_h+b_h$ and $d^{{r}}_\ell +b^{{r}}_\ell <d_\ell +b_\ell $ we know that $o_{i\ell}(\vec{d}^{{{r}}},\vec{b}^{{{r}}})\in S_{{{r}}}(\vec{d},\vec{b})$; thus the inductive assumption implies that moving an empty dock from $i$ to $\ell$ in $(\vec{d}^{{{r}}},\vec{b}^{{{r}}})$ does not yield an improvement. Specifically, we can bound $c_\ell (d^{{r}}_\ell +1,b^{{r}}_\ell )-c_{\ell}(d_\ell^r,b_\ell^r)
\geq c_i(d_i^r,b_i^r)-c_i(d_i^r-1,b_i^r)$. 
Further, the algorithm's choice to carry out $E_{ijh}$ rather than $e_{hj}$ in iteration~$r+1$, i.e., to take an empty dock from $i$ and a bike from $h$ rather than taking both from $h$, implies that $$c_i(d^{{r}}_i,b^{{r}}_i)-c_i(d^{{r}}_i-1,b^{{r}}_i)\geq c_h(d^{{r}}_h+1,b^{{r}}_h-1)-c_h(d^{{r}}_h,b^{{r}}_h-1).$$ Combining these two inequalities again implies the nonnegativity of $c_h(d^{{r}}_h,b^{{r}}_h-1)-c_h(d^{{r}}_h+1,b^{{r}}_h-1)+c_\ell(d^{{r}}_\ell +1,b^{{r}}_\ell )-c_\ell(d^{{r}}_\ell ,b^{{r}}_\ell )$. 
\hfill\halmos

\subsection{Proof of Theorem \ref{thm: koptimal}}\label{appendix:proof_theorem}

In order to prove Theorem \ref{thm: koptimal} we first need the following lemma{, for which we recall the definition of~$S_r(\vec{d},\vec{b})=\left\{(\vec{d}',\vec{b}'):|\vec{d}'+\vec{b}'-\vec{d}-\vec{b}|_1\leq 2r\right\}$ from the previous section.}

\begin{lemma}\label{lemma: iexists}
With $(\vec{d}^\star,\vec{b}^\star)\in S_{{{r}}+1}(\vec{d},\vec{b})$ and $c(\vec{d}^\star,\vec{b}^\star)<c(\vec{d}^{{{r}}+1},\vec{b}^{{{r}}+1})$, and $(\vec{d}^{{{r}}},\vec{b}^{{{r}}})$, $(\vec{d}^{{{r}}+1},\vec{b}^{{{r}}+1})$ as in Lemma \ref{lemma:const_neighborhood}, there exist $j,k$ with $d^\star_j+b^\star_j>d^{{{r}}}_j+b^{{{r}}}_j \geq d_j+b_j$ and $d^\star_k+b^\star_k<d^{{{r}}}_k+b^{{{r}}}_k \leq d_k+b_k$.
\end{lemma}

{\emph{Proof of Lemma. }} We know from Lemma \ref{lemma: neighborhoods} that there exists a solution {$(\vec{d}',\vec{b}')\in N(\vec{d}^{{{r}}+1},\vec{b}^{{{r}}+1})$} that has smaller objective and smaller dock-move distance to $(\vec{d}^\star,\vec{b}^\star)$ than $(\vec{d}^{{{r}}+1},\vec{b}^{{{r}}+1})$ does.  Further, since~$(\vec{d}^{{{r}}+1},\vec{b}^{{{r}}+1})$ is a local optimum within~$S_{{{r}}+1}(\vec{d},\vec{b})$ by Lemma~\ref{lemma:const_neighborhood} and since $(\vec{d}',\vec{b}')\in N(\vec{d}^{{{r}}+1},\vec{b}^{{{r}}+1})$, it follows that $(\vec{d}',\vec{b}')\not\in S_{{{r}}+1}(\vec{d},\vec{b})$. Consider the station $j$ that receives a dock {in the dock-move from $(\vec{d}^{{{r}}+1},\vec{b}^{{{r}}+1})$ to $(\vec{d}',\vec{b}')$}, i.e., $d^{{{r}}+1}_j+b^{{{r}}+1}_j+1=d_j'+b_j'$, and the station $k$ that has a dock taken away, i.e., $d^{{{r}}+1}_k+b^{{{r}}+1}_k-1=d_k'+b_k'$. We prove that the desired inequalities are fulfilled at $j$ by showing that $d_j^\star+b_j^\star\geq d_j'+b_j'=d_j^{r+1}+b_j^{r+1}+1> d_j^r+b_j^r\geq d_j+b_j$; the proof for $k$ is symmetric.

First, we must have $d_j^\star+b_j^\star\geq d_j'+b_j'$  for $(\vec{d}',\vec{b}')$ to be at a closer dock-move distance to~$(\vec{d}^\star,\vec{b}^\star)$ than~$(\vec{d}^{r+1},\vec{b}^{r+1})$.
Next, for $(\vec{d}',\vec{b}')\not\in S_{{{r}}+1}(\vec{d},\vec{b})$, the dock-move from $(\vec{d}^{r+1},\vec{b}^{r+1})$ to $(\vec{d}',\vec{b}')$ must increase the dock-move distance to $(\vec{d},\vec{b})$; this requires in particular that  $d_j'+b_j'>d_j+b_j$.  But then, since $d^{r+1}_j+b^{r+1}_j+1=d_j'+b_j'$, we find $d^{r+1}_j+b^{r+1}_j\geq d_j+b_j$. 

By definition  $(\vec{d}^{r},\vec{b}^{r})$ and $(\vec{d}^{r+1},\vec{b}^{r+1})$ are one dock-move apart, so (i) $|d^{r}_j+b^{r}_j-d^{r+1}_j-b^{r+1}_j|\leq 1$. Also,~(ii)~$(\vec{d}^{r+1},\vec{b}^{r+1})$ is at larger dock-move distance from $(\vec{d},\vec{b})$ than $(\vec{d}^{r},\vec{b}^{r})$ --- else we would have~$(\vec{d}^{r+1},\vec{b}^{r+1})\in S_{r}(\vec{d},\vec{b})$ with no improvement in the period $r+1$ move from $(\vec{d}^{r},\vec{b}^{r})$ to $(\vec{d}^{r+1},\vec{b}^{r+1})$. 
 
Suppose for the sake of contradiction that $d_j^r+b_j^r<d_j+b_j$. Since $d_j^{r+1}+b_j^{r+1}\geq d_j+b_j$, (i) implies that $d^{r+1}_j+b^{r+1}_j=d_j+b_j$, so the dock-move from $(\vec{d}^{r},\vec{b}^{r})$ to $(\vec{d}^{r+1},\vec{b}^{r+1})$ does not increase the dock-move distance to $(\vec{d},\vec{b})$ which contradicts (ii). Thus, we must have $d_j^r+b_j^r\geq d_j+b_j$.
 
 Since $d^{r+1}_j+b^{r+1}_j\geq d_j+b_j$ (ii) also implies that we have $d^{r+1}_j+b^{r+1}_j\geq d_j^r+b_j^r$; otherwise the dock-move from $(\vec{d}^{r},\vec{b}^{r})$ to $(\vec{d}^{r+1},\vec{b}^{r+1})$ cannot increase the dock-move distance to $(\vec{d},\vec{b})$. Combining the above we find that $d^\star_j+b^\star_j>d^{{{r}}}_j+b^{{{r}}}_j \geq d_j+b_j$. \hfill\Halmos

{\emph{Proof of Theorem. }
Define $(\vec{d}^{{{r}}},\vec{b}^{{{r}}})$, $(\vec{d}^{{{r}}+1},\vec{b}^{{{r}}+1})$ as in Lemma \ref{lemma:const_neighborhood} and suppose $(\vec{d}^{{{r}}+1},\vec{b}^{{{r}}+1})$ is not an optimal allocation within $S_{{{r}}+1}(\vec{d},\vec{b})$; let  $(\vec{d}^\star,\vec{b}^\star)\in S_{{{r}}+1}(\vec{d},\vec{b})$ be an allocation} that minimizes the dock-move distance to $(\vec{d}^{{r}},\vec{b}^{{r}})$ among allocations in $S_{{{r}}+1}(\vec{d},\vec{b})$ with smaller objective than $(\vec{d}^{{{r}}+1},\vec{b}^{{{r}}+1})$. By Lemma \ref{lemma: iexists} there exist $j$ and $k$ such that $d^\star_j+b^\star_j>d^{{{r}}}_j+b^{{{r}}}_j \geq d_j+b_j$ and $d^\star_k+b^\star_k<d^{{{r}}}_k+b^{{{r}}}_k \leq d_k+b_k$. Applying Lemma \ref{lemma: exchange} to $(\vec{d}^{{{r}}},\vec{b}^{{{r}}})$ and $(\vec{d}^\star,\vec{b}^\star)$, we find that there exist $(\vec{d}',\vec{b}')\in N(\vec{d}^{{{r}}},\vec{b}^{{{r}}})$ and $(\vec{d}^{\star\star},\vec{b}^{\star\star})\in N(\vec{d}^{\star},\vec{b}^{\star})${, both not necessarily in $S_{{r}}(\vec{d},\vec{b})$, such that 
$$
c(\vec{d}^{{r}},\vec{b}^{{r}})-c(\vec{d'},\vec{b'})
{\stackrel{\star}{\geq}}
c(\vec{d}^{\star\star},\vec{b}^{\star\star})-c(\vec{d}^{\star},\vec{b}^{\star}), \text{ and either}
$$
\begin{enumerate}
    \item going from $(\vec{d}^{{{r}}},\vec{b}^{{{r}}})$ to $(\vec{d}',\vec{b}')$ involves moving a dock from some $i$ to $j$ and going from $(\vec{d}^{\star},\vec{b}^{\star})$ to $(\vec{d}^{\star\star},\vec{b}^{\star\star})$ involves moving a dock from $j$ to that same $i$, or
    \item going from $(\vec{d}^{{{r}}},\vec{b}^{{{r}}})$ to $(\vec{d}',\vec{b}')$ involves moving a dock from $k$ to some $i$ and going from $(\vec{d}^{\star},\vec{b}^{\star})$ to $(\vec{d}^{\star\star},\vec{b}^{\star\star})$ involves moving a dock from that same $i$ to $k$.
\end{enumerate}
For simplicity we assume the first case; the proof is symmetric for the second. 
Notice that $(\vec{d}^{\star\star},\vec{b}^{\star\star})$ is at closer dock-move distance to $(\vec{d}^{{{r}}},\vec{b}^{{{r}}})$  than $(\vec{d}^\star,\vec{b}^\star)$ is. Further, since $d_j^\star+b_j^\star>d_j+b_j$ and the first case involves moving a dock from $j$ to go from $(\vec{d}^{\star},\vec{b}^{\star})$ to $(\vec{d}^{\star\star},\vec{b}^{\star\star})$, we know that this move cannot increase the dock-move distance to $(\vec{d},\vec{b})$ no matter what other station the dock is going to (an increase in distance at the station the dock is moved to would be canceled out by the decrease in distance at $j$). Thus, it either decreases the distance and gives $(\vec{d}^{\star\star},\vec{b}^{\star\star})\in S_{r}(\vec{d},\vec{b})$ or it keeps the distance constant and gives $(\vec{d}^{\star\star},\vec{b}^{\star\star})\in  S_{r+1}(\vec{d},\vec{b})\setminus S_{r}(\vec{d},\vec{b})$. We derive a contradiction from both.

{Suppose first $(\vec{d}^{\star\star},\vec{b}^{\star\star})\in S_{{r}}(\vec{d},\vec{b})$, then the assumption that $(\vec{d}^{{r}},\vec{b}^{{r}})$ is optimal in $S_{{r}}(\vec{d},\vec{b})$ guarantees that $c(\vec{d}^{{r}},\vec{b}^{{r}})\leq c(\vec{d}^{\star\star},\vec{b}^{\star\star})$. Thus, the RHS of the starred inequality, can be bounded below by $c(\vec{d}^{{{r}}},\vec{b}^{{{r}}})-c(\vec{d}^{\star},\vec{b}^{\star})$. For the LHS consider that the algorithm chose to move to  $(\vec{d}^{{{r}}+1},\vec{b}^{{{r}}+1})$ when it could have chosen $(\vec{d'},\vec{b'})$; thus, we must have $c(\vec{d}^{{{r}}+1},\vec{b}^{{{r}}+1})\leq c(\vec{d'},\vec{b'})$. This allows us to bound the LHS of the starred inequality from above by $c(\vec{d}^{{r}},\vec{b}^{{r}})-c(\vec{d}^{{{r}}+1},\vec{b}^{{{r}}+1})$.  {Combining the two bounds} we find $c(\vec{d}^{{{r}}+1},\vec{b}^{{{r}}+1})\leq c(\vec{d}^{\star},\vec{b}^{\star})$, a contradiction. Thus, we cannot have} $(\vec{d}^{\star\star},\vec{b}^{\star\star})\in S_{{r}}(\vec{d},\vec{b})$.

Now, suppose instead that $(\vec{d}^{\star\star},\vec{b}^{\star\star})\in S_{{{r}}+1}(\vec{d},\vec{b})\setminus S_{{{r}}}(\vec{d},\vec{b})$. We first prove that in that case we must have $(\vec{d'},\vec{b'})\in S_{{r}}(\vec{d},\vec{b})$. Observe that at station $j$ we have $d_j'+b_j' = d_j^r+b_j^r+1>d_j+b_j$; thus, it suffices to show that station $i$ fulfills $d_i^r+b_i^r>d_i+b_i$. In order for $(\vec{d'},\vec{b'})$ to be closer to $(\vec{d}^\star,\vec{b}^\star)$ than~$(\vec{d}^{r},\vec{b}^{r})$ is, we must have 
$d_i^{\star\star}+b_i^{\star\star}=d_i^{\star}+b_i^{\star}+1\leq d_i^r+b_i^r$. However, since $(\vec{d}^{\star\star},\vec{b}^{\star\star})$ is at dock-move distance~$r+1$ from $(\vec{d},\vec{b})$ and since~$d_j^{\star}+b_j^{\star}>d_j+b_j$, it must be the case that $d_i^{\star\star}+b_i^{\star\star}>d_i+b_i$ -- else the move from $j$ to $i$, in going from $(\vec{d}^{\star},\vec{b}^{\star})$ to $(\vec{d}^{\star\star},\vec{b}^{\star\star})$, reduces the distance to $(\vec{d},\vec{b})$. But then we can combine the two aforementioned inequalities to get $d_i^r+b_i^r\geq d_i^{\star\star}+b_i^{\star\star}>d_i+b_i$.
Thus, we have $(\vec{d'},\vec{b'})\in S_{{r}}(\vec{d},\vec{b})$, and may use that $(\vec{d}^{{r}},\vec{b}^{{r}})$ is optimal within $S_r(\vec{d},\vec{b})$ to bound the LHS of the starred inequality from above by 0. However, since $(\vec{d}^{\star\star},\vec{b}^{\star\star})$ is closer to $(\vec{d}^{{r}},\vec{b}^{{r}})$ than $(\vec{d}^{\star},\vec{b}^{\star})$ is, we know, by definition of $(\vec{d}^{\star},\vec{b}^{\star})$, that the RHS must be strictly greater 0, which also yields a contradiction. Thus, we cannot have $(\vec{d}^{\star\star},\vec{b}^{\star\star})\in S_{{{r}}+1}(\vec{d},\vec{b})\setminus S_{{{r}}}(\vec{d},\vec{b})$  which concludes the proof of the theorem. 
}
\hfill\halmos





\section{Scaling techniques}
\label{sec:appendix_proof_scaling}

In this section we first prove that the algorithm described in Section \ref{sec: scaling} (see Algorithm \ref{alg:scaling}) runs in polynomial time for Problem \ref{opt:no_ops}, then describe how to adapt it for Problem \ref{opt:adapted} (see Algorithm \ref{alg:scaling_ops}), and finally prove that this adaptation finds an optimal solution to \ref{opt:adapted} in polynomial time too. \\ 
Before we begin, we introduce some notation. We rely, throughout the section, on the inductive assumption that at the end of phase $k$ we have an allocation $(\vec{d}^k,\vec{b}^k)$ that is optimal among solutions at dock-move distance at most $z$ from $(\vec{\bar{d}},\vec{\bar{b}})$ that differ for each $d_i$ and each $b_i$ by a multiple of $\alpha_k=2^{\lfloor\log_2(B+D)\rfloor+1-k}$. Specifically, $(\vec{d}^k,\vec{b}^k)$ denotes an optimal solution to the optimization problem in phase $k$ which we denote P2$(\alpha_k)$, where $P2(\alpha)$ is defined as:

\OneAndAHalfSpacedXI

\begin{equation*}
\begin{aligned}
\mathtt{minimize}_{(\vec{d},\vec{b})\in\mathbb{Z}_0^2} & c(\vec{d},\vec{b}) &\\
\mathtt{s.t.} &  |\vec{d}+\vec{b}|_1 &= D+B,\\
& |\vec{b}|_1 &= B, \\
& \forall i: b_i &\equiv \bar{b}_i  \mod \alpha,\\
& \forall i: d_i+b_i &\equiv \bar{d}_i+\bar{b}_i   \mod \alpha,\\
& {|\vec{\bar{d}}+\vec{\bar{b}}-\vec{d}-\vec{b}|_1\leq 2z} &\, \text{ where }  |\bar{d}+\bar{b}|_1=D+B,
\end{aligned}
\tag{P2($\alpha$)}\label{opt:scaling}
\end{equation*}
\DoubleSpacedXI

Notice the deviation in notation from Section \ref{sec: algorithm} where $(\vec{d}^r,\vec{b}^r)$ denotes the solution found by Algorithm \ref{alg:gradient_descent} in iteration $r$.
Further, we denote by $(\vec{d}^\star,\vec{b}^\star)$ an optimal solution of P2$(\alpha_{k+1})$ that is at minimum dock-move distance to $(\vec{d}^k,\vec{b}^k)$, {and again use $S_r(\vec{d},\vec{b})=\left\{(\vec{d}',\vec{b}'):|\vec{d}'+\vec{b}'-\vec{d}-\vec{b}|_1\leq 2r\right\}$ to denote the set of allocations at dock-move distance at most $r$ to $(\vec{d},\vec{b})$}. 
\\
\noindent
\textbf{Proof of Theorem \ref{thm:scaling}.}
We use the following bound from \cite{shioura2018m} in our proof, and highlight that \cite{shioura2018m} offers no such bound when $z<\infty$ (the scaling algorithms in that paper do not rely on such a bound, whereas our Algorithm \ref{alg:scaling_ops} does).
\begin{lemma}[Theorem 6.5 in \cite{shioura2018m}]\label{lem:shioura_bound}
With $z=\infty$, and $(\vec{d}^k,\vec{b}^k)$, $(\vec{d}^\star,\vec{b}^\star)$ defined as above, we have $|\vec{d^{k}}+\vec{b^{k}})-(\vec{d}^\star+\vec{b}^\star)|_1< 8n\alpha_k$. 
\end{lemma}
In the last phase of Algorithm \ref{alg:scaling}, when $\alpha_k=1$, the algorithm terminates when no gradient-descent step yields improvement; Lemma \ref{lemma: neighborhoods} implies that this only occurs when a globally optimal solution has been found. Thus, it follows that Algorithm \ref{alg:scaling} finds an optimal solution to \ref{opt:no_ops}. What remains to be shown is that the number of iterations required is bounded by $O(n\log(B+D))$. Since there are $O(\log(B+D))$ phases, it suffices to show that each phase consists of at most $O(n)$ iterations. In the first phase at most one iteration can be carried out (as $\alpha_0\geq B+D$ at most one station can have docks removed from it at most once); for each subsequent phase Theorem \ref{thm: koptimal} implies that the number of iterations in phase $k$ is equal to the minimum dock-move distance in phase $k$ (with each move picking $\alpha_{k}$ docks/bikes) between the solution found in phase $k-1$ and a solution that is optimal for phase $k$. Lemma \ref{lem:shioura_bound} bounds this dock-move distance by~$O(n)$.\hfill\halmos
\\\noindent\textbf{Description of Algorithm \ref{alg:scaling_ops}.}
We adapt Algorithm \ref{alg:scaling} in the following way to accommodate the operational constraints (see Algorithm \ref{alg:scaling_ops}): rather than starting phase $k+1$ at the optimal dock allocation found in phase $k$, we move (see Algorithm \ref{alg:find_d_b}) to a new starting point $(\vec{d},\vec{b})$  for phase~$k+1$. Our goal in changing the starting point to $(\vec{d},\vec{b})$ is to find an allocation with the property that for every~$i$ we have either $d_i^\star+b_i^\star\leq d_i+b_i\leq \bar{d}_i+\bar{b}_i$ or $d_i^\star+b_i^\star\geq d_i+b_i\geq \bar{d}_i+\bar{b}_i$, i.e., we want to \emph{sandwich} the number of docks at each station $i$ between $d_i^\star+b_i^\star$ and $\bar{d}_i+\bar{b}_i$. While finding such $(\vec{d},\vec{b})$ may be nontrivial, as we do not know $(\vec{d}^\star,\vec{b}^\star)$, we first focus on the advantage of initiating the gradient-descent algorithm in phase $k+1$ at $(\vec{d},\vec{b})$: 
As long as we have~{$|\vec{d}+\vec{b}|_1=D+B=|\vec{d}^\star+\vec{b}^\star|_1$}, Theorem \ref{thm: koptimal} guarantees that within {$|(\vec{d}+\vec{b})-(\vec{d}^\star+\vec{b}^\star)|_1/(2\alpha_{k+1})$} iterations a solution with objective no worse than $c(\vec{d}^\star,\vec{b}^\star)$ is found. Further, the solution found is at dock-move distance
{$$
|(\vec{d}+\vec{b})-(\vec{d}^\star+\vec{b}^\star)|_1 + |(\vec{d}+\vec{b})-(\vec{\bar{d}}+\vec{\bar{b}})|_1=|(\vec{d}^\star+\vec{b}^\star)-(\vec{\bar{d}}+\vec{\bar{b}})|_1\leq 2z
$$
}
from $(\vec{\bar{d}},\vec{\bar{b}})$, i.e., it fulfills the operational constraints. This implies that, given such $(\vec{d},\vec{b})$ Algorithm~\ref{alg:scaling_ops} finds an optimal solution for phase $k+1$ in {$|(\vec{d}+\vec{b})-(\vec{d}^\star+\vec{b}^\star)|_1/(2\alpha_{k+1})$} iterations. 
Proving that the algorithm runs in polynomial time then requires us to show that, using $(\vec{d}^k,\vec{b}^k)$, (i) we can find such an allocation $(\vec{d},\vec{b})$ while (ii) guaranteeing that  {$|(\vec{d}+\vec{b})-(\vec{d}^\star+\vec{b}^\star)|_1/(2\alpha_{k+1})$} is bounded by a polynomial in $n$.
%
\\\noindent\textbf{Proof of Theorem \ref{thm:scaling_ops}.}
%
Similar to Lemma \ref{lem:shioura_bound} we rely on the following bound.
\begin{lemma}\label{lem:scaling_bound}
For any $z$ we have {$|(\vec{d}^k+\vec{b}^k)-(\vec{d}^\star+\vec{b}^\star)|_1\leq 32n^6(n+4)\alpha_{k+1}$}.
\end{lemma}
We emphasize that the bounds in this section are only meant to prove that Algorithm~\ref{alg:scaling_ops} runs in polynomial time, not to be tight. Thus, we emphasize, in various places, simplicity over  tighter bounds. Now, denoting by $M$ an upper bound on {$|(\vec{d}^k+\vec{b}^k)-(\vec{d}^\star+\vec{b}^\star)|_1$}, we characterize in the next lemma the output of the allocation $(\vec{d},\vec{b})$ returned by Algorithm \ref{alg:find_d_b}. Based on Lemma \ref{lem:scaling_bound} we know that $M\leq 32n^6(n+4)\alpha_{k+1}$. 
\begin{lemma}\label{lem:d_b_found}
Algorithm \ref{alg:find_d_b} returns $(\vec{d},\vec{b})$ such that $| \vec{d}+\vec{b}|_1=D+B$ and for every~$i$
\begin{itemize}
    \item if $d_i^k+b_i^k\geq \bar{d}_i+\bar{b}_i$, we have 
$$
\max\{d_i^k+b_i^k-nM,\bar{d}_i+\bar{b}_i\}\leq
d_i+b_i\leq \max\{d_i^k+b_i^k-M,\bar{d}_i+\bar{b}_i\},$$ 
\item if $d_i^k+b_i^k\leq \bar{d}_i+\bar{b}_i$
we have $$
\min\{d_i^k+b_i^k+M,\bar{d}_i+\bar{b}_i\}\leq
d_i+b_i\leq \min\{d_i^k+b_i^k+nM,\bar{d}_i+\bar{b}_i\}.$$
\end{itemize}
\end{lemma}
We prove both lemmas further below after first applying them to prove Theorem \ref{thm:scaling_ops}.
\\
\emph{Proof of Theorem.} We argue that any $d_i+b_i$ fulfilling the two  inequalities in Lemma \ref{lem:d_b_found} must be sandwiched between $d_i^\star+b_i^\star$ and~$\bar{d}_i+\bar{b}_i$. Suppose we have  $d_i^k+b_i^k\geq \bar{b}_i+\bar{d}_i$, the proof is symmetric for the other direction. Observe that Lemma \ref{lem:scaling_bound} implies that $d_i^k+b_i^k-M\leq d_i^\star+b_i^\star\leq d_i^k+b_i^k+M$. But then we either have $\bar{d}_i+\bar{b}_i= d_i+b_i$, or $\bar{d}_i+\bar{b}_i\leq d_i+b_i\leq d_i^k+b_i^k-M\leq d_i^\star+b_i^\star$, so $d_i+b_i$ is sandwiched as required. At the same time, we can use the inequalities to bound 
$$|(d_i+b_i)-(d^\star_i+b^\star_i)|
\leq
|(d_i+b_i)-(d^k_i+b^k_i)|+
|(d^\star_i+b^\star_i)-(d^k_i+b^k_i)|\leq Mn+M.
$$
Thus, we have also argued that $|(\vec{d}+\vec{b})-(\vec{d}^\star+\vec{b}^\star)|_1/\alpha_k$ is polynomially bounded by $n$.
\\
We have shown that in each phase $O(nM)$ iterations of Algorithm \ref{alg:find_d_b} yield an allocation $(\vec{d},\vec{b})$ from which~$O(Mn^2)$ gradient-descent steps suffice to find the optimal solution for the phase, so Algorithm~\ref{alg:scaling_ops} finds an optimal solution in time polynomial in $n$ and $\log(D+B)$.\hfill\Halmos 

{
\textbf{Proof of Lemma \ref{lem:scaling_bound}. } 
{Our proof of this lemma proceeds in three steps. We first argue that for the purpose of bounding the dock-move distance we may argue about allocations that have the same allocation of docks as $(\vec{d}^k,\vec{b}^k)$ and $(\vec{d}^\star,\vec{b}^\star)$, but different allocations of bikes. Then we consider a particular multiset of dock-moves to get from one allocation to the other. In  Claim~1  we characterize a condition that must be true if a single dock-move appears more than~$n+4$ times in this multiset. Finally, in Claim~2, we show that if a single dock-move appears ``too often'' (more than~$4n^3(n+4)$ times), then that condition is violated. Thus, each individual dock-move cannot occur too often, and since we have $2n(n-1)(n-2)+2n(n-1)\leq 4n^3$ distinct dock-moves (see Definition \ref{def:neighborhood}), we can bound the number of dock-moves between $(\vec{d}^k,\vec{b}^k)$ and~$(\vec{d}^\star,\vec{b}^\star)$ by $4n^3\times 4n^3(n+4)$.}

For notational simplicity we prove the result for the penultimate and the ultimate phase, i.e., when~$\alpha_k=2$ and~$\alpha_{k+1}=1$. The proof for any other phases $k'$ and $k'+1$ is the same with each distance multiplied by $\alpha_{k'+1}$ and each dock-move replaced by $\alpha_{k'+1}$ copies of itself (see definition of $N^\alpha(\cdot,\cdot)$ in Appendix \ref{appendix:pseudocode}).
Recall that we consider a solution $(\vec{d}^k,\vec{b}^k)$ that is optimal in phase $k$ and a solution~$(\vec{d}^\star,\vec{b}^\star)$ that minimizes the dock-move distance to $(\vec{d}^k,\vec{b}^k)$ among optimal solutions for phase $k+1$, i.e., P2$(\alpha_{k+1})$. Though~$(\vec{d}^k,\vec{b}^k)$ is optimal in phase $k$, it is not necessarily bike-optimal in phase $k+1$ (where the number of bikes at each station need not be even). Denote by~$(\vec{d}',\vec{b}')$ an allocation of bikes in phase~$k+1$ that is closest to the dock allocation of~$(\vec{d}^k,\vec{b}^k)$, i.e., it allocates~$d_i^k+b_i^k$ docks to each station $i$, it is bike-optimal for phase $k+1$, and it minimizes~$|\vec{b}^k-\vec{b}'|_1$. Observe that $|\vec{b}^k-\vec{b}'|_1\leq 2n$ (this is immediate by Fact 2 since otherwise there exist stations $i$ and $j$ such that moving 2 bikes from $i$ to $j$ in $(\vec{d}^k,\vec{b}^k)$ yields improvement, which violates optimality of $(\vec{d}^k,\vec{b}^k)$ in phase $k$), and therefore for every~$i$ we have $|b_i^k-b_i'|\leq n$. Also, since $(\vec{d}^k,\vec{b}^k)$ and $(\vec{d}',\vec{b}')$ differ in their allocation of bikes, but not in their allocation of docks, we know
$|(\vec{d^k}+\vec{b^k})-(\vec{d^{k+1}}+\vec{b^{k+1}})|_1=|(\vec{d'}+\vec{b'})-(\vec{d^{k+1}}+\vec{b^{k+1}})|_1$
.
\\
We consider a set of minimal size of dock-moves that lead from $(\vec{d}',\vec{b}')$ to a bike-optimal allocation with $d_i^\star+b_i^\star$ docks at each station $i$. By Lemma \ref{lemma: bike_opt} we can get to a bike-optimal allocation of $d_i^\star+b_i^\star$ docks at each station $i$ by maintaining bike-optimality in each of the $|(\vec{d}^\star+\vec{b}^\star)-(\vec{d}'+\vec{b}')|_1/2$ moves. We relabel $(\vec{d}^\star,\vec{b}^\star)$ to be that allocation.\footnote{Since there may be more than one bike-optimal allocation for the same allocation of docks, this guarantees that there exist $|(\vec{d}^\star+\vec{b}^\star)-(\vec{d}'+\vec{b}')|_1/2$ dock-moves that lead from $(\vec{d}',\vec{b}')$ to $(\vec{d}^\star,\vec{b}^\star)$.}
Now, consider a multiset $L$ \dfedit{of minimum size} that minimizes the number of moves of type either $E_{ijh}$ or~$O_{ijh}$ such that applying each move in $L$ to $(\vec{d}',\vec{b}')$ yields~$(\vec{d}^\star,\vec{b}^\star)$. {Then, by bounding $|L|\leq 16n^6(n+4)$ we simultaneously obtain our desired bound on the dock-move distance between~$(\vec{d^k},\vec{b^k})$ and~$(\vec{d^\star},\vec{b^\star})$.}\\
We argue about a transformation $X$, either a dock-move or the composition of two different dock-moves, that occurs repeatedly in~$L$. When referring to a composition of two different dock-moves occurring repeatedly, we mean that each of them occurs repeatedly.
Consider a transformation~$X$  that occurs more than $n+4$ times in~$L${, e.g., a move of an empty dock from station $i$ to station~$j$ or a move of an empty dock from station $i$ to station~$j$ combined with a move of a full dock from station~$i'$ to station~$j'$. 
We first prove the following condition for any transformation $X$ that appears at least $(n+4)$ times in $L$.}\\
\textbf{Claim 1. }
{For any transformation $X$ we use the notation ~$X^2(\vec{d}^k,\vec{b}^k):=X(X(\vec{d}^k,\vec{b}^k))$ and denote by   $(\vec{d}^{\star\star},\vec{b}^{\star\star})$ the allocation that fulfills~$X(\vec{d}^{\star\star},\vec{b}^{\star\star})=(\vec{d}^{\star},\vec{b}^{\star})$. Then, with} $X$ appearing $\geq n+4$ times in $L$ we have $X^2(\vec{d}^k,\vec{b}^k)\not\in S_z(\vec{\bar{d}},\vec{\bar{b}})$ or $(\vec{d}^{\star\star},\vec{b}^{\star\star})\not\in S_z(\vec{\bar{d}},\vec{\bar{b}})$.\\
{Next, in Claim 2 we show that if, within~$L$, a single dock-move from some station $i$ to some station~$j$ appears at least~$4n^3(n+4)$ times, then there must be a transformation~$Y$ that appears at least~$(n+4)$ times in~$L$ but violates the condition of Claim 1.}
}
\\
{
\textbf{Claim 2.}
\dfedit{If $L$ contains any dock-move from $i$ to $j$, i.e., one of $e_{ij},o_{ij},E_{ijh},$ or $O_{ijh}$ for any $h$, at least $4n^3(n+4)$ times, then $L$ must contain $n+4$ copies of a transformation $Y$ such that, with~$Y(\vec{d}^{\star\star},\vec{b}^{\star\star})=(\vec{d}^{\star},\vec{b}^{\star})$, we have $Y^2(\vec{d}^k,\vec{b}^k)\in S_z(\vec{\bar{d}},\vec{\bar{b}})$ and $(\vec{d}^{\star\star},\vec{b}^{\star\star})\in S_z(\vec{\bar{d}},\vec{\bar{b}})$}.\\
{This leads to a contradiction if a single dock-move appears at least~$4n^3(n+4)$ times. Since, as mentioned above, we have less than $4n^3$ distinct dock-moves, it follows that, with each one appearing at most $4n^3(n+4)$ times, we have $|L|\leq 4n^3\times 4n^3(n+4)=16n^6(n+4)$. \halmos}
}
\\
\emph{Proof of Claim 1.} 
Suppose $X$ appears $\geq n+4$ times in $L$ and we have both $X^2(\vec{d}^k,\vec{b}^k)\in S_z(\vec{\bar{d}},\vec{\bar{b}})$ and~$(\vec{d}^{\star\star},\vec{b}^{\star\star})\in S_z(\vec{\bar{d}},\vec{\bar{b}})$. From this we derive that $X^2(\vec{d}^k,\vec{b}^k)$ is feasible in the penultimate phase and that $(\vec{d}^{\star\star},\vec{b}^{\star\star})$ is feasible in the last phase of the scaling algorithm --- these both follow since dock-moves keep the total number of allocated bikes and docks constant, though the first also requires that $X^2$, by definition, moves docks and bikes in pairs. But then the definitions of $(\vec{d}^k,\vec{b}^k)$ (optimal in phase $k$) and $(\vec{d}^\star,\vec{b}^\star)$ (at minimum distance from $(\vec{d}^k,\vec{b}^k)$ among optimal allocations in phase $k+1$) imply that
$$
c(\vec{d}^k,\vec{b}^k)\leq c(X^2(\vec{d}^k,\vec{b}^k)) \text{ and } c(\vec{d}^\star,\vec{b}^\star)< c(\vec{d}^{\star\star},\vec{b}^{\star\star}).
$$
Thus, if we can argue (similar to Lemma \ref{lemma: exchange}) \dfedit{that we have
$$
0\geq c(\vec{d}^k,\vec{b}^k)- c(X^2(\vec{d}^k,\vec{b}^k)) {\geq} c(\vec{d}^{\star\star},\vec{b}^{\star\star})-c(\vec{d}^\star,\vec{b}^\star)>0,
$$
then we derive a contradiction. Notice that all stations not involved in $X$ cancel out in the differences.\\
We derive this inequality by using the multimodular inequalities to bound the change at $(\vec{d}^k,\vec{b}^k)$ and at~$(\vec{d}^{\star\star},\vec{b}^{\star\star})$ for each station $j$ involved in $X$; in order to apply the inequalities, we first need to understand the relative values of $(d^k_j,b^k_j)$ and at~$(d_j^{\star},b_j^{\star})$.}

\begin{figure}[ht]
    \centering
    \includegraphics[width=.44\textwidth]{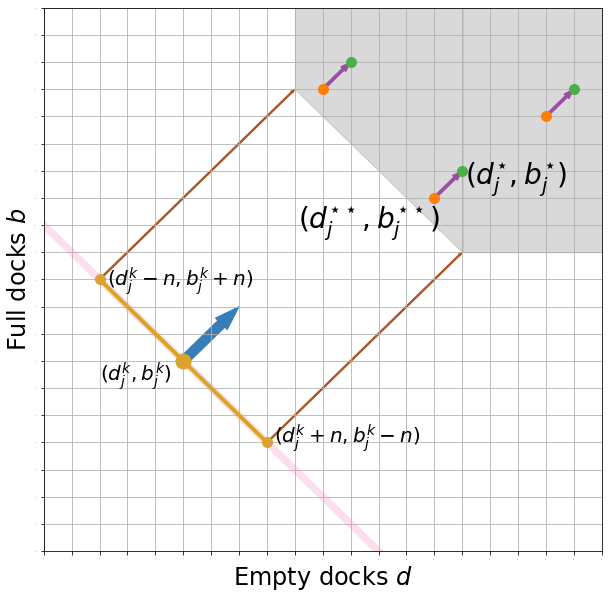}
    \includegraphics[width=.45\textwidth]{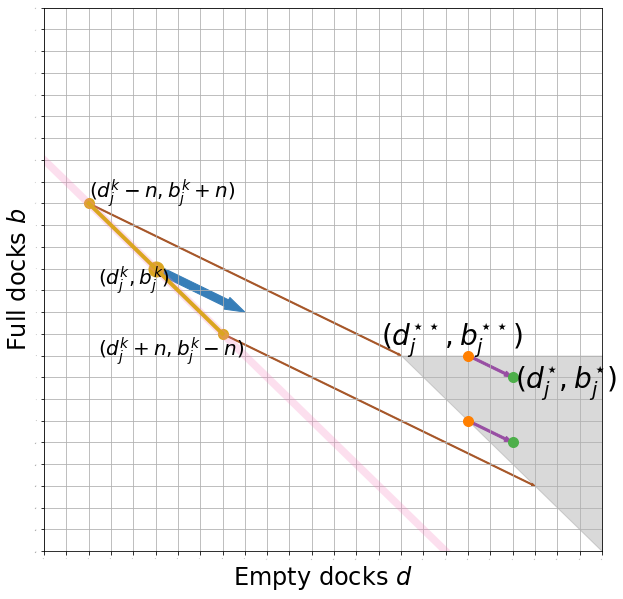}
    \caption{Bounding the improvement of two added docks, one empty and one full, (LHS) and an added empty dock and a bike removal (RHS) at $(d_j^{\star\star},b_j^{\star\star})$ relative to that at $(d_j^k,b_j^k)$.}
    \label{fig:lemma_scaling_extra}
\end{figure}

\dfedit{Suppose $L$ involves a station $j$ receiving an empty dock. In this case $L$ cannot involve a move in which $j$ has either an empty or a full dock taken away: if it did, then $L$ would contain both dock-moves from some $r$ to $j$ and from $j$ to some $s$; two such moves could be replaced by one dock-move from $r$ to $s$ which would yield a multiset $L$ of smaller size such that applying each move in this multiset to $(\vec{d}',\vec{b}')$ yields the same allocation of docks in phase $k+1$ as $(\vec{d}^\star,\vec{b}^\star)$, i.e., it would contradict $L$ being a smallest multiset of dock-moves to get to this allocation (recall from Lemma~\ref{lemma: bike_opt} that bike-optimality can be maintained in all dock-moves so reaching the same allocation of docks in fewer dock-moves suffices for a contradiction). 
Next, $L$ may not include any moves of the kind~$O_{rsj}$ for some~$r$ and~$s$, as we must have either~$O_{ijh}$ or~$o_{ij}$ in~$L$ for~$j$ to receive an empty dock, and as we can replace an occurrence~(i) of~$O_{ijh}$ and~$O_{rsj}$ by~$e_{rj}$ and~$e_{is}$ and~(ii) of $o_{ij}$ and~$O_{rsj}$ by~$o_{is}$ and~$e_{rj}$ to obtain a multiset $L'$ with strictly fewer moves of type $E_{ijh}$ or $O_{ijh}$ than $L$ has. Thus, we have excluded all dock-moves from $L$ that remove docks (empty or full) from~$j$ or add bikes to $j$. It follows that if $L$ includes a dock-move in which $j$ receives an empty dock, then all dock-moves in $L$ that include~$j$ must (i) add empty docks to~$j$, (ii) add full docks to~$j$, or (iii) take away bikes from~$j$ (recall from Definition \ref{def:neighborhood} the types of moves that would yield each of these).
With the same reasoning we can partition the ways in which dock-moves in $L$ affect $j$ into the following six mutually exclusive cases: dock-moves in $L$ affect $j$ through
\begin{enumerate}
    \item added empty and full docks;
    \item added empty docks and removed bikes;
    \item added full docks and added bikes;
    \item removed empty and full docks;
    \item removed empty docks and bikes added;
    \item removed full docks and bikes.
\end{enumerate}
We focus on the first two cases as the others follow through symmetries either along the axis of removing/adding or the axis of empty/full. In particular, the first case is symmetric to the fourth by exchanging added and removed, and the same holds for second/fifth and third/sixth. Further, the second case is symmetric to the third by exchanging the coordinates of empty and full docks. Since the multimodular inequalities are symmetric with respect to the coordinates, the first and second case thus imply the others.\\
Since we focus on the first two cases, notice that in these cases $X$ must have one of the following effects at $j$: (a) one (or two) empty docks added, (b) one (or two) full docks added, (c) one empty and one full dock added, (d) one empty dock added and one bike removed, or (e) one (or two) bikes removed. As (a) and (b) are symmetric with empty/full exchanged, we omit (b). 
For~(a),~(c), and~(d) Figure~\ref{fig:lemma_scaling} displays the relative position of $(d_j^{\star},b_j^{\star})$ and $(d_j^k,b_j^k)$: $(d_j^{\star},b_j^{\star})$ has to be in the shaded region given that dock-moves in $L$ involve at least $n+4$ empty docks being added to $j$, that combining all the moves in $L$ leads from $(\vec{d}',\vec{b}')$ to $(d_j^{\star},b_j^{\star})$, that $d_j'+b_j'=d_j^k+b_j^k$ but the number of bikes can be off by up to $n$ (meaning that $(\vec{d}',\vec{b}')$ is along the yellow line), and that dock-moves in~$L$ do not involve docks being removed from $j$ or bikes (without docks) being added to $j$  --- below we will hone in on the special cases of (c) and (d).
\begin{figure}[ht]
    \centering
    \includegraphics[width=.45\textwidth]{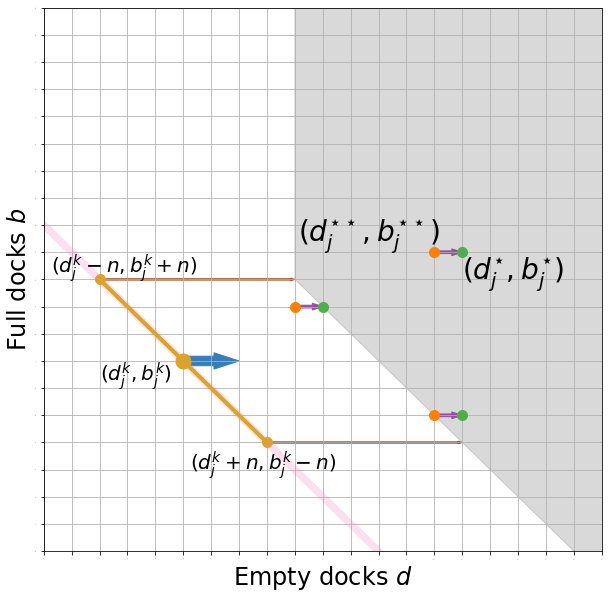}
    \caption{Bounding the improvement of an added empty dock at $(d_j^{\star\star},b_j^{\star\star})$ relative to that at $(d_j^k,b_j^k)$: the yellow line indicates the possible position of $(d_j',b_j')$, the grey shaded region the possible position of $(d_j^{\star\star},b_j^{\star\star})$ given that it involves $n+4$ added empty docks starting at $(d_j',b_j')$.}
    \label{fig:lemma_scaling}
\end{figure}\\
Under~(a) we know that $X$ adds either one or two empty docks to~$j$, and Inequalities (3) and (6) (recall Figures~\ref{fig:inequalities} and~\ref{fig:lemma_3}) guarantee that
$$
\frac{c_j(d_j^k,b_j^k)-c_j(d_j^k+4,b_j^k)}{4}
\geq c_j(d_j^\star-\frac{4}{2},b_j^\star)-c_j(d_j^\star,b_j^\star)
$$ 
for the case of $X$ adding two empty docks to $j$; for the case of $X$ adding one empty dock to $j$, the same inequality holds with all $4$s replaced by $2$s. \\
In case (c) $X$ instead involves an empty and a full dock being added to $j$: then $X$ appearing $n+4$ times in $L$ implies that the feasible region for $(d_j^{\star},b_j^{\star})$, relative to $(d_j^{k},b_j^{k})$ is as displayed in the shaded region in the left plot of Figure \ref{fig:lemma_scaling_extra}, i.e., there are at least $n+4$ more empty and at least~$n+4$ more full docks than at $(d_j',b_j')$. In this case, the improvement of adding two empty and two full docks at $(d_j^{k},b_j^{k})$ is greater equal to twice that of adding one empty and one full dock at~$(d_j^{\star\star},b_j^{\star\star})$, as is clear from Inequalities~(4) and~(6) for the empty docks and~(1) and~(5) for the full docks. Next, in case~(d),~$X$ involves one empty dock being added to~$j$ and one bike being removed. Then, the feasible region would be as shaded in the right plot of Figure \ref{fig:lemma_scaling_extra}, and in this case the improvement of two bikes being removed and two empty docks being added at~$(d_j^k,b_j^k)$ is greater equal to twice that of one empty dock added and one bike removed at $(d_j^{\star\star},b_j^{\star\star})$. Combining Inequalities (3) and (4) we find that the improvement is decreasing in the number of empty docks (horizontal axis), and combining Inequalities (2) and (3) we find that the improvement is increasing in the number of bikes (along the diagonal axis, i.e., it is decreasing as we move to the bottom-right).
Finally, in case (e), Figure \ref{fig:lemma_scaling_bike} displays the relative position of $(d_j^{\star},b_j^{\star})$ and $(d_j^k,b_j^k)$; as the cost of removing a bike is increasing as one goes downwards along the pink diagonal (this is Fact~1 in Appendix~\ref{appendix:bike_opt}) and increasing as one moves to the right -- this is Inequality (3) -- it is smaller equal at~$(d_j^k,b_j^k)$ and at ~$X(d_j^k,b_j^k)$ than it is at~$(d_j^{\star\star},b_j^{\star\star})$.}

\dfedit{
Thus, the improvement of $X^2$ is greater-equal at $(\vec{d}^k,\vec{b}^k)$ than that of~$X$ at~$(\vec{d}^{\star\star},\vec{b}^{\star\star})$, i.e., 
$$
c(\vec{d}^k,\vec{b}^k)- c(X^2(\vec{d}^k,\vec{b}^k)) {\geq} c(\vec{d}^{\star\star},\vec{b}^{\star\star})-c(\vec{d}^\star,\vec{b}^\star),
$$
which is a contradiction to $X^2(\vec{d}^k,\vec{b}^k)\in S_z(\vec{\bar{d}},\vec{\bar{b}})$ and~$(\vec{d}^{\star\star},\vec{b}^{\star\star})\in S_z(\vec{\bar{d}},\vec{\bar{b}})$.
\halmos
} 
\begin{figure}[ht]
    \centering
    \includegraphics[width=.4\textwidth]{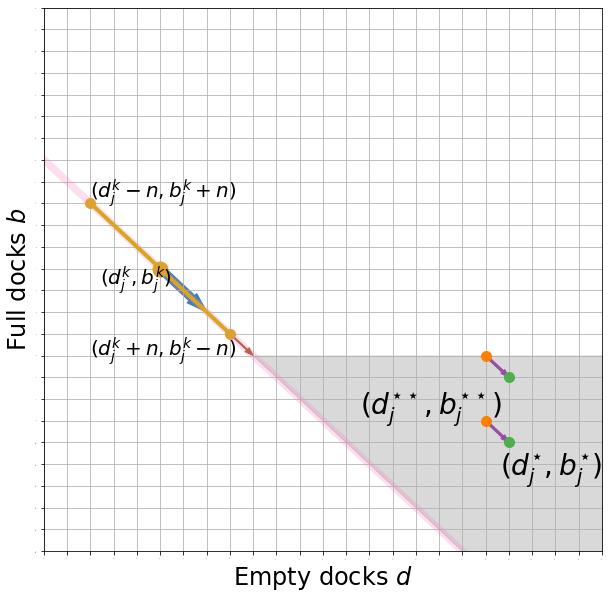}
    \caption{Bounding the cost of a removed bike at $(d_j^{\star\star},b_j^{\star\star})$ relative to that at $(d_j^k,b_j^k)$ when $L$ has $n+4$ dock-moves that remove a bike from $j$.}
    \label{fig:lemma_scaling_bike}
\end{figure}

\emph{Proof of Claim 2.}
In order to prove Claim 2 we start by partitioning stations with $d_r^k+b_r^k\neq d_r^\star+b_r^\star$ into six different sets based on the ordering of the magnitude of $d_r^k+b_r^k$, $d_r^\star+b_r^\star$, and $\bar{d}_r+\bar{b}_r$. Specifically, we define
\begin{itemize}
    \item $S_1=\{r: d_r^k+b_r^k>d_r^\star+b_r^\star\geq  \bar{d}_r+\bar{b}_r\}$
    \item $S_2=\{r:d_r^k+b_r^k>\bar{d}_r+\bar{b}_r>d_r^\star+b_r^\star  \}$
    \item $S_3=\{r:   \bar{d}_r+\bar{b}_r\geq d_r^k+b_r^k>d_r^\star+b_r^\star\}$
    \item $S_4=\{r: d_r^\star+b_r^\star>d_r^k+b_r^k\geq  \bar{d}_r+\bar{b}_r\}$
    \item $S_5=\{r:d_r^\star+b_r^\star>\bar{d}_r+\bar{b}_r>d_r^k+b_r^k  \}$
    \item $S_6=\{r:   \bar{d}_r+\bar{b}_r\geq d_r^\star+b_r^\star>d_r^k+b_r^k\}$
\end{itemize}
Notice first that $|d_r^k+b_r^k-\bar{d}_r-\bar{b}_r|\equiv 0 \mod \alpha_k$ holds for every $r$. It follows that if, e.g., $d_r^k+b_r^k>\bar{d}_r+\bar{b}_r$, then also~$d_r^k+b_r^k-\alpha_k=d_r^k+b_r^k-2\alpha_{k+1}\geq\bar{d}_r+\bar{b}_r$. Let the transformation $X$ denote the dock-move from $i$ to $j$ that occurs at least $4n^3(n+4)$ times. Observe that we must have $i\in S_1\cup S_2\cup S_3$ and $j\in S_4\cup S_5\cup S_6$.

Suppose~$i\in S_1$ and~$j\in S_4$. Then between $(\vec{d}^k,\vec{b}^k)$ and $(\vec{d}^\star,\vec{b}^\star)$ the dock-move distance to~$(\vec{\bar{d}},\vec{\bar{b}})$ is invariant to that move, i.e., when $d_i^k+b_i^k\geq d_i+b_i>d_i^\star+b_i^\star$ and $d_j^k+b_j^k\leq d_j+b_j<d_i^\star+b_i^\star$, then moving a dock from $i$ to $j$ in $(\vec{d},\vec{b})$ keeps the dock-move distance to $(\vec{\bar{d}},\vec{\bar{b}})$ fixed. In particular, since~$(\vec{d}^k,\vec{b}^k)$ and~$(\vec{d}^\star,\vec{b}^\star)$ are both feasible with respect to~$S_z(\vec{\bar{d}},\vec{\bar{b}})$, we then also know that~$X^2(\vec{d}^k,\vec{b}^k)$ and~$(\vec{d}^{\star\star},\vec{b}^{\star\star})$, defined as before,  are feasible with respect to $S_z(\vec{\bar{d}},\vec{\bar{b}})$, giving a contradiction to Claim~1. That same argument holds true when $i\in S_3$ and $j\in S_6$. Furthermore, with $i\in S_2$, regardless of $j\in S_4,S_5,$ or $S_6$ we must have both $X^2(\vec{d}^k,\vec{b}^k)$ and $(\vec{d}^{\star\star},\vec{b}^{\star\star})$ in $S_z(\vec{\bar{d}},\vec{\bar{b}})$, and the same holds true for $j\in S_5$ regardless of $i$. Thus, we cannot have $n+4$ copies of dock-moves from $i$ to $j$ with any of these combinations.

The remaining combinations are $i\in S_1, j\in S_6$ and $i\in S_3, j\in S_4$. We focus on $i\in S_1, j\in S_6$, the argument is symmetric for $i\in S_3, j\in S_4$. 
Notice that with $i\in S_1, j\in S_6$ the fact that $(\vec{d}^k,\vec{b}^k)$ is feasible guarantees that $X^2(\vec{d}^k,\vec{b}^k)\in S_z(\vec{\bar{d}},\vec{\bar{b}})$. However, $(\vec{d}^{\star\star},\vec{b}^{\star\star})$ is at greater dock-move distance to $(\vec{\bar{d}},\vec{\bar{b}})$ than $(\vec{d}^{\star},\vec{b}^{\star})$. Now, if $(\vec{d}^{\star},\vec{b}^{\star})$ was at dock-move distance at most $z-\alpha_{k+1}$ from $(\vec{\bar{d}},\vec{\bar{b}})$, then~$(\vec{d}^{\star\star},\vec{b}^{\star\star})$ would be guaranteed to be at distance at most $z$, and $X$ would violate the conditions in Claim~1. Otherwise, the moves in $L$, in aggregate, cannot be decreasing the dock-move distance to~$(\vec{\bar{d}},\vec{\bar{b}})$; thus, we must have at least as many moves in $L$ that increase the distance as we have moves that decrease the distance (notice that a dock-move from a station in $S_2$ to one in $S_6$ decreases the distance at $(\vec{d}^k,\vec{b}^k)$, but keeps it constant at $(\vec{d}^{\star\star},\vec{b}^{\star\star})$). However, each move from $i$ to~$j$ decreases the distance  to $(\vec{\bar{d}},\vec{\bar{b}})$ by $\alpha_{k+1}$. Since we know that there are at least $4n^3(n+4)$ moves from $i$ to $j$ in~$L$, and since there are less than $4n^3$ moves in total, there must be a dock-move, from some $h$ to some~$\ell$, in $L$ that 
\dfedit{appears repeatedly and increases the dock-move distance to~$(\vec{\bar{d}},\vec{\bar{b}})$, in aggregate across all its appearances,  by at least $(n+4)$}. Thus, we cannot have $h\in S_1,\ell\in S_6$. If~$h$ and $\ell$ are of the combinations excluded in the last paragraph, then we are done, as the moves from $h$ to $\ell$ \dfedit{would then} violate the assumptions of Claim~1. Else, we must have~$h\in S_3,\ell\in S_{\dfedit{4}}$. But then, let $Y$ be the composition of $X$ and that dock-move from $h$ to $\ell$. Now, $Y$ appears at least~$n+4$ times in $L$ and keeps the dock-move distance to $(\vec{\bar{d}},\vec{\bar{b}})$ constant, i.e., it is a transformation that violates the assumptions of Claim~1. {This completes the proof of Claim~2.} 
 \hfill\halmos

\textbf{Proof of Lemma \ref{lem:d_b_found}. }
Notice first that in each iteration Algorithm \ref{alg:find_d_b} carries out one dock-move (of $\alpha_k$ bikes and docks), so the number of bikes and docks remains constant throughout. Further, in each iteration at least one station gets at least $\alpha_k$ closer to being within the required bounds; thus, the algorithm terminates in at most $nM/\alpha_k$ iterations.
Now, define sets $\mathcal{S}^+=\{i:d_i^k+b_i^k\geq\bar{d}_i+\bar{b}_i\}$ and~$\mathcal{S}^-=\{j:d_j^k+b_j^k<\bar{d}_j+\bar{b}_j\}$. 
We want to argue that in each iteration in which the current solution $(\vec{d},\vec{b})$ does not fulfill the conditions of the lemma, the algorithm can update the solution. If the conditions are not fulfilled, one of the following two must exist: 
\begin{enumerate}
    \item a station $i\in \mathcal{S}^+$ with $d_i+b_i>\max\{d_i^k+b_i^k-M,\bar{d}_i+\bar{b}_i\}$ 
    \item a station $j\in \mathcal{S}^-$ with $d_j+b_j<\min\{d_j^k+b_j^k+M,\bar{d}_j+\bar{b}_j\}$.
\end{enumerate}
If both exist, then the algorithm carries out a dock-move from one such $i$ to one such $j$, and the next iteration begins. Suppose the first holds true but the second does not (the proof for the second existing but not the first is symmetric); denote a station~$i\in \mathcal{S}^+$ with $d_i+b_i>\max\{d_i^k+b_i^k-M,\bar{d}_i+\bar{b}_i\}$ as $i'$. Then, either there exists $j\in \mathcal{S}^-$ such that $d_j+b_j<\max\{d_j^k+b_j^k+nM,\bar{d}_j+\bar{b}_j\}$ or no such $j$ exists, i.e., for every $j\in \mathcal{S}^-$ we have $d_j+b_j\geq\max\{d_j^k+b_j^k+nM,\bar{d}_j+\bar{b}_j\}$. In the first case, again, a possible iteration involves a dock-move from $i'$ to $j$, so the algorithm does not terminate. In the second case notice we have either (i) $d_j+b_j=\bar{d}_j+\bar{b}_j$ for every $j\in \mathcal{S}^-$ or (ii) there exists $j'\in \mathcal{S}^-$ such that $d_{j'}+b_{j'}=d_{j'}^k+b_{j'}^k+nM$. If (i) holds we derive from $D+B=\sum_i d_i+b_i=\sum_i\bar{d}_i+\bar{b}_i$ that $\sum_{i\in \mathcal{S}^+}(d_i+b_i)-(\bar{d}_i+\bar{b}_i)=\sum_{j\in \mathcal{S}^-}(\bar{d}_j+\bar{b}_j)-(d_j+b_j)$, where each summand in both sums is non-negative. Thus, we also have $d_i+b_i\leq\bar{d}_i+\bar{b}_i$ for every~$i\in \mathcal{S}^+$, a contradiction to there being an~$i\in \mathcal{S}^+$ with $d_i+b_i>\max\{d_i^k+b_i^k-M,\bar{d}_i+\bar{b}_i\}$. If~(ii) holds, then the Algorithm must have carried out at least $nM/\alpha_k$ iterations already, which is a contradiction to the conditions of the lemma not yet being fulfilled. \hfill\halmos 

\section{Running Time}\label{sec: running_time}

\begin{figure}[h]
\centering
\begin{subfigure}{.5\textwidth}
  \centering
  \includegraphics[width=.75\textwidth]{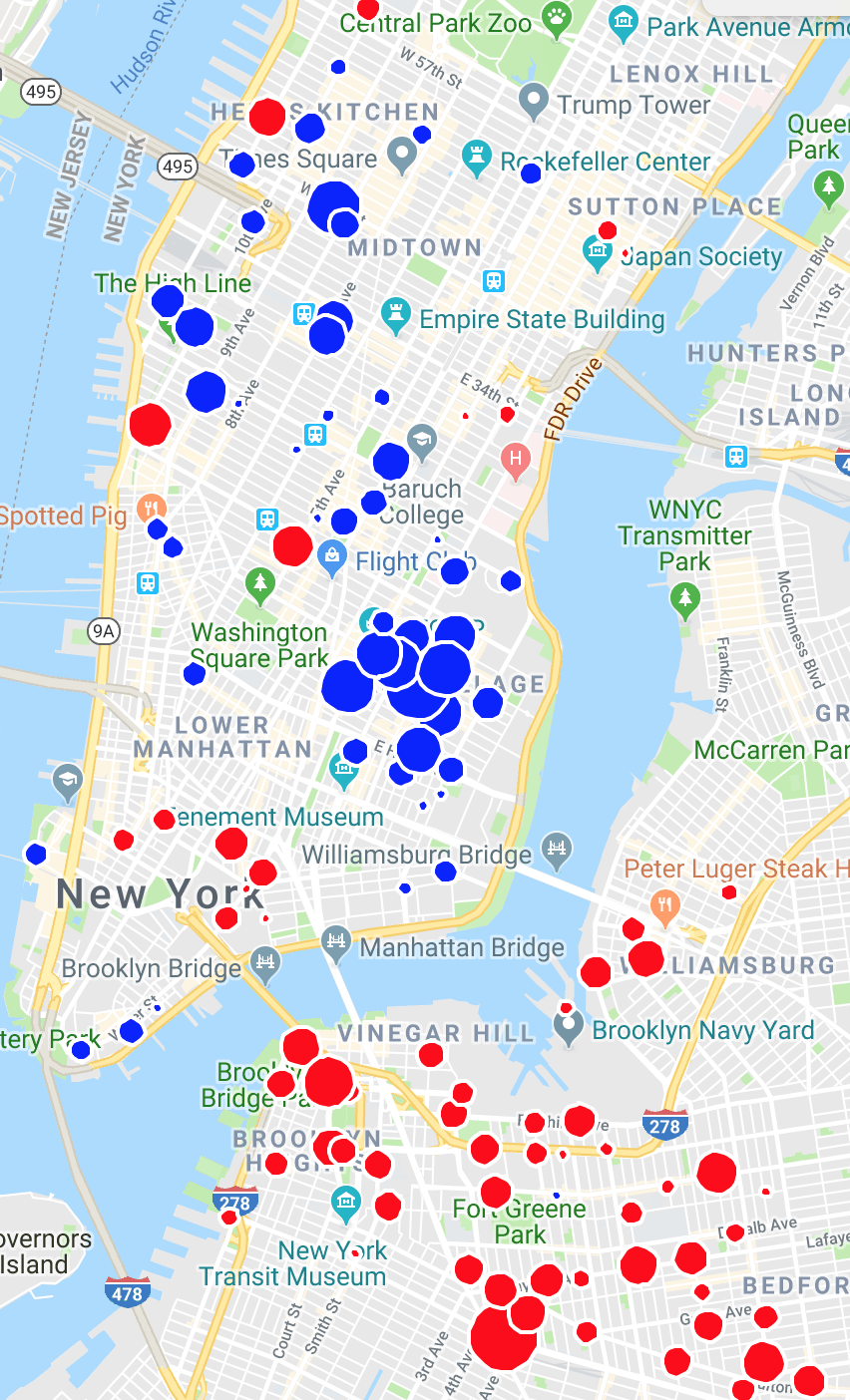}
\end{subfigure}%
\begin{subfigure}{.5\textwidth}
  \centering
  \includegraphics[width=.75\textwidth]{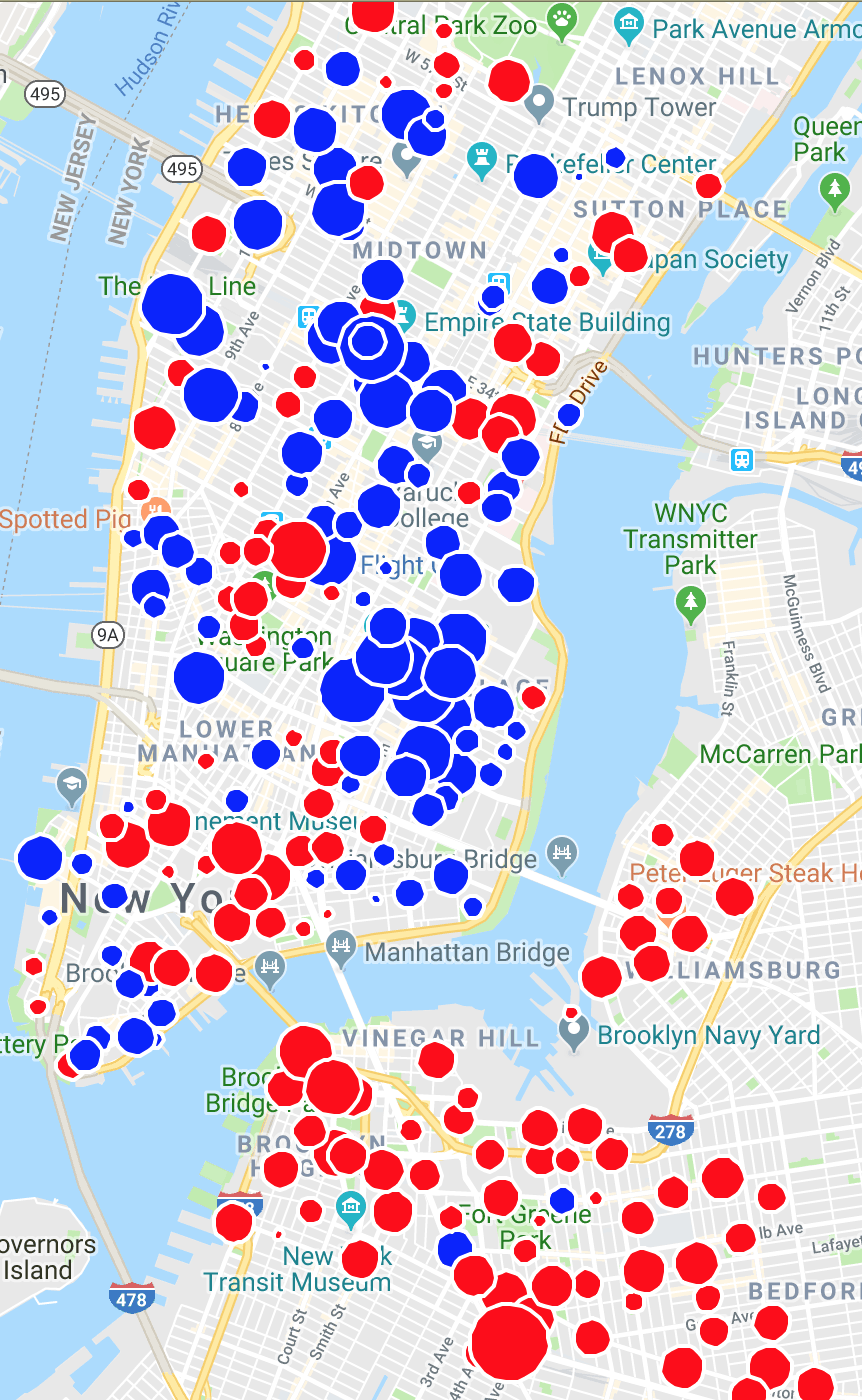}
\end{subfigure}
\caption{Visualization of docks moved by optimal solution in NYC for $z\in\{500,1500\}$; red circles correspond to stations at which docks are taken, blue circles to stations at which they are added.}
\label{fig: map_solution}
\end{figure}

Even though the reallocation of docks is a strategic question, the time to solve
the associated optimization models is not irrelevant for practical considerations. Given the expensive computation of each user dissatisfaction value, an early approach to compute the LP-relaxation of the optimization problem took a weekend to solve (on a high-end laptop). This was due to the time required to set up the LP; once it was set up, it solved quickly. While this is acceptable for a one-off analysis, in practice system operators care about regularly running different analyses that include different demand patterns, different bounds on number of docks moved, and even different bounds on station sizes. Having a fast algorithm allows system operators to run the analysis without our support. We provided them with a Jupyter notebook (\citealt{Kluyver2016aa}) that includes the entire workflow, from estimating the demand profiles to computing the user dissatisfaction functions to running the optimization problem to creating map-based visualizations of the resulting solutions (see Figure~\ref{fig: map_solution}) and does not rely on specialized optimization software like Gurobi or CPLEX. Crucially, this workflow happens, on a MacBook Pro with a 2.8 GHz dual core processor and 8GB of RAM, in a matter of minutes rather than hours or days (see Table \ref{table:running_time}).

\OneAndAHalfSpacedXI
\DoubleSpacedXI

\begin{figure}[ht]
    \centering
\includegraphics[height=.25  \textwidth]{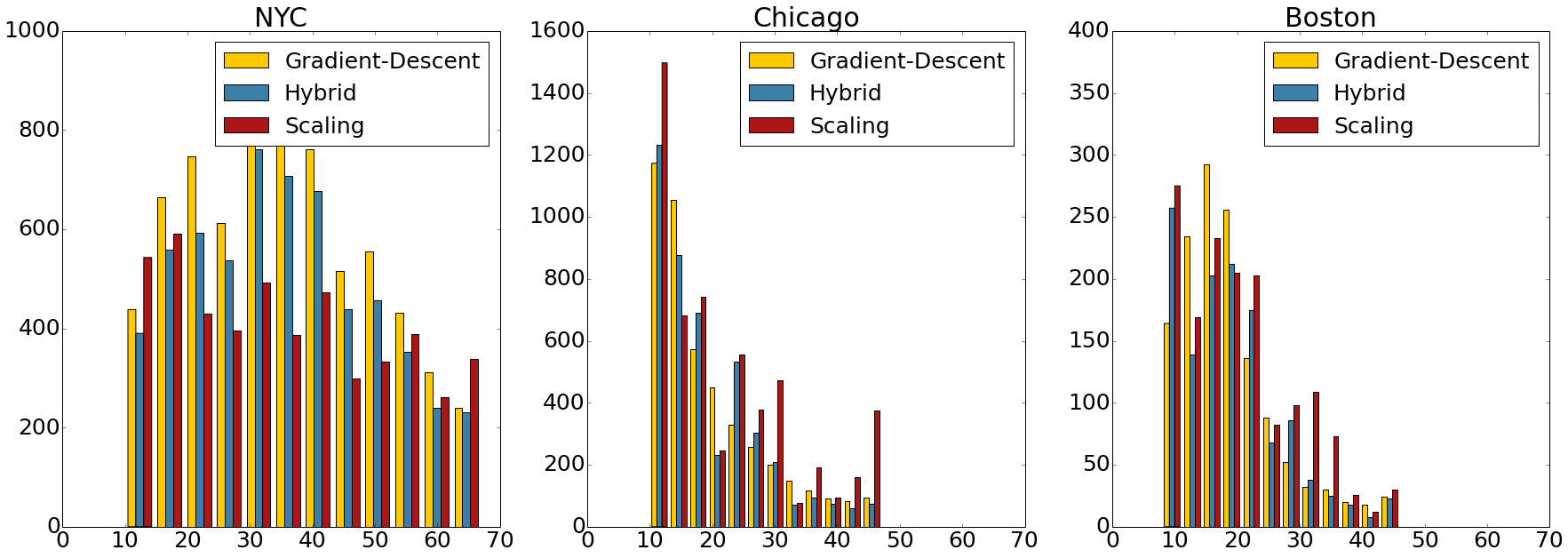}
    \caption{Number of user dissatisfaction functions, grouped by capacity $d+b$, evaluated by each algorithm in each city.}
    \label{fig:CChistogram}
\end{figure}

To conclude this numerical exploration, we now compare the measured running times of the gradient-descent and the scaling algorithm. Given that the running-time of each algorithm is dominated by the computational effort to compute values of the user dissatisfaction functions (the effort for which grows as the cube of the capacity), we only computed values to which the respective algorithm needed access. In  Figure \ref{fig:CChistogram} we plot the number of user dissatisfaction functions that are computed by each algorithm. In Chicago, the scaling algorithm created unnecessary overhead by requiring values for large capacities at many stations that the gradient-descent algorithm did not. This illustrates why the gradient-descent algorithm outperforms the scaling algorithm in both Boston and Chicago (see Table \ref{table:running_time}). In NYC on the other hand, the scaling algorithm performed significantly better than the gradient-descent algorithm. Motivated by this contrast, we implemented a hybrid algorithm (see Algorithm \ref{alg:hybrid}) that only iterates over  3, 2, and 0 as values of $k$, rather than all powers of 2. The hybrid outperforms the gradient-descent algorithm on all three data-sets and outperforms the scaling algorithm on the data-sets from Chicago and Boston. All three algorithms outperform, by 2 orders of magnitude, the linear programming based approach that needs to evaluate every value of the user dissatisfaction functions at all stations before solving.

\begin{table}[ht]
\centering
\begin{tabular}{c|c|c|c|}
\cline{2-4}
                                    & \multicolumn{3}{c|}{Running Time (Minutes)}                                                   \\ \cline{2-4} 
                                    & {Gradient-descent}                        & Hybrid                        & Scaling                       \\ \hline
\multicolumn{1}{|c|}{New York City} & 14.88 & 12.78 & 10.83 \\ \hline
\multicolumn{1}{|c|}{Chicago}       & 5.03  & 4.75  & 6.57  \\ \hline
\multicolumn{1}{|c|}{Boston}        & 1.40  & 1.30  & 1.70  \\ \hline
\end{tabular}
\caption{Comparison of the running times of each of the three algorithms in each of the three cities}
\label{table:running_time}
\end{table}


\OneAndAHalfSpacedXI
\DoubleSpacedXI

\section{Code reproducibility}\label{sec:ec_online}

As part of this paper we publish a Jupyter notebook \dfedit{and} a data set that together allow for the reproducibility of most of our data-driven results. Below we outline the extent to which each of the results in the figures and tables can be reproduced. Since Figures \ref{fig:fill_levels} and \ref{fig: map_solution} are merely included to motivate our results, we restrict our attention to Figures \ref{fig:UDF2D}, \ref{fig:kmoves}, \ref{fig:evaluation}, and \ref{fig:CChistogram}, as well as Tables~\ref{tab:summary}--\ref{table: impact} and \ref{table:running_time}. The extent to which the results can be reproduced is limited by the fact that the data set made public only involves demand estimates from June 2016. As such, results based on other months cannot be reproduced. 

{The data we use for most of our analyses could have been collected from public sources. Specifically, our optimization algorithms are applied to demand data that requires two ingredients: the monthly ridership data that Citi Bike makes public, and (to decensor) timestamps of when stations were empty/full; the latter can be collected from \url{https://gbfs.citibikenyc.com/gbfs/en/station_status.json}. However, since we have not collected this data ourselves (and it cannot be collected retroactively), our data is limited by what Motivate allowed us to publish. That said, by collecting data from that URL at frequent timestamps, e.g., every two minutes, one can create analogous data for future months, and use our scripts to conduct similar analyses.}

\subsubsection*{Figure \ref{fig:UDF2D}. } Both parts of Figure \ref{fig:UDF2D} are exactly reproduced in the notebook (see cell 5 in the notebook).

\subsubsection*{Table \ref{tab:summary}.} The results for each cell in the Table are taken from the notebook (see cells 7, 12, and 17 in the notebook).

\subsubsection*{Figure \ref{fig:kmoves}.} The lines in Figure \ref{fig:kmoves} that correspond to June 2016 are reproduced in the notebook (see cells 10 and 19); since the other lines rely on data from July and August 2016, they are not reproducible using the demand estimates we are making public.

\subsubsection*{Table \ref{table: seasons}.} The three numbers reported in Table \ref{table: seasons} are not reproducible. For March and November 2017, this is due to demand estimates not being shared for these months. But even when restricting the analysis of moving 200 docks only to June 2016, the value in the table is not exactly reproducible with the shared data set. This is because of how we handled, in our analysis, the variability of data available across different months: in each case, we restricted the set of stations to the stations for which demand estimates were available for all months considered. The data we are making public, to reproduce Table \ref{tab:summary} and Figure \ref{fig:kmoves}, is restricted to the stations for which we can provide demand estimates for each of June, July, and August 2016; however, the set of stations for which demand estimates are available for June 2016, March 2017, and November 2017 is a different one (e.g., due to a station being taken offline for some time to allow for street construction); thus, even the cell in Table \ref{table: seasons} that corresponds to June 2016 cannot be reproduced.

\subsubsection*{Figure \ref{fig:evaluation} and Table \ref{table: impact}.}  
Table \ref{table: impact} and Figure \ref{fig:evaluation} cannot be reproduced using the data we are making public. However, we created a synthetic data set, that is being shared in addition to the demand estimates, and demonstrate in the supplement how to use it to create a Figure that closely resembles Figure~\ref{fig:evaluation} (see cell 21).

\subsubsection*{Figure \ref{fig:CChistogram} and Table \ref{table:running_time}.}
Figure \ref{fig:CChistogram} can be reproduced exactly in the notebook (see cell~35). The numbers in Table \ref{table:running_time} are taken from a run of the notebook (see cells 22--30) but the running time can vary between different computational setups.



\section{Connections to discrete convex functions} 
\subsection{$M$-Convex Functions}\label{appendix_M}

In this appendix we first provide the definitions of $M$-convex sets and functions, and then show that our objective with budget constraints is not $M$-convex. For the definitions, it is useful to denote $supp^+(\vec{x}-\vec{y})=\{i:x_i>y_i\}$, $supp^-(\vec{x}-\vec{y})=\{i:x_i<y_i\}$, and $\vec{e}_i$ as the canonical unit vector.

\begin{definition}[$M$-convex set]
A nonempty set of integer points $B\subseteq Z^{2n}$ is defined to be an $M$-convex set if it satisfies $
\forall \vec{x},\vec{y}\in B, i\in supp^+(\vec{x}-\vec{y}), \exists j\in supp^-(\vec{x}-\vec{y}): \vec{x}-\vec{e}_i+\vec{e}_j\in B$.
\end{definition}

\begin{definition}[$M$-convex function]
A function $f$ is $M$-convex if for all $x,y\in dom(f), i\in supp^+(x-y), \exists j\in supp^-(x-y): f(x)+f(y)\geq f(x-e_i+e_j)+f(y+e_i-e_j)$.
\end{definition}

\cite{kaspi2015bike} prove a statement equivalent to $c(\cdot,\cdot)$ being $M$-convex. \cite{murota2004steepest} characterized the minimum of an $M$ convex function as follows to show that Algorithm \ref{algo2} minimizes $M$-convex functions:

\begin{lemma}[\cite{murota1996convexity,murota1998discrete}]\label{lemma:Mconvex}
For an $M$-convex function $f$ and $x\in dom(f)$ we have $f(x)\leq f(y)\;\forall y\in dom(f) \; \textit{ if and only if }\; f(x)\leq f(x-e_i+e_j)\forall i,j$.
\end{lemma}

\OneAndAHalfSpacedXI
\begin{algorithm}
\caption{$M$-convex function minimization, cf. \cite{murota2003discrete}}
\label{algo2}
\begin{algorithmic}[1]
\item Find a vector $x\in dom(f)$
\State Find $i,j$ that minimize $f(x-e_i+e_j)$
\State If $f(x)>f(x-e_i+e_j)$, set $x:=x-e_i+e_j$ and go to 2 
\State Else, return $x$
\end{algorithmic}
\end{algorithm}
\vspace{-.3in}

\DoubleSpacedXI
As the following example shows, the restriction of $c$ to the feasible set (with budget constraints) does not guarantee $M$-convexity, despite both the set and $c$ being $M$-convex.
\begin{example}
Our example consists of three stations $i,j$, and $k$ with demand-profiles:
\begin{eqnarray*}
p_i(-1)=\frac{1}{2},\; p_i(+1,-1)=\frac{1}{2}; &
p_j(+1)=\frac{1}{2}; &
p_k(+1,-1,-1)=1.
\end{eqnarray*}
\end{example}
We consider two solutions. In the first, $i$, $j$, and $k$ each have a dock allocated with $i$ also having a bike allocated, i.e., $b_i'=d_j'=d_k'=1$, whereas $d_i'=b_j'=b_k'=0$ and our budget constraint is $D=2,B=1$. Then $c_i(d_i',b_i')=\frac{1}{2}$, $c_j(d_j',b_j')=0$, and $c_k(d_k',b_k')=1$. In the second solution, $d_i^*=b_k^*=d_k^*=1$, whereas $b_i^*=d_j^*=b_j^*=0$. Thus, we have  $c_i(d_i^*,b_i^*)=\frac{1}{2}$, $c_j(d_j^*,b_j^*)=\frac{1}{2}$, and $c_k(d_k^*,b_k^*)=0$, giving that $1=c(\vec{d}^*,\vec{b}^*)<c(\vec{d}',\vec{b}')=\frac{3}{2}$.
But then the statement of Lemma \ref{lemma:Mconvex} with $y=(\vec{d}^*,\vec{b}^*)$ and $x=(\vec{d},\vec{b})$ implies that, if $c$ is $M$-convex then {either moving an empty dock from $j$ or $k$ to $i$, or moving a full dock from $i$ to $j$ or $k$ must yield a solution better than $c(\vec{d}',\vec{b}')$. This is not the case, and therefore $c$, restricted to the feasible set is not $M$-convex, even though the underlying feasible set is $M$-convex.}

\subsection{Connections to Discrete Midpoint Convex Functions}\label{appendix_dmcf}
In this appendix we show that the constrained optimization problem formulated in Section \ref{sec: model} is not multimodular. To do so, we apply an equivalence proven in \cite{murota2005note} that characterizes a function $f$ as multimodular if and only if there exists an $L^\natural$ convex function $g$ such that $f(x_1,x_2,\ldots,x_n)=g(x_1, x_1+x_2,\ldots,\sum_{i=1}^n x_i)$. While we do not state the explicit definition of $L^\natural$ convex functions here, it was shown by \cite{fujishige2000notes} that $L^\natural$ convex functions fulfill the following discrete midpoint convexity property.

\begin{definition}
A function $g:\mathbb{Z}^n\to\mathbb{R}^n\cup\{+\infty\}$ is called discrete midpoint convex if
$$
g(x)+g(y)\geq g\big(\lceil\frac{x+y}{2}\rceil\big)+g\big(\lfloor\frac{x+y}{2}\rfloor\big).
$$
Here, the floor and ceiling refer to component-wise floor and ceiling.
\end{definition}

We now argue that the function $g$ corresponding to our (constrained) objective $c$ is not discrete midpoint convex. Consider the current allocation (see Section \ref{sec: model}) $\vec{\bar{d}}=(0,1,0,1)$ and $\vec{\bar{b}}=(0,0,0,0)$. As all values for $\vec{b}$ are 0 throughout this construction, we do not restate it from now on. Suppose $z=1$, that is, only one dock is allowed to be moved (and solutions moving more than one are infeasible and thus have value infinity). Then the vector  $\vec{\bar{d}}=(1,0,1,0)$ is not feasible given the constraint (as it would involve moving 2 docks). Now, if $g$ was discrete midpoint convex, then the inequality would state that
\begin{eqnarray*}
f(1,0,0,1)+f(0,1,1,0)=g(1,1,1,2)+g(0,1,2,2)\\\geq g(1,1,2,2)+g(0,1,1,2)=f(1,0,1,0)+f(0,1,0,1).
\end{eqnarray*}
However, both terms on the LHS are feasible whereas the first term on the RHS is not. Thus, the inequality does not hold, $g$ is not discrete midpoint convex, and therefore $f$ is not multimodular.

\end{document}